\documentclass[a4paper,12pt,amsart,frenchb]{article}

\usepackage{amsmath,amsbsy,amsfonts,amssymb}
\usepackage[french]{babel}
\newcommand{\ass}[2]{\vskip0.3cm\noindent
{\bf {#1}}. { \sl {#2}}\vskip0.3cm\noindent
}
  
\begin{document}

 \title{  La conjecture locale de Gross-Prasad pour les repr\'esentations temp\'er\'ees des groupes sp\'eciaux orthogonaux}
\author{J.-L. Waldspurger}
\date{19 novembre 2009}
\maketitle

  {\bf Introduction}

\bigskip

Soit $F$ un corps local non archim\'edien de caract\'eristique nulle. Appelons espace quadratique un couple $(V,q)$ o\`u $V$ est un espace vectoriel de dimension finie sur $F$ et $q$ est une forme bilin\'eaire sym\'etrique et non d\'eg\'en\'er\'ee sur $V$. Fixons deux tels espaces quadratiques $(V,q)$ et $(V',q')$. On note $d$ et $d'$ leurs dimensions et $G$ et $G'$ leurs groupes sp\'eciaux orthogonaux. On suppose $d$ pair et $d'$ impair.  On suppose donn\'e un isomorphisme entre le plus grand des espaces et la somme orthogonale du plus petit et d'un espace quadratique  qui est lui-m\^eme somme orthogonale de plans hyperboliques et d'une droite quadratique $(D,q_{D})$. Le plus petit des deux groupes $G$, $G'$ est alors un sous-groupe du plus grand. On   d\'efinit un \'el\'ement $\nu_{0}\in F^{\times}/F^{\times,2}$ en fixant  un \'el\'ement non nul $v_{0}\in D$  et en posant
$$\nu_{0}=\left\lbrace\begin{array}{cc}q_{D}(v_{0},v_{0})/2,&\text{ si }d>d',\\ -q_{D}(v_{0},v_{0})/2,&\text{ si }d<d'.\\ \end{array}\right.$$  

Soit $\sigma$, resp. $\sigma'$, une repr\'esentation admissible et irr\'eductible de $G(F)$, resp. $G'(F)$. Gross et Prasad ont d\'efini une multiplicit\'e $m(\sigma,\sigma')$ en [GP2]. Rappelons la d\'efinition dans  deux cas extr\^emes. Supposons $d'=d+1$. Fixons des espaces $E_{\sigma}$ et $E_{\sigma'}$ dans lesquels se r\'ealisent les repr\'esentations $\sigma$ et $\sigma'$. On note  $Hom_{G}(\sigma',\sigma)$ l'espace des applications lin\'eaires $l:E_{\sigma'}\to E_{\sigma}$ telles que $l\circ \sigma'(g)=\sigma(g)\circ l$ pour tout $g\in G(F)$. La multiplicit\'e $m(\sigma,\sigma')$ est la dimension de l'espace complexe $Hom_{G}(\sigma',\sigma)$. Supposons maintenant $d=0$. Le groupe $G$ est \'egal \`a $\{1\}$ et la repr\'esentation $\sigma$ dispara\^{\i}t. Le groupe $G'$ est d\'eploy\'e. Fixons un sous-groupe unipotent maximal $U'$ de $G'$ et un caract\`ere r\'egulier $\psi_{U'}$ de $U'(F)$. On note $Hom_{\psi_{U'}}( \sigma',{\mathbb C})$ l'espace des formes lin\'eaires $l$ sur $E_{\sigma'}$ telles que $l\circ \sigma'(u')=\psi_{U'}(u')l$ pour tout $u'\in U'(F)$. La multiplicit\'e $m(\sigma')$ est la dimension de l'espace $Hom_{\psi_{U'}}(\sigma',{\mathbb C})$. La d\'efinition g\'en\'erale est rappel\'ee en 2.1 ci-dessous. D'apr\`es [AGRS] et [GGP] corollaire 15.2, on a toujours $m(\sigma,\sigma')\leq1$.

On se limitera d\'esormais aux repr\'esentations temp\'er\'ees. Rappelons la classification conjecturale de ces repr\'esentations (conjecture de Langlands, pr\'ecis\'ee par Deligne et Lusztig). Notons $W_{F}$ le groupe de Weil de $\bar{F}/F$, o\`u $\bar{F}$ est une cl\^oture alg\'ebrique de $F$ et $W_{DF}=W_{F}\times SL(2,{\mathbb C})$ le groupe de Weil-Deligne.  Notons $Sp(d'-1,{\mathbb C})$ le groupe symplectique  d'un espace symplectique complexe  de dimension $d'-1$. Notons $\Phi_{temp}(G')$ l'ensemble des classes de conjugaison par $Sp(d'-1,{\mathbb C})$ d'homomorphismes $\varphi:W_{DF}\to Sp(d'-1,{\mathbb C})$ qui v\'erifient quelques propri\'et\'es usuelles et qui sont temp\'er\'es, c'est-\`a-dire que leurs restrictions \`a $W_{F}$ sont d'images relativement compactes. Consid\'erons un tel $\varphi$. Poussons-le en un homomorphisme \`a valeurs dans $GL(d'-1,{\mathbb C})$. On peut alors le d\'ecomposer sous la forme
$$(1) \qquad \varphi=\oplus_{i\in I}l_{i}\varphi_{i}$$
o\`u chaque $\varphi_{i}$ est un homomorphisme irr\'eductible de $W_{DF}$ dans un groupe $GL(N(\varphi_{i}),{\mathbb C})$ et $l_{i}$ est sa multiplicit\'e. On note $I^{symp}$ le sous-ensemble des $i\in I$ tels que $N(\varphi_{i})$ est pair et, \`a conjugaison pr\`es, l'image de $\varphi_{i}$ est contenue dans $Sp(N(\varphi_{i}),{\mathbb C})$.  Notons $S_{\varphi}$ le commutant de l'image de $\varphi$ dans $Sp(d'-1,{\mathbb C})$ et $S_{\varphi}^0$ sa composante neutre.Le groupe $S_{\varphi}/S_{\varphi}^0$ est isomorphe \`a $({\mathbb Z}/2{\mathbb Z})^{I^{symp}}$. Notons $z_{\varphi}$ l'image dans ce groupe de l'\'el\'ement central $-1\in Sp(d'-1,{\mathbb C})$. Il s'identifie \`a $(l_{i})_{i\in I^{symp}}\in ({\mathbb Z}/2{\mathbb Z})^{I^{symp}}$. Posons
$$\mu(G')=\left\lbrace\begin{array}{cc}1,&\text{ si }G'\text{ est d\'eploy\'e,}\\ -1,&\text{ sinon,}\\ \end{array}\right.$$
et notons ${\cal E}^{G'}(\varphi)$ l'ensemble des caract\`eres $\epsilon$ du groupe $S_{\varphi}/S_{\varphi}^0$ tels que $\epsilon(z_{\varphi})=\mu(G')$. On conjecture que l'ensemble des classes d'isomorphie de repr\'esentations admissibles irr\'eductibles et temp\'er\'ees de $G'(F)$ est union disjointe de $L$-paquets $\Pi^{G'}(\varphi)$ index\'es par les $\varphi\in \Phi_{temp}(G')$. Pour un tel $\varphi$, on conjecture que l'ensemble $\Pi^{G'}(\varphi)$ est en bijection avec ${\cal E}^{G'}(\varphi)$, on notera $\epsilon\mapsto \sigma'(\varphi,\epsilon)$ cette bijection. Cette param\'etrisation doit \^etre  compatible avec deux types d'endoscopie. D'une part avec l'endoscopie usuelle entre $G'$ et ses groupes endoscopiques elliptiques. Ceux-ci sont des produits $G'_{+}\times G'_{-}$ de groupes sp\'eciaux orthogonaux d\'eploy\'es d'espaces quadratiques de dimensions $d'_{+}$ et $d'_{-}$ impaires et telles que $d'_{+}+d'_{-}=d'+1$. On renvoie \`a 4.2 pour la description des propri\'et\'es que doivent v\'erifier les param\'etrages relativement \`a ce type d'endoscopie. D'autre part, dans le cas o\`u $G'$ est d\'eploy\'e, on veut une compatibilit\'e avec une endoscopie tordue. Usuellement, on consid\`ere $G'$ comme un groupe endoscopique du groupe $GL(d'-1)$ tordu par un automorphisme ext\'erieur. Mais $G'$ est aussi un groupe endoscopique du groupe $GL(d')$ tordu et c'est ce cas d'endoscopie que nous utiliserons. Plus pr\'ecis\'ement,  notons $\theta$ l'automorphisme usuel $g\mapsto {^tg}^{-1}$ de $GL(d')$. Introduisons le groupe $GL(d')\rtimes \{1,\theta\}$ et sa composante non neutre $\tilde{GL}(d')$.  Notons ${\bf 1}$ la repr\'esentation triviale de dimension $1$ de $W_{DF}$. Soit $\varphi\in \Phi_{temp}(G')$, posons $\varphi_{>}=\varphi\oplus {\bf 1}$. La correspondance de Langlands, prouv\'ee par  Harris-Taylor ([HT])  et Henniart ([H]), associe \`a $\varphi_{>}$ une repr\'esentation admissible irr\'eductible $\pi(\varphi_{>})$ de $GL(d',F)$. Elle est temp\'er\'ee et autoduale, donc se prolonge en une repr\'esentation de $GL(d',F)\rtimes \{1,\theta\}$.  On normalise cette extension de fa\c{c}on ad\'equate et on note $\tilde{\pi}(\varphi_{>})$ sa restriction \`a $\tilde{GL}(d',F)$. On d\'efinit son caract\`ere $\Theta_{\tilde{\pi}(\varphi_{>})}$. De m\^eme, pour toute repr\'esentation irr\'eductible $\sigma'$ de $G'(F)$, on note $\Theta_{\sigma'}$ son caract\`ere. On pose
$$\Theta^{G'}(\varphi)=\sum_{\sigma'\in \Pi^{G'}(\varphi)}\Theta_{\sigma'}.$$
On conjecture alors que $\Theta^{G'}(\varphi)$ est une distribution stablement invariante sur $G'(F)$ et qu'il existe un nombre complexe $c^{G'}(\varphi)$ de module $1$ tel que $c^{G'}(\varphi)\Theta_{\tilde{\pi}(\varphi_{>})}$ soit le transfert endoscopique de $\Theta^{G'}(\varphi)$.  

Rappelons plus succinctement la param\'etrisation conjecturale des repr\'esentations temp\'er\'ees de $G(F)$ (on utilise la formulation de [M]). Fixons un "espace quadratique complexe" de dimension $d$, notons $O(d,{\mathbb C})$ et $SO(d,{\mathbb C})$ son groupe orthogonal, resp. sp\'ecial orthogonal. Consid\'erons un homomorphisme $\varphi:W_{DF}\to O(d,{\mathbb C})$. En composant avec le d\'eterminant, on obtient un caract\`ere quadratique de $W_{DF}$ qui est forc\'ement trivial sur $SL(2,{\mathbb C})$ et se restreint en un caract\`ere de $W_{F}$. Par la th\'eorie du corps de classes, celui-ci d\'etermine un \'el\'ement $\delta(\varphi)\in F^{\times}/F^{\times,2}$. On d\'efinit un autre \'el\'ement de ce groupe par $\delta(q)=(-1)^{d/2}det(q)$. Notons $\Phi_{temp}(G)$ l'ensemble des classes de conjugaison par $SO(d,{\mathbb C})$ d'homomorphismes $\varphi:W_{DF}\to O(d,{\mathbb C})$  qui sont temp\'er\'es et tels que $\delta(\varphi)=\delta(q)$. Consid\'erons un tel $\varphi$. En le poussant en un homomorphisme \`a valeurs dans $GL(d,{\mathbb C})$, on peut le d\'ecomposer sous la forme (1). On note $I^{orth}$ le sous-ensemble des $i\in I$ tels que, \`a conjugaison pr\`es, l'image de $\varphi_{i}$ soit contenue dans $O(N(\varphi_{i}),{\mathbb C})$. Notons $S_{\varphi}$ le commutant de l'image de $\varphi$ dans $SO(d,{\mathbb C})$ et $S_{\varphi}^0$ sa composante neutre.Le groupe $S_{\varphi}/S_{\varphi}^0$ est isomorphe \`a un sous-groupe de $({\mathbb Z}/2{\mathbb Z})^{I^{orth}}$. Notons $z_{\varphi}$ l'image dans ce groupe de l'\'el\'ement central $-1\in SO(d,{\mathbb C})$. Il s'identifie \`a $(l_{i})_{i\in I^{symp}}\in ({\mathbb Z}/2{\mathbb Z})^{I^{symp}}$.  Si $\delta(q)=1$, on pose
$$\mu(G)=\left\lbrace\begin{array}{cc}1,&\text{ si }G\text{ est d\'eploy\'e,}\\ -1,&\text{ sinon.}\\ \end{array}\right.$$
Si $\delta(q)\not=1$, le groupe $G$ est toujours quasi-d\'eploy\'e. Dans ce cas, on fixe $\mu(G)\in \{\pm 1\}$.
Notons ${\cal E}^{G}(\varphi)$ l'ensemble des caract\`eres $\epsilon$ du groupe $S_{\varphi}/S_{\varphi}^0$ tels que $\epsilon(z_{\varphi})=\mu(G)$. On conjecture que l'ensemble des classes d'isomorphie de repr\'esentations admissibles irr\'eductibles et temp\'er\'ees de $G(F)$ est union disjointe de $L$-paquets $\Pi^{G}(\varphi)$ index\'es par les $\varphi\in \Phi_{temp}(G)$. Pour un tel $\varphi$, on conjecture que l'ensemble $\Pi^{G}(\varphi)$ est en bijection avec ${\cal E}^{G}(\varphi)$, on notera $\epsilon\mapsto \sigma(\varphi,\epsilon)$ cette bijection. Ces param\'etrages doivent \^etre compatibles avec l'endoscopie entre $G$ et ses groupes endoscopiques elliptiques. Ceux-ci sont des produits $G_{+}\times G_{-}$ de groupes sp\'eciaux orthogonaux quasi-d\'eploy\'es d'espaces quadratiques  de dimensions paires, leurs dimensions et discriminants devant satisfaire certaines \'egalit\'es. D'autre part, dans le cas o\`u $\mu(G)=1$, les param\'etrages doivent \^etre compatibles avec l'endoscopie entre $G$ et le groupe $GL(d)$ tordu.

Revenons un instant sur  la d\'efinition de $\mu(G)$. Comme on le sait, pour un discriminant $\delta$ fix\'e, il y a deux classes d'isomorphisme d'espaces quadratiques $(V,q)$ de dimension $d$ et de discriminant $\delta(q)=\delta$ (du moins si $d\geq4$). Pour $\delta=1$, les groupes sp\'eciaux orthogonaux de ces deux espaces sont diff\'erents: l'un est d\'eploy\'e et l'autre n'est pas quasi-d\'eploy\'e. Par contre, si $\delta\not=1$, les deux espaces ont le m\^eme groupe sp\'ecial orthogonal, qui est quasi-d\'eploy\'e. La division traditionnelle entre groupes quasi-d\'eploy\'es et groupes non quasi-d\'eploy\'es n'est pas assez fine pour notre propos et le signe $\mu(G)$ s'introduit dans les preuves pour distinguer les deux classes d'espaces quadratiques dont $G$ est le groupe sp\'ecial orthogonal. Remarquons que, selon le choix de $\mu(G)$, on obtient des param\'etrages des m\^emes paquets $\Pi^G(\varphi)$ par des ensembles de caract\`eres diff\'erents. La relation entre ces param\'etrages est facile \`a expliciter, cf. 4.6

Pour que les conjectures aient un sens, il faut d\'efinir pr\'ecis\'ement les notions de transfert, c'est-\`a-dire fixer des facteurs de transfert. C'est ce que l'on fait en 1.7 et 1.8. D'autre part, notre \'enonc\'e des conjectures laisse place \`a des constantes non pr\'ecis\'ees (par exemple la constante $c^{G'}(\varphi)$ ci-dessus). Un r\'esultat pr\'eliminaire est que, quitte \`a modifier les param\'etrages, on peut pr\'eciser ces constantes (lemme 4.8). Une fois cela fait, les param\'etrages sont enti\`erement d\'etermin\'es pour le groupe $G'$ et le sont presque pour le groupe $G$, le "presque" provenant du probl\`eme d\'elicat de la distinction entre conjugaison par le groupe sp\'ecial orthogonal et conjugaison par le groupe orthogonal tout entier, cf. ci-dessous. Les param\'etrages pour le groupe $G'$ sont ind\'ependants de l'espace $(V,q)$. Par contre, ceux pour le groupe $G$ d\'ependent, sinon de l'espace $(V',q')$, du moins de l'\'el\'ement $\nu_{0}\in F^{\times}/F^{\times,2}$ que l'on a d\'efini ci-dessus. Cela parce que cet \'el\'ement nous sert \`a normaliser les facteurs de transfert.

Posons encore une d\'efinition. Pour deux entiers naturels $N$ et $N'$, avec $N'$ pair, soient $\varphi:W_{DF}\to O(N,{\mathbb C})$ et $\varphi':W_{DF}\to Sp(N',{\mathbb C})$ deux homomorphismes temp\'er\'es. Par la correspondance de Langlands, on leur associe des repr\'esentations irr\'eductibles $\pi(\varphi)$ de $GL(N,F)$ et $\pi(\varphi')$ de $GL(N',F)$. On pose
$$E(\varphi,\varphi')=(\delta(\varphi),-1)_{F}^{N'/2}\epsilon(1/2,\pi(\varphi)\times \pi(\varphi'),\psi_{F}).$$
Le premier terme est un symbole de Hilbert. Le second est le facteur $\epsilon$ de Jacquet, Piatetski-Shapiro et Shalika, d\'efini \`a l'aide d'un caract\`ere $\psi_{F}$ de $F$. Il ne d\'epend pas de ce caract\`ere et $E(\varphi,\varphi')$ est un \'el\'ement de $\{\pm 1\}$.

Cela \'etant, soient $\varphi\in \Phi_{temp}(G)$ et $\varphi'\in \Phi_{temp}(G')$. On d\'ecompose $\varphi$ et $\varphi'$ sous la forme (1), en ajoutant des $'$ aux notations concernant $\varphi'$. On d\'efinit un caract\`ere $\boldsymbol{\epsilon}'$ de $S_{\varphi'}/S_{\varphi'}^0=({\mathbb Z}/2{\mathbb Z})^{(I')^{symp}}$ par
$$\boldsymbol{\epsilon}'((e_{i'})_{i'\in (I')^{symp}})=\prod_{i'\in (I')^{symp}}E(\varphi,\varphi'_{i'})^{e_{i'}}.$$
On d\'efinit un caract\`ere $\boldsymbol{\epsilon}$ de $({\mathbb Z}/2{\mathbb Z})^{I^{orth}}$, que l'on restreint en un caract\`ere de $S_{\varphi}/S_{\varphi}^0$, par
$$\boldsymbol{\epsilon}((e_{i})_{i\in I^{orth}})=\prod_{i\in I^{orth}}E(\varphi_{i},\varphi')^{e_{i}}.$$
Dans le cas o\`u $\delta(q)=1$, on v\'erifie que $G$ et $G'$ sont simultan\'ement d\'eploy\'es o\`u non quasi-d\'eploy\'es. Autrement dit $\mu(G)=\mu(G')$. Dans le cas o\`u $\delta(q)\not=1$, on choisit d\'esormais $\mu(G)=\mu(G')$. On v\'erifie que
$$\boldsymbol{\epsilon}(z_{\varphi})=\boldsymbol{\epsilon}'(z_{\varphi'})=E(\varphi,\varphi').$$
Si $E(\varphi,\varphi')=\mu(G')$, on a $\boldsymbol{\epsilon}\in {\cal E}^G(\varphi)$ et $\boldsymbol{\epsilon}'\in {\cal E}^{G'}(\varphi')$. Si $E(\varphi,\varphi')=-\mu(G')$, on a $\boldsymbol{\epsilon}\not\in {\cal E}^G(\varphi)$ et $\boldsymbol{\epsilon}'\not\in {\cal E}^{G'}(\varphi')$.

On admet la validit\'e des conjectures ci-dessus. On en utilise la forme pr\'ecis\'ee  par le lemme 4.8. On admet aussi certains r\'esultats de la formule des traces locale tordue afin d'assurer la validit\'e des r\'esultats de [W3]. Notre r\'esultat principal est le th\'eor\`eme 4.9(i)  dont voici l'\'enonc\'e.

\ass{Th\'eor\`eme}{Soient $\varphi\in \Phi_{temp}(G)$ et $\varphi'\in \Phi_{temp}(G')$. Si $E(\varphi,\varphi')=-\mu(G')$, on a $m(\sigma,\sigma')=0$ pour tous $\sigma\in \Pi^G(\varphi)$, $\sigma'\in \Pi^{G'}(\varphi')$. Si $E(\varphi,\varphi')=\mu(G')$, on a 
$$m(\sigma(\varphi,\boldsymbol{\epsilon}),\sigma'(\varphi',\boldsymbol{\epsilon}'))=1$$
 et $m(\sigma,\sigma')=0$ pour tous $\sigma\in \Pi^G(\varphi)$, $\sigma'\in \Pi^{G'}(\varphi')$ tels que $(\sigma,\sigma')\not=(\sigma(\varphi,\boldsymbol{\epsilon}),\sigma'(\varphi',\boldsymbol{\epsilon}'))$.}

C'est la conjecture 6.9 de [GP2], limit\'ee aux repr\'esentations temp\'er\'ees. On a admis beaucoup de choses. En ce qui concerne la formule des traces locale tordue, on a tout lieu de penser qu'elle sera d\'emontr\'ee prochainement. Les conjectures de classification sont beaucoup plus profondes. Celles concernant le groupe $G'$ sont annonc\'ees par Arthur, du moins si $G'$ est d\'eploy\'e. Pour le groupe $G$, dans le cas o\`u celui-ci est quasi-d\'eploy\'e, Arthur annonce un r\'esultat un peu moins pr\'ecis, o\`u on ne distingue pas deux repr\'esentations conjugu\'ees par le groupe orthogonal tout entier. Nous ignorons si la forme plus pr\'ecise \'enonc\'ee ci-dessus pourra se d\'eduire de son r\'esultat. Mais le th\'eor\`eme ci-dessus reste valable, sous une forme affaiblie, si l'on admet des conjectures affaiblies conformes aux annonces d'Arthur ((ii) du th\'eor\`eme 4.9).

 Un mot sur la d\'emonstration. Soient $\varphi$, $\varphi'$ comme dans l'\'enonc\'e ci-dessus et soient $\sigma\in \Pi^G(\varphi)$, $\sigma'\in \Pi^{G'}(\varphi')$. Dans [W1], on a calcul\'e $m(\sigma,\sigma')$ par une formule int\'egrale o\`u interviennent les caract\`eres de $\sigma$ et $\sigma'$. En utilisant l'endoscopie ordinaire entre $G$, $G'$ et leurs groupes endoscopiques, on peut exprimer ces caract\`eres \`a l'aide de caract\`eres stables (c'est-\`a-dire les $\Theta^{G'}(\varphi')$ ci-dessus), le prix \`a payer \'etant que ces caract\`eres ne vivent pas sur les groupes de d\'epart, mais sur leurs groupes endoscopiques (c'est le principe de base de l'endoscopie). Par endoscopie tordue, ces caract\`eres stables s'expriment \`a l'aide de caract\`eres sur des groupes tordus $\tilde{GL}(d'_{+})$, $\tilde{GL}(d'_{-})$ etc.... On obtient ainsi une expression de $m(\sigma,\sigma')$ en termes de tels caract\`eres. Dans [W3], on a d\'emontr\'e une formule, parall\`ele \`a celle de [W1], qui calcule des facteurs $\epsilon$ de paires de repr\'esentations de groupes lin\'eaires en termes du m\^eme genre de caract\`eres. Il s'av\`ere  qu'\`a l'aide de cette formule, on peut transformer l'expression obtenue pour $m(\sigma,\sigma')$ en une autre o\`u  les int\'egrales de caract\`eres tordus disparaissent et sont remplac\'ees par des facteurs $\epsilon$ de paires. Il est alors ais\'e d'en d\'eduire le th\'eor\`eme, puisque celui-ci dit justement que $m(\sigma,\sigma')$ se calcule \`a l'aide de tels facteurs.
 
 Je remercie vivement C. Moeglin. Sans ses explications, cet article serait faux.

\section{Param\'etrages}
\bigskip

\subsection{Notations g\'en\'erales}

Soit $F$ un corps local non archim\'edien de caract\'eristique nulle. On note $\vert .\vert _{F}$ sa valeur absolue usuelle. On fixe une cl\^oture alg\'ebrique $\bar{F}$ de $F$ et un caract\`ere $\psi_{F}:F\to {\mathbb C}^{\times}$, continu et non trivial.  On note $(\alpha,\beta)_{F}$ le symbole de Hilbert quadratique de deux \'el\'ements $\alpha,\beta\in F^{\times}$. Pour toute extension quadratique $E$ de $F$, on note $\tau_{E/F}$ l'unique \'el\'ement non trivial du groupe de Galois $Gal(E/F)$ et $sgn_{E/F}$ le caract\`ere quadratique de $F^{\times}$ dont le noyau est le groupe des normes $Norm_{E/F}(E^{\times})$. 

 Pour tout espace topologique $X$ localement compact et totalement discontinu, on note $C_{c}^{\infty}(X)$ l'espace des fonctions de $X$ dans ${\mathbb C}$ qui sont localement constantes et \`a support compact. Consid\'erons un groupe $G$ qui agit sur un ensemble $X$. Pour un sous-ensemble $Y$ de $X$, on note $Norm_{G}(Y)$ le sous-groupe des \'el\'ements de $G$ qui conservent $Y$ et $Z_{G}(Y)$ le sous-groupe des \'el\'ements de $G$ qui fixent tout point de $Y$.
 
 Soit $G$ un groupe r\'eductif d\'efini sur $F$. Tous les sous-groupes que nous consid\'ererons seront suppos\'es d\'efinis sur $F$. On note $\mathfrak{g}$  l'alg\`ebre de Lie de $G$. On note $Nil(\mathfrak{g}(F))$ l'ensemble des orbites nilpotentes dans $\mathfrak{g}(F)$. On note $Temp(G(F))$ l'ensemble des classes d'isomorphie de repr\'esentations admissibles irr\'eductibles et temp\'er\'ees de $G(F)$. Pour un \'el\'ement $\pi\in Temp(G(F))$, on note $E_{\sigma}$ un espace (complexe) dans lequel $\sigma$ se r\'ealise.

\bigskip

\subsection{Mesures}

Pour tout tore $T$ d\'efini sur $F$, on note $A_{T}$ le plus grand sous-tore de $T$ d\'eploy\'e sur $F$. On munit $A_{T}(F)$ de la mesure de Haar pour laquelle le plus grand sous-groupe compact de $A_{T}(F)$ est de mesure $1$. On munit $A_{T}(F)\backslash T(F)$ de la mesure de Haar de masse totale $1$. On munit $T(F)$ de la mesure de Haar telle que, pour $f\in C_{c}^{\infty}(T(F))$, on ait l'\'egalit\'e
$$\int_{T(F)}f(t)dt\,=\int_{A_{T}(F)\backslash T(F)}\int_{A_{T}(F)}f(a\bar{t})da\,d\bar{t}.$$
On dit que $T$ est anisotrope si $A_{T}=\{1\}$. 

Soit $G$ un groupe r\'eductif connexe d\'efini sur $F$. On note $G_{reg}$ le sous-ensemble des \'el\'ements semi-simples fortement r\'eguliers de $G$ et $G_{reg}(F)/conj$ l'ensemble des classes de conjugaison par $G(F)$ dans $G_{reg}(F)$. On dit que deux \'el\'ements  de $G_{reg}(F)$ sont stablement conjugu\'es si et seulement s'ils sont conjugu\'es par un \'el\'ement de $G(\bar{F})$. On note $G_{reg}(F)/stconj$ l'ensemble des classes de conjugaison stable dans $G_{reg}(F)$. Il y a un diagramme naturel
$$\begin{array}{ccccc}&&G_{reg}(F)&&\\ &\phi_{G}\,\swarrow&&\searrow\,\phi^{st}_{G}&\\G_{reg}(F)/conj&&\stackrel{\phi^{st}_{G/conj}}{\to}&&G_{reg}(F)/stconj\\ \end{array}$$
On peut munir  l'ensemble $G_{reg}(F)/conj$ d'une structure de vari\'et\'e analytique sur $F$ et d'une mesure caract\'eris\'ees par la propri\'et\'e suivante.  Soit $\gamma\in G_{reg}(F)$, notons $G_{\gamma}$ son commutant dans $G$, qui est un sous-tore maximal de $G$. Soit $\omega$ un voisinage ouvert de $1$ dans $G_{\gamma}(F)$. Alors, si $\omega$ est assez petit, l'application $t\mapsto \phi_{G}(t\gamma)$ de $\omega$ dans $G_{reg}(F)/conj$ est un isomorphisme analytique de $\omega$ sur un voisinage ouvert de $\phi_{G}(\gamma)$ dans $G_{reg}(F)/conj$ et cet isomorphisme pr\'eserve les mesures. On peut de la m\^eme fa\c{c}on munir l'ensemble $G_{reg}(F)/stconj$ d'une structure de vari\'et\'e analytique sur $F$ et d'une mesure. Alors l'application $\phi^{st}_{G/conj}$ du diagramme ci-dessus est un rev\^etement analytique qui pr\'eserve localement les mesures. 

Pour toute fonction $f$ d\'efinie sur $G_{reg}(F)$ et invariante par conjugaison, resp. par conjugaison stable, on note encore $f$ la fonction sur $G_{reg}(F)/conj$, resp. $G_{reg}(F)/stconj$, qui s'en d\'eduit.

Pour un \'el\'ement semi-simple $\gamma$ de $G(F)$, on pose
$$D^G(\gamma)=\vert det((1-ad(\gamma))_{\vert \mathfrak{g}(F)/\mathfrak{g}_{\gamma}(F)})\vert _{F}.$$

Munissons $G(F)$ d'une mesure de Haar.  Soit $f\in C_{c}^{\infty}(G(F))$.  
 Pour $\gamma\in G_{reg}(F)$, on pose
$$\Phi(\gamma,f)=\int_{G_{\gamma}(F)\backslash G(F)}f(g^{-1}\gamma g)dg.$$
 Pour $y\in G_{reg}(F)/stconj$, on d\'efinit $\Phi^{st}(y,f)=\sum_{x}\Phi(x,f)$, o\`u $x$ parcourt la fibre de $\phi_{G/conj}$ au-dessus de $y$. Avec ces d\'efinitions, la formule de Weyl prend l'une ou l'autre des formes suivantes:
$$\int_{G(F)}f(g)dg=\int_{G_{reg}(F)/conj}\Phi(x,f)D^G(x)dx=\int_{G_{reg}(F)/stconj}\Phi^{st}(y,f)D^G(y)dy.$$

Soit $(M,\tilde{M})$ un groupe tordu d\'efini sur $F$. Cela signifie que $M$ est un groupe r\'eductif connexe d\'efini sur $F$, que $\tilde{M}$ est une vari\'et\'e alg\'ebrique d\'efinie sur $F$ telle que $\tilde{M}(F)$ soit non vide et qu'il y a deux actions \`a droite et \`a gauche de $M$ sur $\tilde{M}$, not\'ees
$$\begin{array}{ccc}M\times \tilde{M}\times M&\to& \tilde{M}\\ (m,\tilde{m},m')&\mapsto&m\tilde{m}m',\\ \end{array}$$
de sorte que, pour chacune des actions, $\tilde{M}$ soit un espace principal homog\`ene sous $M$. Pour $\tilde{m}\in \tilde{M}$, on note $\theta_{\tilde{m}}$ l'automorphisme de $M$ tel que $\tilde{m}m=\theta_{\tilde{m}}(m)\tilde{m}$ pour tout $m\in M$.  On note $Z_{M}(\tilde{m})$ le sous-groupe des points fixes de $\theta_{\tilde{m}}$, c'est-\`a-dire le groupe des $m\in M$ tels que $m\tilde{m}m^{-1}=\tilde{m}$. On note $M_{\tilde{m}}$ la composante neutre de $Z_{M}(\tilde{m})$. On dit que $\tilde{m}$ est semi-simple fortement r\'egulier si $Z_{M}(\tilde{m})$ est ab\'elien et $M_{\tilde{m}}$ est un tore. On note $\tilde{M}_{reg}$ l'ensemble des \'el\'ements fortement r\'eguliers de $\tilde{M}$ et $\tilde{M}_{reg}(F)/conj$ l'ensemble des classes de conjugaison par $M(F)$ dans $\tilde{M}_{reg}(F)$, en appelant conjugaison l'action $(m,\tilde{m})\mapsto m\tilde{m}m^{-1}$ de $M$ sur $\tilde{M}$. On dit que deux \'el\'ements de $\tilde{M}_{reg}(F)$ sont stablement conjugu\'es si et seulement s'ils sont conjugu\'es par un \'el\'ement de $M(\bar{F})$. On note $\tilde{M}_{reg}(F)/stconj$ l'ensemble des classes de conjugaison stable dans $\tilde{M}_{reg}(F)$. Il y a un diagramme naturel
$$\begin{array}{ccccc}&&\tilde{M}_{reg}(F)&&\\ &\phi_{\tilde{M}}\,\swarrow&&\searrow\,\phi^{st}_{\tilde{M}}&\\\tilde{M}_{reg}(F)/conj&&\stackrel{\phi^{st}_{\tilde{M}/conj}}{\to}&&\tilde{M}_{reg}(F)/stconj\\ \end{array}$$

  On peut munir  l'ensemble $\tilde{M}_{reg}(F)/conj$ d'une structure de vari\'et\'e analytique sur $F$ et d'une mesure caract\'eris\'ees par la propri\'et\'e suivante.  Soient $\tilde{\gamma}\in \tilde{M}_{reg}(F)$ et $\omega$ un voisinage ouvert de $1$ dans $M_{\tilde{\gamma}}(F)$. Alors, si $\omega$ est assez petit, l'application $t\mapsto \phi_{\tilde{M}}(t\tilde{\gamma})$ de $\omega$ dans $\tilde{M}_{reg}(F)/conj$ est un isomorphisme analytique de $\omega$ sur un voisinage ouvert de $\phi_{\tilde{M}}(\tilde{\gamma})$ dans $\tilde{M}_{reg}(F)/conj$ et cet isomorphisme pr\'eserve les mesures. On munit de m\^eme l'ensemble $\tilde{M}_{reg}(F)/stconj$ d'une structure de vari\'et\'e analytique sur $F$ et d'une mesure. L'application $\phi^{st}_{\tilde{M}/conj}$ est un rev\^etement analytique qui pr\'eserve localement les mesures.
  
  Pour un \'el\'ement semi-simple $\tilde{\gamma}\in \tilde{M}(F)$, posons
  $$D^{\tilde{M}}(\gamma)=\vert det((1-\theta_{\tilde{\gamma}}))_{\vert \mathfrak{m}(F)/\mathfrak{m}_{\tilde{\gamma}}(F)})\vert _{F}.$$

  Munissons $M(F)$ d'une mesure de Haar.  On en d\'eduit une mesure sur $\tilde{M}(F)$ en consid\'erant cet ensemble comme un espace principal homog\`ene pour l'action \`a gauche, resp. \`a droite, de $M(F)$, et on obtient en fait la m\^eme mesure quelle que soit l'action choisie. Soit $\tilde{f}\in C_{c}^{\infty}(\tilde{M}(F))$.    Pour $\tilde{\gamma}\in \tilde{M}_{reg}(F)$, posons
  $$\Phi(\tilde{\gamma},\tilde{f})=[Z_{M}(\tilde{\gamma})(F):M_{\tilde{\gamma}}(F)]^{-1}\int_{M_{\tilde{\gamma}}(F)\backslash M(F)}\tilde{f}(m^{-1}\tilde{\gamma} m)dm.$$
 Alors, avec les m\^emes conventions de notations que dans le cas non tordu, la formule de Weyl prend la forme:
$$\int_{\tilde{M}(F)}\tilde{f}(\tilde{m})d\tilde{m}=\int_{\tilde{M}_{reg}(F)/conj}\Phi(\tilde{x},\tilde{f})D^{\tilde{M}}(\tilde{x})d\tilde{x}.$$

\bigskip

\subsection{Espace de param\`etres}

  On appelle ici extension alg\'ebrique de $F$ un sous-corps de $\bar{F}$ contenant $F$. Notons $\underline{\Xi}$ l'ensemble des familles $\xi=(I,(F_{\pm i})_{i\in I}, (F_{i})_{i\in I},(y_{i})_{i\in I})$ v\'erifiant les conditions suivantes:

$\bullet$ $I$ est un sous-ensemble fini de ${\mathbb N}$;

$\bullet$ pour tout $i\in I$, $F_{\pm i}$ est une extension finie de $F$ et $F_{i}$ est une $F_{\pm i}$-alg\`ebre commutative de dimension $2$; c'est-\`a-dire $F_{i}$ est une extension quadratique de $F_{\pm i}$ ou $F_{i}=F_{\pm i}\oplus F_{\pm i}$; on note $\tau_{i}$ l'unique automorphisme non trivial de $F_{\pm i}$-alg\`ebre de $F_{i}$;

$\bullet$ pour tout $i\in I$, $y_{i}$ est un \'el\'ement de $F_{i}^{\times}$ tel que $y_{i}\tau_{i}(y_{i})=1$.

Pour une telle famille, on note $I^*$ l'ensemble des $i\in I$ tels que $F_{i}$ soit un corps. Pour $i\in I$, on fixe $\delta_{i}\in F_{\pm i}^{\times}$ tel que $F_{i}=F_{\pm i}(\sqrt{\delta_{i}})$, avec la convention que $\delta_{i}$ appartient au groupe des carr\'es $F_{\pm i}^{\times,2}$ si $i\not\in I^*$. On note $\delta_{\xi}$ l'image de $\prod_{i\in I}Norm_{F_{\pm i}/F}(\delta_{i})$ dans $F^{\times}/F^{\times,2}$. Cet \'el\'ement ne d\'epend \'evidemment pas des choix des $\delta_{i}$. On pose 
$$d_{\xi}=\sum_{i\in I}[F_{i}:F].$$
 
Un isomorphisme entre deux \'el\'ements 
$$\xi=(I,(F_{\pm i})_{i\in I}, (F_{i})_{i\in I},(y_{i})_{i\in I})\text{ et }\xi'=(I',(F'_{\pm i'})_{i'\in I'}, (F'_{i'})_{i'\in I'},(y'_{i'})_{i'\in I'})$$
 de $\underline{\Xi}$ est une famille $(\iota,(\iota_{\pm i})_{i\in I},(\iota_{i})_{i\in I})$, o\`u 

$\bullet$ $\iota:I\to I'$ est une bijection;

$\bullet$ pour tout $i\in I$, $\iota_{\pm i}:F_{\pm i}\to F'_{\pm \iota(i)}$ est un isomorphisme d\'efini sur $F$ et $\iota_{i}:F_{i}\to F'_{\iota(i)}$ est un isomorphisme compatible avec $\iota_{\pm i}$ en un sens \'evident;

$\bullet$ pour tout $i\in I$, $\iota_{i}(y_{i})=y'_{\iota(i)}$.

En particulier, tout \'el\'ement $\xi$ poss\`ede un automorphisme "identit\'e". On dit que $\xi$ est r\'egulier si l'identit\'e est le seul automorphisme de $\xi$. On note $\underline{\Xi}_{reg}$ l'ensemble des \'el\'ements r\'eguliers de $\underline{\Xi}$ et $\Xi_{reg}$ le quotient de $\underline{\Xi}_{reg}$ obtenu en identifiant les \'el\'ements isomorphes. En pratique, on notera un \'el\'ement de $\Xi_{reg}$ comme un \'el\'ement de $\underline{\Xi}_{reg}$ qui le repr\'esente.

  Soit  $\xi=(I,(F_{\pm i})_{i\in I}, (F_{i})_{i\in I},(y_{i})_{i\in I})\in \Xi_{reg}$. Introduisons le tore $T_{\xi}$ d\'efini sur $F$ tel que $T_{\xi}(F)=\prod_{i\in I}\{t_{i}\in F_{i}^{\times}; t_{i}\tau_{i}(t_{i})=1\}$. Soit $\omega$ un voisinage ouvert de l'unit\'e dans $T_{\xi}(F)$. Si $\omega$ est assez petit, l'application
$$(t_{i})_{i\in I}\mapsto ( I,(F_{\pm i})_{i\in I}, (F_{i})_{i\in I},(y_{i}t_{i})_{i\in I})$$
est une injection de $\omega$ dans $\Xi_{reg}$. On peut munir  et on munit $\Xi_{reg}$ d'une structure de vari\'et\'e analytique sur $F$ et d'une mesure de sorte que cette injection soit un isomorphisme qui pr\'eserve les mesures.  

Si on fixe un entier pair $d$, ou un \'el\'ement $\delta\in F^{\times}/F^{\times,2}$, ou les deux, on d\'efinit les variantes $\Xi_{d}$, $\Xi_{\delta}$, $\Xi_{d,\delta}$, $\Xi_{reg,d}$ etc... des objets pr\'ec\'edents en se limitant aux $\xi$ tels que $d_{\xi}=d$, ou $\delta_{\xi}=\delta$, ou $d_{\xi}=d$ et $\delta_{\xi}=\delta$.

Soit $\xi=(I,(F_{\pm i})_{i\in I}, (F_{i})_{i\in I},(y_{i})_{i\in I})$ un \'el\'ement de $\Xi_{reg}$. On pose
$$C(\xi)=\prod_{i\in I^*}F_{\pm i}^{\times}/Norm_{F_{i}/F_{\pm i}}(F_{i}^{\times}).$$
Ce groupe s'identifie naturellement \`a $\{\pm 1\}^{I^*}$. On note $C(\xi)^1$ le sous-groupe qui s'identifie au sous-groupe des \'el\'ements $(\epsilon_{i})_{i\in I^*}\in \{\pm 1\}^{I^*}$ tels que $\prod_{i\in I^*}\epsilon_{i}=1$. Si $I^*\not=\emptyset$, on note $C(\xi)^{-1}=C(\xi)\setminus C(\xi)^1$.

Pour $i\in I$, notons $\Gamma(y_{i})$ l'ensemble des \'el\'ements $\gamma_{i}\in F_{i}^{\times}$ tels que $\gamma_{i}\tau_{i}(\gamma_{i})^{-1}=y_{i}$. On pose
$$\Gamma_{pair}(\xi)=\prod_{i\in I^*}\Gamma(y_{i})/Norm_{F_{i}/F_{\pm i}}(F_{i}^{\times}),$$
$$\Gamma_{imp}(\xi)=(F^{\times}/F^{\times,2})\times \prod_{i\in I^*}\Gamma(y_{i})/Norm_{F_{i}/F_{\pm i}}(F_{i}^{\times}).$$

En pratique,  pour un \'el\'ement $c=(c_{i})_{i\in I^*}$ de $C(\xi)$, on identifiera chaque $c_{i}$ \`a un \'el\'ement de $F_{\pm i}^{\times}$ qui le repr\'esente. On fera de m\^eme pour les \'el\'ements de $\Gamma_{pair}(\xi)$ ou $\Gamma_{imp}(\xi)$.

\subsection{Param\'etrage des classes de conjugaison}
Soient $V$ un espace vectoriel de dimension finie $n$ sur $F$ et $q$ une forme bilin\'eaire sym\'etrique non d\'eg\'en\'er\'ee sur $V$ (on appellera le couple $(V,q)$ un espace quadratique).  On d\'efinit son discriminant $\delta(q)=(-1)^{[n/2]}det(q)\in F^{\times}/F^{\times,2}$. On note $G$ le groupe sp\'ecial orthogonal du couple $(V,q)$.

 Supposons d'abord $n$ impair. Soient $\xi=(I,(F_{\pm i})_{i\in I}, (F_{i})_{i\in I},(y_{i})_{i\in I})\in \Xi_{reg,n-1}$ et $c=(c_{i})_{i\in I}\in C(\xi)$. Posons
$$V_{\xi}=\oplus_{i\in I}F_{i}$$
et munissons cet espace de la forme bilin\'eaire sym\'etrique et non d\'eg\'en\'er\'ee $q_{\xi,c}$ d\'efinie par
$$q_{\xi,c}(\sum_{i\in I}v_{i},\sum_{i\in I}v'_{i})=\sum_{i\in I}trace_{F_{i}/F}(c_{i}\tau_{i}(v_{i})v'_{i}).$$
Consid\'erons la condition: il existe une droite $D$ sur $F$ et une forme bilin\'eaire sym\'etrique non d\'eg\'en\'er\'ee $q_{D}$ sur $D$ telles que les espaces quadratiques $(V,q)$ et $(D\oplus V_{\xi},q_{D}\oplus q_{\xi,c})$ soient isomorphes. On montre qu'il existe une unique classe $C(\xi)^{\epsilon}\subset C(\xi)$  (avec $\epsilon=\pm 1$ d\'ependant de $\xi$ et $(V,q)$) telle que cette condition soit v\'erifi\'ee si et seulement si $c$ appartient \`a cette classe. Supposons-la v\'erifi\'ee. La droite quadratique $(D,q_{D})$ est alors unique \`a isomorphisme pr\`es. Fixons un isomorphisme entre $(V,q)$ et $(D\oplus V_{\xi},q_{D}\oplus q_{\xi,c})$. Introduisons l'\'el\'ement $x\in G(F)$ qui, modulo cet isomorphisme, agit sur chaque $F_{i}$ par multiplication par $y_{i}$ et sur $D$ par l'identit\'e. C'est un \'el\'ement de $G_{reg}(F)$ dont la classe de conjugaison est bien d\'etermin\'ee. Inversement, toute classe de conjugaison dans $G_{reg}(F)$ est obtenue par ce proc\'ed\'e. On obtient ainsi une application de $G_{reg}(F)/conj$ dans $\Xi_{reg,n-1}$ qui se factorise en un diagramme
$$ \begin{array}{ccccc}G_{reg}(F)/conj&&\stackrel{\phi^{st}_{G/conj}}{\to}&&G_{reg}(F)/stconj\\ &p_{G}\,\searrow&& \swarrow\,p^{st}_{G}&\\ &&\Xi_{reg,n-1}&&\\ \end{array}$$
L'application $p^{st}_{G}$ est un isomorphisme analytique qui pr\'eserve les mesures. La fibre en un point $\xi$ de l'application $p_{G}$ s'identifie \`a une classe  $C(\xi)^{\epsilon}\subset C(\xi)$.   

Supposons maintenant que $n$ est pair et $n\geq2$.  Posons $\delta=\delta(q) $.  Soient $\xi=(I,(F_{\pm i})_{i\in I}, (F_{i})_{i\in I},(y_{i})_{i\in I})\in \Xi_{reg,n,\delta}$ et $c=(c_{i})_{i\in I}\in C(\xi)$. D\'efinissons $V_{\xi}$ et $q_{\xi,c}$ comme pr\'ec\'edemment. Consid\'erons la condition: les espaces quadratiques $(V,q)$ et $(V_{\xi},q_{\xi,c})$ sont isomorphes. Il existe une unique classe  $C(\xi)^{\epsilon}\subset C(\xi)$  telle que cette condition soit v\'erifi\'ee si et seulement si $c$ appartient \`a cette classe. Supposons-la v\'erifi\'ee et fixons un isomorphisme entre nos deux espaces quadratiques. Introduisons l'\'el\'ement $x\in G(F)$ qui, modulo cet isomorphisme, agit sur chaque $F_{i}$ par multiplication par $y_{i}$. C'est un \'el\'ement de $G_{reg}(F)$. Il y a deux diff\'erences avec le cas $n$ impair. D'une part, c'est seulement la classe de conjugaison de $x$ par le groupe orthogonal tout entier qui est bien d\'etermin\'ee. Or cette classe se d\'ecompose en deux classes de conjugaison par $G(F)$. D'autre part, on n'obtient par ce proc\'ed\'e que les classes d'\'el\'ements $x\in G_{reg}(F)$ qui n'ont pas de valeurs propres \'egales \`a $\pm 1$. On \'elimine ce deuxi\`eme probl\`eme en modifiant la d\'efinition de $G_{reg}$.

\ass{Convention}{Pour un groupe sp\'ecial orthogonal $G$ d'un espace quadratique de dimension paire, on note dor\'enavant $G_{reg}$ l'ensemble des \'el\'ements semi-simples r\'eguliers de $G$ qui n'ont aucune valeur propre \'egale \`a $\pm 1$.}

On obtient alors le diagramme
$$\begin{array}{ccccc}G_{reg}(F)/conj&&\stackrel{\phi^{st}_{G/conj}}{\to}&&G_{reg}(F)/stconj\\ &p_{G}\,\searrow&& \swarrow\,p^{st}_{G}&\\ &&\Xi_{reg,n,\delta}&&\\ \end{array}$$
L'application $p^{st}_{G}$ est un rev\^etement qui pr\'eserve localement les mesures. Ses fibres sont d'ordre $2$ et form\'ees de deux classes de conjugaison stable qui sont conjugu\'ees par l'action du groupe orthogonal. Soit $\xi\in \Xi_{reg,n,\delta}$ et $y\in G_{reg}(F)/stconj$ tel que $p^{st}_{G}(y)=\xi$. La fibre de $\phi^{st}_{G/conj}$ au-dessus de $y$ s'identifie \`a une classe $C(\xi)^{\epsilon}\subset C(\xi)$.

 {\bf Variante.} Notons $G^+$ le groupe orthogonal de $(V,q_{V})$. On d\'efinit les ensembles $G_{reg}(F)/conj^+$ et $G_{reg}(F)/stconj^+$ en rempla\c{c}ant dans les d\'efinitions conjugaison par $G$ par conjugaison par $G^+$. On a alors un diagramme
 $$\begin{array}{ccccc}G_{reg}(F)/conj^+&&\stackrel{\phi^{st}_{G/conj^+}}{\to}&&G_{reg}(F)/stconj^+\\ &p^+_{G}\,\searrow&& \swarrow\,p^{st+}_{G}&\\ &&\Xi_{reg,n,\delta}&&\\ \end{array}$$
 o\`u cette fois, $p^{st+}_{G}$ est un isomorphisme.

Pour simplifier certaines constructions ult\'erieures, il convient de consid\'erer le cas $n=0$. Dans ce cas, on pose $\delta=1$ et $G_{reg}=G=\{1\}$. Les diagrammes ci-dessus se trivialisent, tous les ensembles y figurant ayant un seul \'el\'ement.

Que $d$ soit pair ou impair, notons $G(F)_{ani}$ l'ensemble des \'el\'ements de $G_{reg}(F)$ dont le commutant est un tore anisotrope. Il s'agit des \'el\'ements appel\'es usuellement elliptiques, sauf dans le cas o\`u $n=2$ et $G$ est d\'eploy\'e, auquel cas tout \'el\'ement est elliptique alors que $G(F)_{ani}$ est vide. Notons $\Xi^*$ l'ensemble des \'el\'ements $\xi=(I,(F_{\pm i})_{i\in I}, (F_{i})_{i\in I},(y_{i})_{i\in I})\in \Xi$ tels que $I=I^*$. On d\'efinit aussi $\Xi^*_{reg}$, $\Xi^*_{n}$ etc... Dans les diagrammes ci-dessus, on peut remplacer les ensembles $G_{reg}(F)$ par $G(F)_{ani}$  et les $\Xi_{reg,n}$ etc... par $\Xi^*_{reg,n}$ etc... Les diagrammes obtenus conservent les m\^emes propri\'et\'es.

 Soit  $U$ un espace vectoriel sur $F$ de dimension finie $n$. On note $M$ le groupe des automorphismes $F$-lin\'eaires de $U$. On note $\tilde{M}$ l'ensemble des formes bilin\'eaires non d\'eg\'en\'er\'ees sur $U$. Le groupe $M$ agit \`a droite et \`a gauche sur $\tilde{M}$ par la formule
$$(m\tilde{m}m')(u,u')=\tilde{m}(m^{-1}u,m'u')$$
pour $m,m'\in M$, $\tilde{m}\in \tilde{M}$ et $u,u'\in U$. Le couple $(M,\tilde{M})$ est un groupe tordu. 

Supposons d'abord $n$ pair. Soit $\xi=(I,(F_{\pm i})_{i\in I}, (F_{i})_{i\in I},(y_{i})_{i\in I})\in \Xi_{reg,n}$ et $\gamma=(\gamma_{i})_{i\in I}\in \Gamma_{pair}(\xi)$. Introduisons le m\^eme espace $V_{\xi}$ que ci-dessus et fixons un isomorphisme $F$-lin\'eaire de $U$ sur $V_{\xi}$. Introduisons l'\'el\'ement $\tilde{x}\in \tilde{M}(F)$ d\'efini, modulo cet isomorphisme, par la formule
$$(1) \qquad \tilde{x}(\sum_{i\in I}u_{i},\sum_{i\in I}u'_{i})=\sum_{i\in I}trace_{F_{i}/F}(\tau_{i}(u_{i})u'_{i}\gamma_{i}).$$
C'est un \'el\'ement de $\tilde{M}_{reg}(F)$ dont la classe de conjugaison est uniquement d\'etermin\'ee. Inversement, toute classe semi-simple r\'eguli\`ere est obtenue par ce proc\'ed\'e. On obtient un diagramme
$$ \begin{array}{ccccc}\tilde{M}_{reg}(F)/conj&&\stackrel{\phi^{st}_{\tilde{M}/conj}}{\to}&&\tilde{M}_{reg}(F)/stconj\\ &p_{\tilde{M}}\,\searrow&& \swarrow\,p^{st}_{\tilde{M}}&\\ &&\Xi_{reg,n}&&\\ \end{array}$$
 La fibre de $p_{\tilde{M}}$ au-dessus de $\xi\in \Xi_{reg,n}$ s'identifie \`a $\Gamma_{pair}(\xi)$.
L'application $p^{st}_{\tilde{M}}$ est un isomorphisme analytique  mais ne pr\'eserve pas localement les mesures. En effet, fixons $\xi=(I,(F_{\pm i})_{i\in I}, (F_{i})_{i\in I},(y_{i})_{i\in I})\in \Xi_{reg,n}$ et construisons $\tilde{x}$ comme ci-dessus.  Le tore $M_{\tilde{x}}$ est \'egal \`a $T_{\xi}$. Soit $\omega$ un voisinage ouvert assez petit de l'unit\'e dans $T_{\xi}(F)$. La mesure sur $\tilde{M}_{reg}(F)/stconj$ a \'et\'e d\'efinie de sorte que l'application $t\mapsto \tilde{x}t$ de $\omega$ dans $\tilde{M}_{reg}(F)/stconj$ pr\'eserve les mesures. Changer $\tilde{x}$ en $\tilde{x}t$ revient \`a changer $\gamma_{i}$ en $\gamma_{i}t_{i}$ dans la formule (1). D'apr\`es la relation $y_{i}=\gamma_{i}\tau_{i}(\gamma_{i})^{-1}$, on a l'\'egalit\'e $p^{st}_{\tilde{M}}(\tilde{x}t)=(I,(F_{\pm i})_{i\in I}, (F_{i})_{i\in I},(y_{i}t_{i}^2)_{i\in I})$ 
Or la mesure sur $\Xi_{reg,n}$ a \'et\'e d\'efinie de sorte que l'application $t=(t_{i})_{i\in I}\mapsto (I,(F_{\pm i})_{i\in I}, (F_{i})_{i\in I},(y_{i}t_{i})_{i\in I})$ pr\'eserve les mesures.  Le jacobien de l'application $p^{st}_{\tilde{M}}$ est donc le m\^eme que celui de l'application $t\mapsto t^2$ de $T_{\xi}(F)$ dans lui-m\^eme, c'est-\`a-dire $\vert 2\vert _{F}^{n/2}$.  

Supposons maintenant $n$ impair. $\xi=(I,(F_{\pm i})_{i\in I}, (F_{i})_{i\in I},(y_{i})_{i\in I})\in \Xi_{reg,n-1}$ et $\gamma=(\gamma_{D},(\gamma_{i})_{i\in I})\in \Gamma_{imp}(\xi)$.  Posons $D=F$ et fixons un isomorphisme de $U$ sur $D\oplus V_{\xi}$. Introduisons l'\'el\'ement $\tilde{x}\in \tilde{M}(F)$ d\'efini, modulo cet isomorphisme, par la formule
$$\tilde{x}(u_{D}+\sum_{i\in I}u_{i},u'_{D}+\sum_{i\in I}u'_{i})=\gamma_{D}u_{D}u'_{D}+\sum_{i\in I}trace_{F_{i}/F}(\tau_{i}(u_{i})u'_{i}\gamma_{i}).$$ 
C'est un \'el\'ement de $\tilde{M}_{reg}(F)$ dont la classe de conjugaison est uniquement d\'etermin\'ee. Inversement, toute classe semi-simple r\'eguli\`ere est obtenue par ce proc\'ed\'e. On obtient le  diagramme
$$ \begin{array}{ccccc}\tilde{M}_{reg}(F)/conj&&\stackrel{\phi^{st}_{\tilde{M}/conj}}{\to}&&\tilde{M}_{reg}(F)/stconj\\ &p_{\tilde{M}}\,\searrow&& \swarrow\,p^{st}_{\tilde{M}}&\\ &&\Xi_{reg,n-1}&&\\ \end{array}$$
La fibre de $p_{\tilde{M}}$ au-dessus de $\xi\in \Xi_{reg,n-1}$ s'identifie \`a $\Gamma_{imp}(\xi)$. L'application $p^{st}_{\tilde{M}}$ est un isomorphisme analytique dont le jacobien vaut  $\vert 2\vert _{F}^{(n-1)/2}$ d'apr\`es le m\^eme calcul que ci-dessus.

Que $n$ soit pair ou impair, notons $\tilde{M}(F)_{ani}$ l'ensemble des \'el\'ements $\tilde{\gamma}\in\tilde{M}_{reg}(F)$  tels que $M_{\tilde{\gamma}}$ soit un tore anisotrope. Dans les diagrammes ci-dessus, on peut remplacer $\tilde{M}_{reg}(F)$ par $\tilde{M}(F)_{ani}$ et les ensembles $\Xi_{reg,n}$ ou $\Xi_{reg,n-1}$ par $\Xi^*_{reg,n}$ ou $\Xi^*_{reg,n-1}$. Les diagrammes obtenus conservent les m\^emes propri\'et\'es.

\bigskip

\subsection{D\'efinition de fonctions sur l'ensemble de param\`etres}

Soit $\xi=(I,(F_{\pm i})_{i\in I},(F_{i})_{i\in I},(y_{i})_{i\in I})\in \Xi_{reg}$.  Pour tout $i\in I$, notons $\Phi_{i}$ l'ensemble des homomorphismes non nuls de $F$-alg\`ebres de $F_{i}$ dans $\bar{F}$. On d\'efinit un polyn\^ome 
$$P_{\xi}(T)=\prod_{i\in I}\prod_{ \phi\in \Phi_{i}}(T-\phi(y_{i})).$$
 Introduisons un espace quadratique $(V_{\xi},q_{\xi})$ tel que $dim(V_{\xi})=d_{\xi}$ et $\delta(q_{\xi})=\delta_{\xi}$. Notons $G_{\xi}$ son groupe sp\'ecial orthogonal. Alors $\xi$ param\`etre deux classes de conjugaison stable dans 
 $G_{\xi,reg}(F)$. Fixons-en une et un \'el\'ement $t$ de cette classe. On pose
 $$\Delta(\xi)=\vert det((1-t)_{\vert V_{\xi}}\vert _{F}.$$
 Ce terme ne d\'epend ni du choix de l'espace $(V_{\xi},q_{\xi})$, ni du choix de $t$.  En fait, il est clair que
 $$\Delta(\xi)=P_{\xi}(1).$$
 Pour un entier $n\geq d_{\xi}$, fixons un espace quadratique $(Z,q_{Z})$ tel que $dim(Z)=n-d_{\xi}$. Notons $G$ le groupe sp\'ecial orthogonal de la somme orthogonale  $(V_{\xi},q_{\xi})\oplus (Z,q_{Z})$. Alors $t$ est un \'el\'ement semi-simple de $G(F)$. Il n'est plus r\'egulier, mais on a n\'eanmoins d\'efini $D^G(t)$. On pose
 $$D^n(\xi)=D^G(t).$$
Ici encore, ce terme ne d\'epend pas des choix. On pose simplement $D(\xi)=D^{d_{\xi}}(\xi)$. On v\'erifie l'\'egalit\'e
$$(1) \qquad D^n(\xi)=D(\xi)\Delta(\xi)^{n-d_{\xi}}.$$

Soient $\xi_{+}=(I_{+}, (F_{\pm i})_{i\in I_{+}},(F_{i})_{i\in I_{+}},(y_{i})_{i\in I_{+}})\in \Xi_{reg}$ et $\xi_{-}=(I_{-}, (F_{\pm i})_{i\in I_{-}},(F_{i})_{i\in I_{-}},(y_{i})_{i\in I_{-}})\in \Xi_{reg}$, o\`u l'on suppose, ainsi qu'il est loisible, que $I_{+}\cap I_{-}=\emptyset$. Posons  $\xi=\xi_{+}\sqcup \xi_{-}$ et $I=I_{+}\sqcup I_{-}$. Supposons $\xi$ r\'egulier.  Soit $c=(c_{i})_{i\in I}\in C(\xi)$.  On note $P'_{\xi}(T)$ le polyn\^ome d\'eriv\'e du polyn\^ome $P_{\xi}$. Pour tout $i\in I$, posons
 $$  C_{i}=(-1)^{d_{\xi}/2} c_{i}P'_{\xi}(y_{i})P_{\xi}(-1)y_{i}^{1-d_{\xi}/2}(y_{i}-1)^{-1}(y_{i}+1).$$
 C'est un \'el\'ement de $F_{\pm i}$. Pour $\nu\in F^{\times}$, on pose
 $$ \Delta_{\nu}(\xi_{+},\xi_{-},c)=\prod_{i\in I^*_{-}}sgn_{F_{i}/F_{\pm i}}(\nu C_{i}).$$
 
     Soit $\xi=(I,(F_{\pm i})_{i\in I},(F_{i})_{i\in I},(y_{i})_{i\in I})\in \Xi_{reg}$.  Soit $\gamma=(
 \gamma_{i})_{i\in I}\in \Gamma_{pair}(\xi)$. Pour $i\in I$, posons
 $$C_{i}=-\gamma_{i}^{-1}P_{\xi}(1)P'_{\xi}(y_{i})y_{i}^{1-d_{\xi}/2}(y_{i}-1).$$
 Soit maintenant $\gamma=(\gamma_{D},(\gamma_{i})_{i\in I})\in \Gamma_{imp}(\xi)$ (la composante $\gamma_{D}$ appartient \`a $F^{\times}/F^{\times,2}$). Pour $i\in I$, posons
  $$C_{i}=\gamma_{D}\gamma_{i}^{-1}P_{\xi}(1)P'_{\xi}(y_{i})y_{i}^{1-d_{\xi}/2}(y_{i}-1).$$
 Dans les deux cas, on pose
 $$\Delta(\xi,\gamma)=\prod_{i\in I^*}sgn_{F_{i}/F_{\pm i}}(C_{i}).$$

\bigskip

\subsection{Endoscopie}

Soit $G$ un groupe r\'eductif connexe d\'efini sur $F$. Fixons une forme quasi-d\'eploy\'ee $\underline{G}$ de $G$ et un torseur int\'erieur $\psi_{G}:G\to \underline{G}$. On fixe une fonction $u:Gal(\bar{F}/F)\to \underline{G}$ telle que $ \psi_{G}\circ\sigma(g)=u(\sigma)\sigma\circ\psi_{G}(g)u(\sigma)^{-1}$ pour tous $g\in G$ et $\sigma\in Gal(\bar{F}/F)$. Le compos\'e de $u$ et de l'application naturelle de $\underline{G}$ dans son groupe adjoint $\underline{G}_{ad}$ est un cocycle. Nous consid\'erons le cas favorable o\`u $u$ lui-m\^eme est un cocycle, \`a valeurs dans $\underline{G}$. Fixons  une paire de Borel \'epingl\'ee de $\underline{G}$, d\'efinie sur $F$. C'est-\`a-dire que l'on fixe un sous-groupe de Borel $\underline{B}$ d\'efini sur $F$, un sous-tore maximal $\underline{T}$ de $\underline{B}$ d\'efini sur $F$. Notons $\Pi$ l'ensemble des racines simples de $\underline{T}$ associ\'e \`a $\underline{B}$. Pour $\alpha\in \Pi$, on fixe un \'el\'ement non nul $E_{\alpha}$ de l'espace radiciel associ\'e \`a $\alpha$ et on suppose que la famille $(E_{\alpha})_{\alpha\in \Pi}$ est stable par l'action de $Gal(\bar{F}/F)$. Remarquons que $N=\sum_{\alpha\in \Pi}E_{\alpha}$ est un \'el\'ement nilpotent r\'egulier de $\underline{\mathfrak{g}}(F)$ et la classe de conjugaison de la paire de Borel \'epingl\'ee est d\'etermin\'ee par celle de $N$. Il nous suffit en fait de fixer la classe de conjugaison de $N$. Consid\'erons enfin une donn\'ee endoscopique $(H,s,{^L\xi})$ de $G$, cf. [LS] 1.2. On peut d\'efinir une application, la correspondance endoscopique, dont l'ensemble de d\'epart est un sous-ensemble de $H_{reg}(F)/stconj$ et l'ensemble d'arriv\'ee est  $G_{reg}(F)/stconj$. Elle pr\'eserve localement les mesures.  Pour $y\in H_{reg}(F)/stconj$ et $x\in G_{reg}(F)/conj$ tels que $\phi^{st}_{G/conj}{x}$ corresponde \`a $y$, on d\'efinit le facteur de transfert $\Delta_{H,G}(y,x)$ (nous en supprimons le terme $\Delta_{IV}$). Pour $f\in C_{c}^{\infty}(G(F))$ et $f^H\in C_{c}^{\infty}(H(F))$, on dit que $f^H$ est un transfert de $f$ si on a l'\'egalit\'e
$$D^H(y)^{1/2}\Phi^{st}(y,f^H)=\sum_{x}\Delta_{H,G}(y,x)D^G(x)^{1/2}\Phi(x,f)$$
pour presque tout $y\in H_{reg}(F)/stconj$, o\`u $x$ parcourt l'ensemble des \'el\'ements de $G_{reg}(F)/conj$ tels que $\phi^{st}_{G/conj}(x)$ soit l'image de $y$ dans $G_{reg}(F)/stconj$ (c'est l'ensemble vide si $y$ n'a pas d'image dans $G_{reg}(F)/stconj$). 

{\bf Remarque.} Pour \^etre pr\'ecis, il conviendrait d'introduire l'ensemble des \'el\'ements  $G$-r\'eguliers dans $H$ et de dire que l'\'egalit\'e doit \^etre v\'erifi\'ee pour tout $y$ $G$-r\'egulier au lieu de dire comme ci-dessus qu'elle l'est presque partout.

Soient $\Theta$ une distribution sur $G(F)$ invariante par conjugaison et $\Theta^H$ une distribution sur $H(F)$ invariante par conjugaison stable. On dit que $\Theta$ est un transfert de $\Theta^H$ si et seulement si $\Theta(f)=\Theta^H(f^H)$ pour tout couple de fonctions $(f,f^H)$ tel que $f^H$ soit un transfert de $f$. Supposons $\Theta$ et $\Theta^H$ localement int\'egrable, ce qui permet de les identifier \`a des fonctions d\'efinies presque partout. Les formules d'int\'egration du paragraphe 1.2 montrent que $\Theta$ est un transfert de $\Theta^H$ si et seulement si on a l'\'egalit\'e
$$(1) \qquad \Theta(x)D^G(x)^{1/2}=\sum_{y}D^H(y)^{1/2}\Delta_{H,G}(y,x)\Theta^H(y)$$
pour tout $x\in G_{reg}(F)/conj$, o\`u $y$ parcourt l'ensemble des \'el\'ements de $H_{reg}(F)/stconj$ tels que $\phi^{st}_{G/conj}(x)$ soit l'image de $y$ par la correspondance endoscopique.

  La th\'eorie de l'endoscopie s'adapte au cas d'un groupe tordu $(M,\tilde{M})$. On ne l'utilisera que sous des hypoth\`eses simplificatrices. On suppose $M$ d\'eploy\'e. On fixe un \'el\'ement $ \tilde{\theta}\in \tilde{M}(F)$,  on suppose qu'il existe une paire de Borel \'epingl\'ee de $M$, d\'efinie sur $F$ et conserv\'ee par $\theta_{\tilde{\theta}}$. On fixe une telle paire. L'\'el\'ement $N$ construit comme ci-dessus est un \'el\'ement nilpotent r\'egulier de $\mathfrak{m}_{\tilde{\theta}}(F)$ et c'est sa classe de conjugaison par $M_{\tilde{\theta}}(F)$ qui compte. Consid\'erons une donn\'ee endoscopique $(H,s,{^L\xi})$ de $(M,\tilde{M})$, cf. [KS] 2.1. On d\'efinit une correspondance endoscopique, qui est une application dont l'ensemble de d\'epart est un sous-ensemble de $H_{reg}(F)/stconj$ et l'ensemble d'arriv\'ee est $\tilde{M}(F)/stconj$. Il y a une difficult\'e avec les mesures. Soient $h\in H_{reg}(F)$ et $\tilde{m}\in \tilde{M}(F)$ tels que $\phi^{st}_{\tilde{M}}(\tilde{m})$ soit l'image de $\phi^{st}_{H}(h)$. Notons $T^{\flat}=M_{\tilde{m}}$, $T^H=H_{h}$ , $T$ le commutant de $T^{\flat}$ dans $M$ et $\theta=\theta_{\tilde{m}}$. Posons
$$d(\tilde{\gamma})=\vert det((1-\theta)_{\vert \mathfrak{t}/\mathfrak{t}^{\flat}}\vert _{F}$$
et
 $$D_{0}^{\tilde{M}}(\tilde{\gamma})=\vert det((1-\theta)_{\vert \mathfrak{m}/\mathfrak{t}}\vert _{F}=D^{\tilde{M}}(\tilde{\gamma})d(\tilde{\gamma})^{-1}.$$
Il y a une application naturelle de $T^{\flat}(F)$ dans $T^H(F)$ qui est localement un isomorphisme. On a d\'efini des mesures sur $T^{\flat}(F)$ et $T^H(F)$. D'apr\`es [KS] 5.5, il convient de remplacer la mesure sur $T^H(F)$ par celle telle que le jacobien de l'application pr\'ec\'edente soit \'egal \`a $d(\tilde{\gamma})$. On change de fa\c{c}on correspondante la mesure sur $H_{reg}(F)/stconj$. Alors  la correspondance endoscopique ci-dessus est de jacobien $d(\tilde{y})^{-1}$ en un point $y\in H_{reg}(F)/stconj$ d'image $\tilde{y} \in \tilde{M}_{reg}/stconj$.  On d\'efinit un facteur de transfert $\Delta_{H,\tilde{M}}$. Pour $\tilde{f}\in C_{c}^{\infty}(\tilde{M}(F))$ et $f^H\in C_{c}^{\infty}(H(F))$, on dit que $f^H$ est un transfert de $\tilde{f}$ si on a l'\'egalit\'e
$$D^H(y)^{1/2}\Phi^{st}(y,f^H)=\sum_{\tilde{x}}\Delta_{H,\tilde{M}}(y,\tilde{x})D_{0}^{\tilde{M}}(\tilde{x})^{1/2}\Phi(\tilde{x},\tilde{f})$$
pour presque tout $y\in H_{reg}(F)/stconj$, o\`u $\tilde{x}$ parcourt l'ensemble des \'el\'ements de $\tilde{M}_{reg}(F)/conj$ tels que $\phi^{st}_{\tilde{M}/conj}(\tilde{x})$ soit l'image de $y$ dans $\tilde{M}_{reg}(F)/stconj$ (cf. [KS] 5.5). On d\'efinit la notion de transfert de distributions comme ci-dessus. Soit $\tilde{\Theta}$ une distribution sur $\tilde{M}(F)$ invariante par conjugaison et $\Theta^H$ une distribution sur $H(F)$ invariante par conjugaison stable.  Supposons-les localement int\'egrables. Les normalisations de mesures ci-dessus et les formules d'int\'egration de 1.2 montrent que $\tilde{\Theta}$ est un transfert de $\Theta^H$ si et seulement si on a l'\'egalit\'e
$$(2) \qquad \tilde{\Theta}(\tilde{x})D_{0}^{\tilde{M}}(\tilde{x})^{1/2}=\sum_{y}D^H(y)^{1/2}\Delta_{H,\tilde{M}}(y,\tilde{x})\Theta^H(y)$$
pour tout $\tilde{x}\in \tilde{M}_{reg}(F)/conj$, o\`u $y$ parcourt l'ensemble des \'el\'ements de $H_{reg}(F)/stconj$ tels que $\phi^{st}_{\tilde{M}/conj}(\tilde{x})$ soit l'image de $y$ par la correspondance endoscopique.

On a modifi\'e ci-dessus nos mesures pour traduire les d\'efinitions de [KS]. Mais l'\'egalit\'e (2) ne fait intervenir aucune mesure. Dans la suite, on revient aux mesures d\'efinies en 1.2, en consid\'erant (2) comme la d\'efinition du transfert de la distribution $\Theta^H$.

 \bigskip
 
 \subsection{Endoscopie pour les groupes sp\'eciaux orthogonaux}
 
 Soient $(V',q')$, $(V'_{+},q'_{+})$ et $(V'_{-},q'_{-})$  trois espaces quadratiques sur $F$. Notons  $d'$, $d'_{+}$ et $d'_{-}$ leurs dimensions et $G'$, $G'_{+}$ et $G'_{-}$ leurs groupes sp\'eciaux orthogonaux. On suppose
 
 $\bullet$ $d'$, $d'_{+}$ et $d'_{-}$ impairs;
 
 $\bullet$ $d'_{+}+d'_{-}=d'+1$;
 
 $\bullet$ $G'_{+}$ et $G'_{-}$ sont d\'eploy\'es.

 Comme on le sait,  $G'_{+}\times G'_{-}$ est un groupe endoscopique de $G'$.

{\bf Remarque.} La notion correcte est celle de donn\'ee endoscopique et non de groupe endoscopique. Ici, les termes que l'on doit ajouter pour obtenir une donn\'ee endoscopique sont presque \'evidents, on renvoie \`a [W4] 1.8 pour plus de d\'etails.

Pour d\'efinir des facteurs de transfert, on doit fixer quelques donn\'ees suppl\'ementaires,  \`a savoir une forme  int\'erieure quasi-d\'eploy\'ee  $\underline{G}'$ de  $G'$, un torseur int\'erieur  $\psi_{G'}:G'\to \underline{G}'$ et un cocycle $u':Gal(\bar{F}/F)\to \underline{G}'$. Supposons $G'$ d\'eploy\'e.   On pose $(\underline{V}',\underline{q}')=(V',q')$,   $\underline{G}'=G'$,  on prend pour torseur int\'erieur l'identit\'e et pour cocycle le cocyle trivial. Supposons maintenant $G'$ non d\'eploy\'e. On introduit l'espace quadratique  $(\underline{V}',\underline{q}')$  qui a m\^emes dimension et discriminant que   $(V',q')$ mais un indice de Witt oppos\'e.  Le groupe sp\'ecial orthogonal   $\underline{G}'$ de cet espace est d\'eploy\'e.  On fixe un isomorphisme $\bar{F}$-lin\'eaire $\beta:V'\otimes_{F}\bar{F}\to \underline{V}'\otimes_{F}\bar{F}$ qui identifie les prolongements $\bar{F}$-bilin\'eaires de $q'$ et $\underline{q}'$.  Pour  $g'\in G'$, on pose $\psi_{G'}(g')=\beta g'\beta^{-1}$. Pour $\sigma\in Gal(\bar{F}/F)$, on pose $u'(\sigma)=\beta\sigma(\beta)^{-1}$.  On a a priori besoin de fixer  une classe de conjugaison d'\'el\'ement nilpotent r\'egulier dans   $\underline{\mathfrak{g}}'(F)$, mais il n'y en a qu'une.

 Avec  ces donn\'ees, on peut d\'efinir la notion de correspondance entre classes de conjugaison stable semi-simples r\'eguli\`eres dans $G'(F)$ et dans $G'_{+}(F)\times G'_{-}(F)$, et d\'efinir un facteur de transfert   $\Delta_{G'_{+}\times G'_{-},G'}$.  La correspondance entre classes de conjugaison stable s'explicite ais\'ement. Il y a une application naturelle
$$\begin{array}{ccc}\Xi_{d'_{+}-1}\times \Xi_{d'_{-}-1}&\to&\Xi_{d'-1}\\ (\xi_{+},\xi_{-})&\mapsto& \xi_{+}\sqcup \xi_{-}.\\ \end{array}$$
Via les applications $p_{G'_{+}}^{st}$, $p_{G'_{-}}^{st}$ et $p_{G'}^{st}$, c'est la correspondance endoscopique (encore une fois, pour \^etre correct, il faut se restreindre aux \'el\'ements $G'$-r\'eguliers).

Soient $y'_{+}\in G'_{+,reg}(F)/stconj$, $y'_{-}\in G'_{-,reg}(F)/stconj$ et $x'\in G'_{reg}(F)/conj$. Posons $\xi_{+}=p_{G'_{+}}^{st}(y_{+})$, $\xi_{-}=p_{G'_{-}}^{st}(y_{-})$, $\xi=\xi_{+}\sqcup\xi_{-}$, suposons $p_{G'}(x')=\xi$. Soit $c\in C(\xi)$ qui param\`etre $x'$. Alors, avec la notation introduite en 1.5,
$$(1) \qquad \Delta_{G'_{+}\times G'_{-},G'}((y'_{+},y'_{-}),x')=\Delta_{-2\delta(q')}(\xi_{+},\xi_{-},c).$$
 Cf. [W4] proposition 1.10. Le $\eta$ de cette r\'ef\'erence vaut ici $(-1)^{(d'-1)/2}\delta(q')$.
 
 {\bf Remarque.} Soit $\alpha\in F^{\times}$, rempla\c{c}ons $q'$ par $\alpha q'$. Cela ne change pas le groupe $G'$. Pour $y'_{+}$, $y'_{-}$ et $x'$ comme ci-dessus, le param\`etre $\xi$ est inchang\'e. En se reportant au paragraphe 1.5, on voit que l'\'el\'ement $c=(c_{i})_{i\in I^*}$ est remplac\'e par $(\alpha c_{i})_{i\in I^*}$. On a aussi $\delta(\alpha q')=\alpha\delta(q')$. On conclut que le facteur de transfert ne change pas.
 
  Soient $(V,q)$, $(V_{+},q_{+})$ et $(V_{-},q_{-})$  trois espaces quadratiques sur $F$. Notons  $d$, $d_{+}$ et $d_{-}$ leurs dimensions et $G$, $G_{+}$ et $G_{-}$ leurs groupes sp\'eciaux orthogonaux. On suppose
 
 $\bullet$ $d$, $d_{+}$ et $d_{-}$ pairs;
 
 $\bullet$ $d_{+}+d_{-}=d$;
 
 $\bullet$ $G_{+}$ et $G_{-}$ sont quasi-d\'eploy\'es.

Le groupe $G_{+}\times G_{-}$ est un groupe endoscopique de $G$. On doit fixer des donn\'ees suppl\'ementaires. Pour cela, 

{\bf on fixe pour la suite de l'article un \'el\'ement $\nu_{0}\in F^{\times}$.}

On introduit l'espace quadratique $(Z_{1},q_{1,-\nu_{0}})$ suivant: $Z_{1}=F$ et $q_{1,-\nu_{0}}(\lambda,\lambda')=-2\nu_{0}\lambda\lambda'$. Consid\'erons la condition:

{\bf (QD)  le groupe sp\'ecial orthogonal de $(V,q)\oplus (Z_{1},q_{1,-\nu_{0}})$ est d\'eploy\'e. }

Il est facile de l'expliciter. Si $\delta(q)=1$, elle est \'equivalente \`a chacune des conditions suivantes: $G$ est d\'eploy\'e; $(V,q)$ est somme orthogonale de plans hyperboliques. Supposons $\delta(q)\not=1$. Posons $E=F(\sqrt{\delta(q)})$. Alors $(V,q)$ est somme orthogonale de plans hyperboliques et de l'espace $E$ muni d'une forme $(\lambda\lambda')\mapsto trace_{E/F}(\tau_{E/F}(\lambda)\lambda'\eta)$, pour un \'el\'ement $\eta\in F^{\times}$. La condition (QD) \'equivaut alors \`a $sgn_{E/F}(\eta\nu_{0})=1$.

Rappelons que les orbites nilpotentes r\'eguli\`eres dans $\mathfrak{g}(F)$ sont param\'etr\'ees par un sous-ensemble ${\cal N}$ de $F^{\times}/F^{\times,2}$, cf. [W1] 7.1. Sous l'hypoth\`ese (QD), on a ${\cal N}=F^{\times}/F^{\times,2}$ si $\delta(q)=1$, tandis que ${\cal N}=\eta Norm_{E/F}(E^{\times})/F^{\times,2}$ si $\delta(q)\not=1$, avec les notations ci-dessus.

Si (QD) est v\'erifi\'ee, on pose $(\underline{V},\underline{q})=(V,q)$, $\underline{G}=G$ (ce groupe est quasi-d\'eploy\'e).  On prend pour torseur int\'erieur l'identit\'e et pour cocycle le cocyle trivial. Supposons maintenant  que (QD) ne soit pas v\'erifi\'ee. On introduit l'espace quadratique  $(\underline{V},\underline{q})$  qui a m\^emes dimension et discriminant que   $(V,q)$ mais un indice de Witt oppos\'e.  Cet espace v\'erifie (QD).  On fixe un isomorphisme $\bar{F}$-lin\'eaire $\beta:V\otimes_{F}\bar{F}\to \underline{V}\otimes_{F}\bar{F}$ qui identifie les prolongements $\bar{F}$-bilin\'eaires de $q$ et $\underline{q}$. Comme en dimension impaire, on en d\'eduit un torseur int\'erieur $\psi_{G}:G\to \underline{G}$ et un cocycle $u:Gal(\bar{F}/F)\to \underline{G}$.  On doit encore fixer une orbite nilpotente r\'eguli\`ere dans $\underline{\mathfrak{g}}(F)$.  D'apr\`es ce qui pr\'ec\`ede, l'\'el\'ement $\nu_{0}$ param\`etre une telle orbite ${\cal O}_{\nu_{0}}$ et c'est celle que l'on choisit.

Il y a une application naturelle
$$\begin{array}{ccc}\Xi_{d_{+},\delta(q_{+})}\times \Xi_{d_{-},\delta(q_{-})}&\to&\Xi_{d,\delta(q)}\\ (\xi_{+},\xi_{-})&\mapsto& \xi_{+}\sqcup \xi_{-}.\\ \end{array}$$
Via les applications $p_{G_{+}}^{st,+}$, $p_{G_{-}}^{st,+}$ et $p_{G}^{st,+}$,  elle d\'efinit  une correspondance entre classes de conjugaison stable au sens de la conjugaison par le groupe orthogonal tout entier, autrement dit entre couples de classes de conjugaison stable. La correspondance endoscopique raffine cette premi\`ere correspondance. Pour la d\'ecrire plus commod\'ement, introduisons l'ensemble $\Lambda(V)$ form\'e des classes de conjugaison par $G(\bar{F})$ dans l'ensemble des lagrangiens d\'efinis sur $\bar{F}$ de $V\otimes_{F}\bar{F}$. Il a deux \'el\'ements. Le groupe orthogonal $G^+(\bar{F})$ agit naturellement sur cet ensemble, un \'el\'ement de d\'eterminant $-1$ \'echangeant les deux \'el\'ements. Consid\'erons le produit $G_{reg}(F)\times \Lambda(V)$. Disons que deux \'el\'ements $(g_{1},L_{1})$ et $(g_{2},L_{2})$ de ce produit sont stablement conjugu\'es si et seulement s'il existe $g\in G^+(\bar{F})$ tel que $g_{2}=gg_{1}g^{-1}$ et $L_{2}=g(L_{1})$. Notons $\Lambda G_{reg}(F)/stconj^+$ l'ensemble des classes de conjugaison stable et $\phi_{\Lambda G}^{st,+}:G_{reg}(F)\times \Lambda(V)\to \Lambda G_{reg}(F)/stconj^+$ l'application naturelle.  Il y a une application naturelle de cet ensemble dans $G_{reg}(F)/stconj^+$ qui, \`a $\phi_{\Lambda G}^{st,+}(g,L)$, associe la classe  $\phi_{G}^{st,+}(g)$ de $g$ dans   l'ensemble d'arriv\'ee. D'autre part

(2)  si l'on fixe $L\in \Lambda(V)$, l'application $g\mapsto(g,L)$ se quotiente en une bijection de $G_{reg}(F)/stconj$ sur $\Lambda G_{reg}(F)/stconj^+$. 

Soit $(g,L)\in G_{reg}(F)\times \Lambda(V)$. On peut repr\'esenter $L$ par un lagrangien, encore not\'e $L$, tel que $g(L)=L$. Notons ${\cal E}(g,L)$ l'ensemble des valeurs propres de $g$ dans $L$ (ce sont des \'el\'ements de $\bar{F}$). Quand on change de repr\'esentant $L$, l'ensemble ${\cal E}(g,L)$ peut changer. Un  tel changement consiste \`a remplacer un nombre pair de valeurs propres par leurs inverses. Disons que deux ensembles de valeurs propres sont \'equivalents s'ils diff\`erent par un tel changement. La classe d'\'equivalence de ${\cal E}(g,L)$ est alors bien d\'etermin\'ee et ne d\'epend que de  $\phi_{\Lambda G}^{st,+}(g,L)$. On remarque que, si $L'$ d\'esigne l'autre \'el\'ement de $\Lambda(V)$, les classes ${\cal E}(g,L')$ et ${\cal E}(g,L)$ sont distinctes: on passe de l'une \`a l'autre en rempla\c{c}ant un nombre impair de valeurs propres par leurs inverses. Cela \'etant, l'application pr\'ec\'edemment d\'efinie entre couples de conjugaison stable se rel\`eve en une application
$$E: \Lambda G_{+,reg}(F)/stconj^+\times \Lambda G_{-,reg}(F)/stconj^+\to \Lambda G_{reg}(F)/stconj^+$$
d\'efinie ainsi.  Soit $(g_{+},\Lambda_{+}) \in G_{+,reg}(F)\times \Lambda(V_{+})$ et $(g_{-},\Lambda_{-})\in G_{-,reg}(F)\times \Lambda(V_{-})$. Soit $g\in G_{reg}$ tel que $\phi_{G}^{st,+}(g)$ corresponde \`a $(\phi_{G_{+}}^{st,+}(g_{+}),\phi_{G_{-}}^{st,+}(g_{-}))$. Il y a un unique $L\in \Lambda(V)$ tel que ${\cal E}(g,L)$ soit \'equivalent \`a la r\'eunion disjointe  ${\cal E}(g_{+},L_{+})\sqcup{\cal E}(g_{-},L_{-})$. On choisit $L$ ainsi. On d\'efinit   $E(\phi_{\Lambda G_{+}}^{st,+}(g_{+},L_{+}),\phi_{\Lambda G_{-}}^{st,+}(g_{-},L_{-})) =\phi_{\Lambda G}^{st,+}(g,L)$.

 Comme on l'a dit, une donn\'ee  endoscopique  est plus riche que la seule donn\'ee des groupes. Elle contient une identification sur $\bar{F}$ entre un tore maximal de $G_{+}\times G_{-}$ et un tore maximal de $G$, modulo l'action du groupe de Weyl de $G$. On a vu  ci-dessus qu'un \'el\'ement de $\Lambda(V)$ d\'eterminait une classe d'\'equivalence d'ensembles de valeurs propres. Il d\'etermine de m\^eme une classe d'\'equivalence d'ensembles de caract\`eres d'un tore. De l'identification des tores se d\'eduit une application $\Lambda(V_{+})\times \Lambda(V_{-})\to \Lambda(V)$. Soit alors $(\Lambda_{+},\Lambda_{-})\in \Lambda(V_{+})\times \Lambda(V_{-})$, notons $\Lambda$ son image par cette application. En utilisant (2) appliqu\'e respectivement \`a $\Lambda_{+}$, $\Lambda_{-}$ et $\Lambda$, l'application $E$ devient une application
 $$G_{+,reg}(F)/stconj\times G_{-,reg}(F)/stconj\to G_{reg}(F)/stconj.$$
 On voit qu'elle ne d\'epend pas du choix de $(\Lambda_{+},\Lambda_{-})$. C'est la correspondance endoscopique cherch\'ee.

   Soient $y_{+}\in G_{+,reg}(F)/stconj$, $y_{-}\in G_{-,reg}(F)/stconj$ et $x\in G_{reg}(F)/conj$. Posons $\xi_{+}=p_{G_{+}}^{st}(y_{+})$, $\xi_{-}=p_{G_{-}}^{st}(y_{-})=\xi_{-}$, $\xi=\xi_{+}\sqcup \xi_{-}$, supposons que   la classe de conjugaison stable de $x$ corresponde \`a $(y_{+},y_{-})$, a fortiori $p_{G}(x)=\xi$. Soit $c\in C(\xi)$ qui param\`etre $x$. Alors
 $$(3) \qquad \Delta_{G_{+}\times G_{-},G}((y_{+},y_{-}),x)= \Delta_{\nu_{0}\delta(q)}(\xi_{+},\xi_{-},c).$$
Cf. [W4] proposition 1.10. Le $\eta$ de cette r\'ef\'erence vaut ici $2(-1)^{d/2}\nu_{0}\delta(q)$.

{\bf Remarque.} Soit $\alpha\in F^{\times}$. Rempla\c{c}ons $q$ par $\alpha q$. Le groupe $G$ reste inchang\'e. Si l'on conservait le m\^eme cocycle et si l'on rempla\c{c}ait $\nu_{0}$ par $\alpha\nu_{0}$, on conserverait le m\^eme facteur de transfert. Mais, avec nos d\'efinitions, le cocycle peut changer et, par contre, on veut conserver le m\^eme $\nu_{0}$. Donc le facteur de transfert n'est plus le m\^eme. Pr\'ecis\'ement, on voit qu'il est multipli\'e par $(\delta(q_{-}),\alpha)_{F}$.

 \bigskip
 
 \subsection{Endoscopie pour les groupes lin\'eaires tordus}

  On consid\`ere un espace $U$  sur $F$ de dimension $n$  et on introduit le groupe tordu $(M,\tilde{M})$  correspondant, cf. 1.4. Introduisons d'autre part un espace quadratique $(V,q)$, dont on note $G$ le groupe sp\'ecial orthogonal. On suppose $dim(V)=n$ et  $G$   quasi-d\'eploy\'e. Le groupe $G$ est un groupe endoscopique de $\tilde{M}$. Il convient ici de pr\'eciser les plongements de $L$-groupes.  On renvoie pour cela \`a [W4] 1.8. Si $n$ est pair, c'est dans cette r\'ef\'erence  le cas du groupe lin\'eaire tordu avec $d=n$ pair, $d^-=d$, $d^+=0$. Si $n$ est impair, c'est le cas du groupe lin\'eaire tordu avec $d=n$ impair, $d^-=d$, $d^+=0$ et $\chi=1$. Pour d\'efinir des facteurs de transfert, on doit fixer un \'el\'ement de base de $\tilde{M}(F)$  (l'\'el\'ement not\'e $\tilde{\theta}$ en 1.6). On fixe une base $(u_{j})_{j=1,...,n}$ de $U$ et on prend  pour point-base l'\'el\'ement $\boldsymbol{\theta}_{n}\in \tilde{M}(F)$ d\'efini par
  $$\boldsymbol{\theta}_{n}(u_{j},u_{k})=(-1)^{k+[(n+1)/2]}\delta_{j,n+1-k},$$
  o\`u le dernier terme est le symbole de Kronecker. On doit aussi choisir une classe de conjugaison d'\'el\'ement nilpotent r\'egulier dans $\mathfrak{m}_{\boldsymbol{\theta}_{n}}(F)$. Si $n$ est impair, le groupe $M_{\boldsymbol{\theta}_{n}}$ est sp\'ecial orthogonal "impair", il n'y a qu'une seule classe possible. Si $n$ est pair, le groupe $M_{\boldsymbol{\theta}_{n}}$ est symplectique, il y a le choix. On prend la classe de l'\'el\'ement $N$ tel que $N(u_{j})=-u_{j-1}$ pour $j=2,...,n$ et $N(u_{1})=0$. On peut alors d\'efinir une correspondance entre classes de conjugaison stable d'\'el\'ements semi-simples r\'eguliers dans $G(F)$ et dans $\tilde{M}(F)$, ainsi qu'un facteur de transfert $\Delta_{G,\tilde{M}}$. La correspondance entre classes de conjugaison stable s'explicite ais\'ement. Supposons $n$ impair. On a les bijections
 $$G_{reg}(F)/stconj\stackrel{p_{G}^{st}}{\to}\Xi_{reg,n-1}\stackrel{p_{\tilde{M}}^{st}}{\leftarrow}\tilde{M}_{reg}(F)/stconj$$
 et la correspondance est celle qui vient de ces bijections. Supposons maintenant $n$ pair.  On d\'efinit une bijection $\xi\mapsto -\xi$ de $\Xi_{reg}$ dans lui-m\^eme: si $\xi=(I,(F_{\pm i})_{i\in I},(F_{i})_{i\in I},(y_{i})_{i\in I})$, on pose $-\xi=(I,(F_{\pm i})_{i\in I},(F_{i})_{i\in I},(-y_{i})_{i\in I})$. On a les applications
  $$G_{reg}(F)/stconj\stackrel{p_{G}^{st}}{\to}\Xi_{reg,n,\delta(q)}\stackrel{\xi\mapsto -\xi}{\to}\Xi_{reg,n}\stackrel{p_{\tilde{M}}^{st}}{\leftarrow}\tilde{M}_{reg}(F)/stconj$$
  et la correspondance est celle qui vient de ces applications.

 Soient $y\in G_{reg}(F)/stconj$ et $\tilde{x}\in \tilde{M}_{reg}(F)/conj$. Posons $\xi=p_{\tilde{M}}(\tilde{x})$. Supposons d'abord $n$ impair et $p^{st}_{G}(y)=\xi$. Soit $\gamma\in \Gamma_{imp}(\xi)$ qui param\`etre la classe de conjugaison de $\tilde{x}$. Alors, avec la notation de 1.5,
 $$(1) \qquad \Delta_{G,\tilde{M}}(y,\tilde{x})=\Delta(\xi,\gamma).$$
 Supposons maintenant $n$ pair et $p^{st}_{G}(y)=-\xi$. Soit $\gamma\in \Gamma_{pair}(\xi)$ qui param\`etre la classe de conjugaison de $\tilde{x}$. Alors
 $$(2) \qquad \Delta_{G,\tilde{M}}(y,\tilde{x})=\Delta(\xi,\gamma).$$
    Cf. [W4] proposition 1.10.   Dans le cas $n$ pair, le $\eta$ de cette r\'ef\'erence vaut $(-1)^{n/2}$.  On doit remplacer les $\xi $ et $y_{i}$ par $-\xi$ et $-y_{i}$ car, dans [W4], ces termes param\'etraient  $y$ et non pas la classe de conjugaison stable de $\tilde{x}$.

\bigskip

\section{Quelques formules revisit\'ees}

\bigskip

\subsection{Multiplicit\'es pour les groupes sp\'eciaux orthogonaux}

Soient $(V,q)$ et $(V',q')$ deux espaces quadratiques sur $F$ de dimensions respectives $d$ et $d'$. On note $G$ et $G'$ leurs groupes sp\'eciaux orthogonaux. On suppose $d$ pair et $d'$ impair.   Pour un entier $r\in {\mathbb N}$ introduisons un espace  $Z_{2r+1}$ de dimension $2r+1$ sur $F$, muni d'une base $(z_{j})_{j=-r,...,r}$. Pour $\nu\in F^{\times}$, d\'efinissons une forme bilin\'eaire sym\'etrique $q_{2r+1,\nu}$ sur $Z_{2r+1}$ par
$$q_{2r+1,\nu}(\sum_{i=-r,...,r}\lambda_{j}z_{j},\sum_{j=-r,...,r}\lambda'_{j}z_{j})=2\nu \lambda_{0}\lambda'_{0}+\sum_{j=1,...,r}(\lambda_{-j}\lambda'_{j}+\lambda_{j}\lambda'_{-j}).$$
On a fix\'e $\nu_{0}\in F^{\times}$ en 1.7.  Supposons

$\bullet$ si $d<d'$, $(V',q')$ est isomorphe \`a $(V,q)\oplus (Z_{d'-d},q_{d'-d,-\nu_{0}})$;

$\bullet$ si $d>d'$, $(V,q)$ est isomorphe \`a $(V',q')\oplus (Z_{d-d'},q_{d-d',\nu_{0}})$.

On fixe de tels isomorphismes. Pour fixer les id\'ees, supposons $d'>d$, les constructions \'etant similaires dans l'autre cas. Alors $G$ s'identifie \`a un sous-groupe de $G'$. On pose $r=(d'-d-1)/2$. On note $P$ le sous-groupe parabolique de $G'$   qui conserve le drapeau de sous-espaces
$$Fz_{r}\subset Fz_{r}\oplus Fz_{r-1}\subset...\subset Fz_{r}\oplus...\oplus Fz_{1}.$$
On note $U$ son radical unipotent et on d\'efinit un caract\`ere $\psi_{U}$ de $U(F)$ par la formule
$$\psi_{U}(u)=\psi_{F}(\sum_{j=0,...,r-1}q(uz_{j},z_{-j-1})).$$
Le groupe $G$ est inclus dans $P$ et $\psi_{U}$ est invariant par conjugaison par $G(F)$. 

Soient $\sigma\in Temp(G(F))$ et $\sigma'\in Temp(G'(F))$. On note $Hom_{G,\psi_{F}}(\sigma',\sigma)$ l'espace des applications lin\'eaires $\varphi:E_{\sigma'}\to E_{\sigma}$ telles que
$$\varphi(\sigma'(ug)e)=\psi_{U}(u)\sigma(g)\varphi(e)$$
pour tous $u\in U(F)$,  $g\in G(F)$ et $e\in E_{\sigma'}$. D'apr\`es [AGRS] (compl\'et\'e par [W5]) et [GGP] corollaire 15.2, cet espace est de dimension au plus $1$. On pose
$$m(\sigma,\sigma')=m(\sigma',\sigma)=dim(Hom_{G,\psi_{F}}(\sigma',\sigma)).$$
Ce nombre est ind\'ependant de $\psi_{F}$.

Dans [W2], on a calcul\'e cette dimension par une formule int\'egrale. Rappelons-la. Notons $\underline{{\cal T}}$ l'ensemble des sous-tores $T$ de $G$ (on suppose encore $d'>d$) pour lesquels il existe une d\'ecomposition orthogonale $V=W_{T}\oplus V_{T}$ de sorte que les conditions (1) \`a (4) ci-dessous soient v\'erifi\'ees. On   pose $V'_{T}=V_{T}\oplus Z_{2r+1}$ et on note $G_{T}$, $G_{V_{T}}$ et $G_{V'_{T}}$ les groupes sp\'eciaux orthogonaux de $W_{T}$, $V_{T}$ et $V'_{T}$.

(1) $T$ est anisotrope.

(2) $dim(W_{T})$ est paire.

(3) $T$ est inclus dans $G_{T}$ et c'en est un sous-tore maximal.

(4) Les groupes  $G_{V_{T}}$ et $G_{V'_{T}}$ sont quasi-d\'eploy\'es.

Pour $T\in \underline{{\cal T}}$, on pose
$$W(G,T)=Norm_{G(F)}(T)/Z_{G(F)}(T).$$
Pour $t\in T(F)$, on pose $\Delta(t)=\vert det((1-t)_{\vert W_{T}})\vert _{F}$. Supposons $t$ en position g\'en\'erale. Alors $G'_{t}=T\times G_{V'_{T}}$ et $G_{t}=T\times G_{V_{T}}$. La repr\'esentation $\sigma'$ admet un caract\`ere $\Theta_{\sigma'}$.  Rappelons dans ce contexte un r\'esultat de Harish-Chandra. Pour ${\cal O}\in Nil(\mathfrak{g}'_{t}(F))$, on d\'efinit l'int\'egrale orbitale $J_{{\cal O}}$ sur $C_{c}^{\infty}(\mathfrak{g}'_{t}(F))$ et sa transform\'ee de Fourier. Celle-ci est une distribution localement int\'egrable, donc s'identifie \`a une fonction $\hat{J}_{{\cal O}}$. Alors il existe une famille $(c_{\sigma',{\cal O}}(t))_{{\cal O}\in Nil(\mathfrak{g}'_{t}(F))}$ de nombres complexes de sorte que, pour $X$ r\'egulier dans un voisinage assez petit de $0$ dans $\mathfrak{g}'_{t} (F)$,  on ait l'\'egalit\'e
$$(5) \qquad \Theta_{\sigma'}(texp(X))=\sum_{{\cal O}\in Nil(\mathfrak{g}'_{t}(F))}c_{\sigma',{\cal O}}(t)\hat{J}_{{\cal O}}(X).$$
 Evidemment, pour que ceci ait un sens pr\'ecis, on doit normaliser les mesures et la transformation de Fourier utilis\'ee, on renvoie pour cela \`a [W1] 1.2 et 2.6. On a $Nil(\mathfrak{g}'_{t}(F))=Nil(\mathfrak{g}_{V'_{T}}(F))$. Puisque le groupe $G_{V'_{T}}$ est quasi-d\'eploy\'e et que $dim(V'_{T})$ est impaire, il y a une unique orbite nilpotente r\'eguli\`ere. Notons-la ${\cal O}_{reg}$. On pose $c_{\sigma'}(t)=c_{\sigma',{\cal O}_{reg}}(t)$. Cela d\'efinit une fonction $c_{\sigma'}$ presque partout sur $T(F)$. Le m\^eme proc\'ed\'e s'applique pour d\'efinir une fonction $c_{\sigma}$, \`a ceci pr\`es que, parce que $dim(V_{T})$ est paire, il peut y avoir plusieurs orbites r\'eguli\`eres dans $Nil(\mathfrak{g}_{V_{T}}(F))$. Pr\'ecisons ce point. Si $dim(V_{T})\leq 2$, il n'y a encore qu'une orbite r\'eguli\`ere ${\cal O}_{reg}$ et on d\'efinit $c_{\sigma}(t)=c_{\sigma,{\cal O}_{reg}}(t)$. Supposons $dim(V_{T})\geq 4$. L'hypoth\`ese que $G_{V'_{T}}$ est quasi-d\'eploy\'e  implique qu'il existe un \'el\'ement $v_{0}\in V_{T}$ tel que $q(v_{0},v_{0})=2\nu_{0}$. Fixons un tel \'el\'ement et notons $G'_{0}$ le groupe sp\'ecial orthogonal de l'hyperplan orthogonal \`a $Fv_{0}$. La m\^eme hypoth\`ese implique que $G'_{0}$ est d\'eploy\'e. Il y a donc une unique orbite nilpotente r\'eguli\`ere dans $\mathfrak{g}'_{0}(F)$. Fixons un point de cet orbite. Alors sa $G_{V_{T}}(F)$-orbite est encore r\'eguli\`ere dans $\mathfrak{g}_{V_{T}}(F)$, on la note ${\cal O}_{\nu_{0}}$ (c'est l'orbite param\'etr\'ee par $\nu_{0}$ dans le param\'etrage \'evoqu\'e en 1.7). On pose $c_{\sigma}(t)=c_{\sigma,{\cal O}_{\nu_{0}}}(t)$.
 
 Le groupe $G(F)$ agit par conjugaison sur $\underline{{\cal T}}$. Fixons un ensemble de repr\'esentants ${\cal T}$ des orbites. Posons
 $$m_{geom}(\sigma,\sigma')=m_{geom}(\sigma',\sigma)=\sum_{T\in {\cal T}}\vert W(G,T)\vert ^{-1}\int_{T(F)}c_{\sigma}(t)c_{\sigma'}(t)D^{G}(t)\Delta(t)^rdt.$$
 Cette expression est absolument convergente et le th\'eor\`eme principal de [W2] affirme que l'on a l'\'egalit\'e
 $$m_{geom}(\sigma,\sigma')=m(\sigma,\sigma').$$
 
 {\bf Remarque.} Dans les d\'efinitions de [W1], c'est l'orbite ${\cal O}_{-\nu_{0}}$ et non ${\cal O}_{\nu_{0}}$  qui intervient. Ce n'est plus le cas ici car on a gliss\'e le signe $-$ dans les hypoth\`eses concernant le plongement de $(V,q)$ dans $(V',q')$.
 
 \bigskip
 
 \subsection{Reformulation de $m_{geom}(\sigma,\sigma')$}
 
 Notons ${\cal D}(q,q')$ l'ensemble des couples $({\bf d},\boldsymbol{\delta})\in {\mathbb N}\times F^{\times}/F^{\times,2}$ qui v\'erifient les conditions suivantes:
 
 (1) ${\bf d}$ est pair et $0\leq {\bf d}\leq inf(d,d')$;
 
 (2) si $G$ ou $G'$ n'est pas quasi-d\'eploy\'e, ${\bf d}\geq 2$; si $G$ et $G'$ sont quasi-d\'eploy\'es et ${\bf d}=0$, alors $\boldsymbol{\delta}=1$;
 
 (3) si ${\bf d}=inf(d,d')$ (ce qui entra\^{\i}ne par parit\'e ${\bf d}=d<d'$), $\boldsymbol{\delta}=\delta(q)$;
 
 (4) si ${\bf d}=2$, $\boldsymbol{\delta}\not=1$.
 
 Pour $({\bf d},\boldsymbol{\delta})\in {\cal D}(q,q')$, consid\'erons les espaces quadratiques $(W,q_{W})$ tels que $dim(W)={\bf d}$ et $\delta(q_{W})=\boldsymbol{\delta}$. Si ${\bf d}=0$, il n'y en a qu'un, \'evidemment. Si ${\bf d}\not=0$, il y en a deux, \`a \'equivalence pr\`es, que l'on distingue selon leurs indices de Witt. D'autre part, pour un espace quadratique $(W,q_{W})$ consid\'erons la condition suivante:
 
 $(H)$ si $d<d'$, il existe un espace quadratique $(V_{1},q_{1})$ tel que $(V,q)$ soit isomorphe \`a $(W,q_{W})\oplus (V_{1},q_{1})$ et que les groupes sp\'eciaux orthogonaux de $(V_{1},q_{1})$ et de $(V_{1},q_{1})\oplus (Z_{2r+1},q_{2r+1,-\nu_{0}})$ soient quasi-d\'eploy\'es; si $d>d'$, il existe un espace quadratique $(V_{1},q_{1})$ tel que $(V',q')$ soit isomorphe \`a $(W,q_{W})\oplus (V_{1},q_{1})$ et que les groupes sp\'eciaux orthogonaux de $(V_{1},q_{1})$ et de $(V_{1},q_{1})\oplus (Z_{2r+1},q_{2r+1,\nu_{0}})$ soient quasi-d\'eploy\'es.
 
 Montrons que
 
 (5) pour tout $({\bf d},\boldsymbol{\delta})\in {\cal D}(q,q')$, il existe \`a \'equivalence pr\`es un unique espace quadratique $(W,q_{W})$ qui v\'erifie la condition $(H)$ ainsi que les \'egalit\'es $dim(W)={\bf d}$, $\delta(q_{W})=\boldsymbol{\delta}$.
 
 On suppose pour fixer les id\'ees $d<d'$, le cas oppos\'e \'etant similaire. Si ${\bf d}=0$, l'espace nul convient d'apr\`es la condition (2). Si ${\bf d}=d>0$, un seul espace convient, \`a savoir $(W,q_{W})=(V,q)$. Supposons $0<{\bf d}<d$. Montrons d'abord que la condition $(H)$ \'equivaut \`a
 
 $(H')$ il existe un espace quadratique $(V'_{1},q'_{1})$ tel que son groupe sp\'ecial orthogonal soit d\'eploy\'e et que $(V',q')$ soit isomorphe \`a $(W,q_{W})\oplus (V'_{1},q'_{1})$.
  
 Rappelons que pour un espace quadratique $(V'_{1},q'_{1})$ de dimension impaire, les trois propri\'et\'es suivantes sont \'equivalentes: son groupe sp\'ecial orthogonal est quasi-d\'eploy\'e; il est d\'eploy\'e; l'espace $(V'_{1},q'_{1})$ est somme orthogonale d'une droite et de plans hyperboliques. Si $(H)$ est v\'erifi\'e, $(H')$ l'est aussi car l'espace $(V'_{1},q'_{1})=(V_{1},q_{1})\oplus (Z_{2r+1},q_{2r+1,-\nu_{0}})$  convient. Inversement, supposons $(H')$ v\'erifi\'ee. Puisque ${\bf d}<d$, l'espace $(V'_{1},q'_{1})$ peut se d\'ecomposer en somme orthogonale d'un espace de dimension $d-{\bf d}-1$ et de $r+1$ plans
hyperboliques. A fortiori, il existe un espace quadratique $(V_{1},q_{1})$ de dimension $d-{\bf d}$ tel que $(V'_{1},q'_{1})=(V_{1},q_{1})\oplus (Z_{2r+1},q_{2r+1,-\nu_{0}})$. Cette \'egalit\'e entra\^{\i}ne que $(V_{1},q_{1})$ est somme orthogonale de plans dont au plus un n'est pas hyperbolique. Donc le groupe sp\'ecial orthogonal de $(V_{1},q_{1})$ est quasi-d\'eploy\'e. D'autre part, d'apr\`es le th\'eor\`eme de Witt, l'\'egalit\'e $(V',q')=(W,q_{W})\oplus (V'_{1},q'_{1})$ entra\^{\i}ne $(V,q)=(W,q_{W})\oplus (V_{1},q_{1})$. Cela v\'erifie $(H)$.

Cela \'etant, soit $(W,q_{W})$ ayant les dimension et discriminant requis. A \'equivalence pr\`es, il existe un unique espace quadratique $(V'_{1},q'_{1})$ tel que son groupe sp\'ecial orthogonal soit d\'eploy\'e et  que l'espace $(V'_{2},q'_{2})=(W,q_{W})\oplus (V'_{1},q'_{1})$ ait le m\^eme discriminant que $(V',q')$:  c'est la somme de plans hyperboliques et d'une droite sur laquelle la forme est uniquement d\'etermin\'ee par cette condition. Le fait que $(V'_{2},q'_{2})$ soit isomorphe \`a $(V',q')$ se lit sur les indices de Witt. Quand on remplace $(W,q_{W})$ par l'espace qui a m\^eme dimension et discriminant, mais l'indice de Witt oppos\'e, cela ne change pas $(V'_{1},q'_{1})$, mais cela change l'indice de Witt de $(V'_{2},q'_{2})$. Il y a donc un et un seul choix de $(W,q_{W})$ pour lequel l'indice de Witt de $(V'_{2},q'_{2})$ est \'egal \`a celui de $(V',q')$. Cela prouve (5).

Pour tout $({\bf d},\boldsymbol{\delta})\in {\cal D}(q,q')$, on note $(W_{{\bf d},\boldsymbol{\delta}},q_{{\bf d},\boldsymbol{\delta}})$ l'espace quadratique dont (1) affirme l'existence, que l'on r\'ealise comme un sous-espace de $(V',q') $ ou $(V,q)$ conform\'ement \`a la propri\'et\'e $(H)$.

 A tout $T\in \underline{{\cal T}}$ est associ\'e un espace quadratique $(W_{T},q_{T})$, o\`u $q_{T}$ est la restriction de $q$ \`a $W_{T}$. Montrons que:
 
 (6) l'ensemble des classes d'\'equivalence d'espaces $(W_{T},q_{T})$ quand $T$ d\'ecrit $\underline{\cal T}$ est \'egal \`a celui des classes d'\'equivalence d'espaces $(W_{{\bf d},\boldsymbol{\delta}},q_{{\bf d},\boldsymbol{\delta}})$ quand $({\bf d},\boldsymbol{\delta})$ d\'ecrit ${\cal D}(q,q')$. 
 
 Le second ensemble est inclus dans le premier: le groupe sp\'ecial orthogonal d'un espace $(W_{{\bf d},\boldsymbol{\delta}},q_{{\bf d},\boldsymbol{\delta}})$  contient au moins un sous-tore maximal anisotrope (gr\^ace \`a (4)) et un tel tore appartient \`a $\underline{\cal T}$. Montrons que le premier ensemble est inclus dans le second. D'apr\`es 2.1(4), tout espace $(W_{T},q_{T})$ v\'erifie $(H)$. Il suffit donc de prouver que le couple form\'e de la dimension de $W$ et du discriminant de $q_{T}$ appartient \`a ${\cal D}(q,q')$, ce qui est imm\'ediat. Cela prouve (6).
 
D'apr\`es (6), on peut supposer que, pour tout $T\in {\cal T}$, l'espace $(W_{T},q_{T})$ est l'un des espaces $(W_{{\bf d},\boldsymbol{\delta}},q_{{\bf d},\boldsymbol{\delta}})$. Cela d\'ecompose ${\cal T}$ en union disjointe de sous-ensembles ${\cal T}_{{\bf d},\boldsymbol{\delta}}$. Fixons maintenant $({\bf d},\boldsymbol{\delta})\in {\cal D}(q,q')$ et consid\'erons la contribution de l'ensemble ${\cal T}_{{\bf d},\boldsymbol{\delta}}$ \`a $m_{geom}(\sigma,\sigma')$. Supposons d'abord $0<{\bf d}<inf(d,d')$ et, pour fixer les id\'ees, supposons $d<d'$. Notons $G^+_{{\bf d},\boldsymbol{\delta}}$ le groupe orthogonal de $(W_{{\bf d},\boldsymbol{\delta}},q_{{\bf d},\boldsymbol{\delta}})$ et $G_{{\bf d},\boldsymbol{\delta}}$ le groupe sp\'ecial orthogonal. L'ensemble ${\cal T}_{{\bf d},\boldsymbol{\delta}}$ est un ensemble de repr\'esentants des  sous-tores maximaux et anisotropes de $G_{{\bf d},\boldsymbol{\delta}}  $ pour la relation d'\'equivalence suivante: deux tores sont \'equivalents s'ils sont conjugu\'es par un \'el\'ement de $G(F)$. On v\'erifie que cette relation \'equivaut \`a la conjugaison par un \'el\'ement de $G^+_{{\bf d},\boldsymbol{\delta}}  (F)$. De m\^eme, pour un tel tore $T$, on voit que le nombre d'\'el\'ements de $W(G,T)$ est le m\^eme que celui de $W(G^+_{{\bf d},\boldsymbol{\delta}}  ,T)$, cet ensemble \'etant d\'efini de fa\c{c}on \'evidente. Pour une fonction $f$ sur $G_{{\bf d},\boldsymbol{\delta}}  (F)_{ani}$, invariante par conjugaison par $G^+_{{\bf d},\boldsymbol{\delta}}  (F)$, on a alors l'\'egalit\'e
$$\sum_{T\in {\cal T}_{{\bf d},\boldsymbol{\delta}}  }\vert W(G,T)\vert ^{-1}\int_{T(F)}f(t)dt=\int_{G_{{\bf d},\boldsymbol{\delta}}  (F)_{ani}/conj^+}f(x)dx.$$
On a d\'efini les fonctions $c_{\sigma}$ et $c_{\sigma'}$  sur les \'el\'ements de ${\cal T}_{{\bf d},\boldsymbol{\delta}}  $, mais la m\^eme d\'efinition permet de les d\'efinir sur l'ensemble $G_{{\bf d},\boldsymbol{\delta}}  (F)_{ani}$ tout entier. Elles sont invariantes par l'action de $G^+_{{\bf d},\boldsymbol{\delta}}  (F)$: cela r\'esulte du fait que la conjugaison par un \'el\'ement de $G^+_{{\bf d},\boldsymbol{\delta}}(F)$ peut se r\'ealiser par celle d'un \'el\'ement de $G(F)$ et du fait que toute orbite nilpotente r\'eguli\`ere de l'alg\`ebre de Lie d'un groupe sp\'ecial orthogonal est invariante par l'action du groupe orthogonal tout entier. Il en est de m\^eme des fonctions $D^G$ et $\Delta$. La contribution de ${\cal T}_{{\bf d},\boldsymbol{\delta}}  $ \`a $m_{geom}(\sigma,\sigma')$ peut donc s'\'ecrire
$$(7) \qquad \int_{G_{{\bf d},\boldsymbol{\delta}}  (F)_{ani}/conj^+}c_{\sigma}(x)c_{\sigma'}(x)D^G(x)\Delta(x)^rdx.$$

Si ${\bf d}=0$, la contribution de ${\cal T}_{{\bf d},\boldsymbol{\delta}}$ a trivialement la m\^eme forme. Consid\'erons maintenant le cas ${\bf d}=d<d'$. Dans ce cas, on voit que ce n'est plus la conjugaison par $G^+_{{\bf d},\boldsymbol{\delta}}(F)$ qui intervient mais seulement celle par $G_{{\bf d},\boldsymbol{\delta}}(F)$. La fonction $c_{\sigma}$ n'a d'ailleurs plus de raison d'\^etre invariante par l'action du groupe orthogonal. On obtient dans ce cas 
 $$\int_{G_{{\bf d},\boldsymbol{\delta}}  (F)_{ani}/conj}c_{\sigma}(x)c_{\sigma'}(x)D^G(x)\Delta(x)^rdx.$$
 Mais, si l'on d\'efinit une fonction $c_{\sigma}$ sur $G_{{\bf d},\boldsymbol{\delta}} (F)_{ani}/conj^+$ par $c_{\sigma}(x)=c_{\sigma}(x')+c_{\sigma}(x'')$, o\`u $x'$ et $x''$ sont les deux \'el\'ements de $G_{{\bf d},\boldsymbol{\delta}} (F)_{ani}/conj$ d'image $x$, on retrouve la formule (7).
 
 Posons 
 $${\cal C}(q,q')=\sqcup_{({\bf d},\boldsymbol{\delta})\in {\cal D}(q,q')}G_{{\bf d},\boldsymbol{\delta}} (F)_{ani}/conj^+$$
 et
 $$\Xi^*(q,q')=\sqcup_{({\bf d},\boldsymbol{\delta})\in {\cal D}(q,q')}\Xi^*_{reg,{\bf d},\boldsymbol{\delta}}.$$
 D'apr\`es 1.4, ${\cal C}(q,q')$ est un rev\^etement de $\Xi^*(q,q')$. Avec les d\'efinitions pr\'ec\'edentes, on obtient l'\'egalit\'e
 $$(8) \qquad m_{geom}(\sigma,\sigma')=\int_{{\cal C}(q,q')}c_{\sigma}(x)c_{\sigma'}(x)D^G(x)\Delta(x)^rdx.$$
 
 \bigskip
 
 \subsection{Facteurs $\epsilon$ de paires}

   Soient $U$ et $U'$ deux espaces vectoriels sur $F$ de dimensions $d$ et $d'$. On suppose $d$ pair et $d'$ impair. On note $M$ et $M'$ leurs groupes lin\'eaires et on introduit les ensembles $\tilde{M}$ et $\tilde{M}'$ comme en 1.2. Posons $r=(d'-d-1)/2$ si $d'>d$, $r=(d-d'-1)/2$ si $d>d'$. Fixons un isomorphisme $U'=U\oplus Z_{2r+1}$ si $d'>d$, $U=U'\oplus Z_{2r+1}$ si $d>d'$, avec la notation introduite en 2.1. Fixons un \'el\'ement $\nu_{1}\in F^{\times}$. Si $d'>d$, on plonge $\tilde{M}$ dans $\tilde{M}'$ en identifiant toute forme bilin\'eaire $\tilde{m}$ sur $U$ \`a la "somme directe" de $\tilde{m}$ et de la forme $q_{2r+1,-\nu_{1}}$ sur $Z_{2r+1}$. Si $d>d'$, on plonge $\tilde{M}'$ dans $\tilde{M}$ en identifiant toute forme bilin\'eaire $\tilde{m}'$ sur $U'$ \`a la somme directe de $\tilde{m}'$ et de la forme $q_{2r+1,\nu_{1}}$ sur $Z_{2r+1}$.
 
   Soit $\pi\in Temp(M(F))$, supposons $\pi$ autoduale, c'est-\`a-dire isomorphe \`a sa contragr\'ediente $\pi^{\vee}$. Soit de m\^eme $\pi'\in  Temp(M'(F))$, supposons $\pi'$ autoduale. On note $\omega_{\pi}$ et $\omega_{\pi'}$ les caract\`eres centraux de $\pi$ et $\pi'$ et $\epsilon(s,\pi\times \pi',\psi_{F})$ le facteur $\epsilon$ d\'efini par Jacquet, Piatetski-Shapiro et Shalika. On pose
 $$\epsilon_{\nu_{1}}(\pi,\pi')=\epsilon_{\nu_{1}}(\pi',\pi)=\omega_{\pi}((-1)^{(d'-1)/2}2\nu_{1})\omega_{\pi'}((-1)^{1+d/2}2\nu_{1})\epsilon(1/2,\pi\times \pi',\psi_{F}).$$
 
 {\bf Remarque.} Ici encore, si $d'>d$, la formule n'est pas la m\^eme que celle de [W3] en ce qui concerne les signes. Mais on a aussi modifi\'e dans ce cas le plongement de $\tilde{M}$ dans $\tilde{M}'$.
 
 Dans [W3], on a calcul\'e ce nombre $\epsilon_{\nu_{1}}(\pi,\pi')$ par une formule int\'egrale que nous allons rappeler. On peut  prolonger la repr\'esentation $\pi$ \`a $\tilde{M}(F)$. C'est-\`a-dire que l'on peut d\'efinir une application $\tilde{\pi}$ de $\tilde{M}(F)$ dans le groupe des automorphismes lin\'eaires de $E_{\pi}$ de sorte que $\tilde{\pi}(m\tilde{m}m')=\pi(m)\tilde{\pi}(\tilde{m})\pi(m')$ pour tous $m,m'\in M(F)$, $\tilde{m}\in \tilde{M}(F)$. On peut supposer $\tilde{\pi}$ unitaire. Alors $\tilde{\pi}$ est uniquement d\'efinie \`a une homoth\'etie pr\`es de rapport un nombre complexe de module $1$. Pour d\'eterminer enti\`erement $\tilde{\pi}$, on utilise la th\'eorie des mod\`eles de Whittaker. Fixons une base $(u_{j})_{j=1,...,d}$ de $U$, qui permet d'identifier $M$ au groupe matriciel $GL_{d}$. On note $N$ le groupe des matrices unipotentes triangulaires sup\'erieures. On  appelle fonctionnelle de Whittaker sur $E_{\pi}$ toute forme lin\'eaire $\phi:E_{\pi}\to {\mathbb C}$ telle que $\phi(\pi(n)e)=\psi_{F}(\sum_{j=1,...,d-1}n_{j,j+1})\phi(e)$ pour tous $n\in N(F)$ et $e\in E_{\pi}$. On sait que l'espace de ces fonctionnelles est une droite.  On a d\'efini un \'el\'ement $\boldsymbol{\theta}_{d}$ de $\tilde{M}(F)$  en 1.8. On v\'erifie que l'application $\phi\mapsto \phi\circ\tilde{\pi}(\boldsymbol{\theta}_{d})$ conserve la droite des fonctionnelles de Whittaker. On peut alors normaliser $\tilde{\pi}$ de sorte que cette application soit l'identit\'e. On v\'erifie  que cette normalisation  d\'epend de $\psi_{F}$ mais pas de la base $(u_{j})_{j=1,...,d}$ choisie. On prolonge de la m\^eme fa\c{c}on la repr\'esentation $\pi'$ en une repr\'esentation $\tilde{\pi}'$ de $\tilde{M}'(F)$, que l'on normalise de m\^eme.
  
 Supposons pour fixer les id\'ees $d<d'$.  On va d\'efinir un ensemble $\underline{{\cal T}}$ de "sous-tores" de $\tilde{M}$.  Consid\'erons une d\'ecomposition $U=W_{T}\oplus U_{T}$ telle que $dim(W_{T})$ soit paire. Notons $(M_{T},\tilde{M}_{T})$ le groupe tordu associ\'e \`a $W_{T}$. Soit $(T,\tilde{T}_{W})$ un sous-tore tordu maximal de $(M_{T},\tilde{M}_{T})$. Une fa\c{c}on de d\'efinir cette notion est de dire qu'il existe un \'el\'ement fortement r\'egulier $\tilde{m}\in \tilde{M}_{T}(F)$ de sorte que, en notant $T^{\flat}$ le commutant de $\tilde{m}$ dans $M_{T}$,  $T$ soit le commutant dans $M_{T}$ de $T^{\flat}$ et que $\tilde{T}_{W}$ soit \'egal \`a $T\tilde{m}$.  Remarquons que $T^{\flat}$ ne d\'epend que de $\tilde{T}_{W}$ et pas de $\tilde{m}$. On suppose   $T^{\flat}$ anisotrope. Soit $q_{M,T}$ une forme bilin\'eaire sym\'etrique non d\'eg\'en\'er\'ee sur $U_{T}$. On note $\tilde{T}$ l'ensemble des \'el\'ements de $\tilde{M}$ qui sont somme directe d'un \'el\'ement de $\tilde{T}_{W}$ sur $W_{T}$ et de l'\'el\'ement $q_{M,T}$ sur $U_{T}$. On note $U'_{T}=U_{T}\oplus Z_{2r+1}$ et $q_{M',T}$ la somme directe de $q_{M,T}$ et de $q_{2r+1,-\nu_{1}}$. On suppose que les groupes sp\'eciaux orthogonaux des formes $q_{M,T}$ et $q_{M',T}$ sont quasi-d\'eploy\'es. L'ensemble $\underline{\cal T}$ est l'ensemble des "sous-tores" $\tilde{T}$ v\'erifiant toutes ces conditions. Remarquons qu'\`a tout tel \'el\'ement $\tilde{T}$ est associ\'e un sous-tore maximal $T$ de $M_{T}$ ainsi qu'un automorphisme de ce tore que l'on note $\theta$: c'est l'\'el\'ement $\theta_{\tilde{t}}$ pour n'importe quel \'el\'ement $\tilde{t}\in \tilde{T}$. On note $T_{\theta}$ la composante neutre du sous-groupe des points fixes de $\theta$ dans $T$ (c'est le tore not\'e $T^{\flat}$ ci-dessus). On note $\tilde{T}(F)_{/\theta}$ le quotient de $\tilde{T}(F)$ par la relation d'\'equivalence: deux \'el\'ements  sont \'equivalents si et seulement s'ils sont conjugu\'es par un \'el\'ement de $T(F)$. On munit $\tilde{T}(F)_{/\theta}$ d'une mesure de sorte que, pour tout $\tilde{t}\in \tilde{T}(F)_{/\theta}$, l'application $t\mapsto t\tilde{t}$ de $T_{\theta}(F)$ dans $\tilde{T}(F)_{/\theta}$ conserve localement les mesures. Consid\'erant que le groupe $M$ agit par conjugaison sur $\tilde{M}$, on d\'efinit le normalisateur $Norm_{M}(\tilde{T})$. On note $G_{q_{M,T}}$ le groupe sp\'ecial orthogonal de $(U_{T},q_{M,T})$. On pose
 $$W(M,\tilde{T})=Norm_{M}(\tilde{T})(F)/(T(F)\times G_{q_{M,T}}(F)).$$
 Ce groupe est fini. Les repr\'esentations $\tilde{\pi}$ et $\tilde{\pi}'$ admettent des caract\`eres, lesquels ont localement des d\'eveloppements comme dans le cas des groupes non tordus ( [C] th\'eor\`eme 3). On peut copier les d\'efinitions de 2.1 (en rempla\c{c}ant $\nu_{0}$ par $\nu_{1}$) pour obtenir des fonctions $c_{\tilde{\pi}}$ et $c_{\tilde{\pi}'}$ presque partout sur $\tilde{T}(F)$. Il reste \`a d\'efinir une fonction $\Delta$ sur $\tilde{T}(F)$. Pour cela, soit $\tilde{t}$ un \'el\'ement fortement r\'egulier de $ \tilde{T}(F)$. Notons $\xi \in \Xi$ son image par l'application $p_{\tilde{M}_{T}}$ de 1.4. On pose $\Delta(\tilde{t})= \Delta(\xi)$.
 
  Le groupe $M(F)$ agit par conjugaison sur $\underline{{\cal T}}$. On fixe un ensemble de repr\'esentants ${\cal T}$ des classes d'\'equivalence. Posons
 $$\epsilon_{geom,\nu_{1}}(\pi,\pi')=\epsilon_{geom,\nu_{1}}(\tilde{\pi}',\tilde{\pi})=\sum_{\tilde{T}\in {\cal T}}\vert  2\vert _{F}^{r^2+r+rdim(U_{T})}\vert W(M,\tilde{T})\vert ^{-1}$$
 $$\int_{\tilde{T}(F)_{/\theta}}c_{\tilde{\pi}}(\tilde{t})c_{\tilde{\pi}'}(\tilde{t})D^{\tilde{M}}(\tilde{t})\Delta(\tilde{t})^rd\tilde{t}.$$
 Cette expression est absolument convergente et le th\'eor\`eme principal de [W3] affirme l'\'egalit\'e
 $$\epsilon_{geom,\nu_{1}}(\pi,\pi')=\epsilon_{\nu_{1}}(\pi,\pi').$$
 Pour d\'emontrer ce th\'eor\`eme, on a besoin de r\'esultats issus de la formule des traces locale tordue. Cette formule n'existe pas encore dans la litt\'erature. On a admis les r\'esultats en question.
 
 \bigskip
 
 \subsection{Reformulation de $\epsilon_{geom,\nu_{1}}(\pi,\pi')$}
 
 Notons ${\cal D}(d,d')$ l'ensemble des entiers  ${\bf d}$   qui v\'erifient  la condition
 
 (1) ${\bf d}$ est pair et $0\leq {\bf d}\leq inf(d,d')$.
 
 Pour tout \'el\'ement ${\bf d}$ de cet ensemble,  posons
 $${\cal M}({\bf d})=\left\lbrace\begin{array}{cc}F^{\times}/F^{\times,2},&\text{ si }{\bf d}<d,\\ \{-\nu_{1}\},&\text{ si }{\bf d}=d.\\ \end{array}\right.$$ 
 Supposons $d<d'$ (des consid\'erations similaires s'appliquent au cas $d>d'$ en \'echangeant les r\^oles de $M$ et $M'$ et en changeant $\nu_{1}$ en $-\nu_{1}$). Pour ${\bf d}\in {\cal D}(d,d')$ et $\mu\in {\cal M}({\bf d})$, fixons une d\'ecomposition $U=W_{{\bf d}}\oplus U_{{\bf d}}$, o\`u $dim(W_{{\bf d}})={\bf d}$. Une preuve similaire \`a celle de 2.2(5) montre qu'\`a \'equivalence pr\`es,il existe une unique forme bilin\'eaire sym\'etrique non d\'eg\'en\'er\'ee $q_{{\bf d},\mu}$ sur $U_{{\bf d}}$ telle que, d'une part le groupe sp\'ecial orthogonal de $(U_{{\bf d}},q_{{\bf d},\mu})$ soit quasi-d\'eploy\'e, d'autre part la forme $q_{{\bf d},\mu}\oplus q_{2r+1,-\nu_{1}}$ sur $U_{{\bf d}}\oplus Z_{2r+1}$ soit \'equivalente \`a la somme de plans hyperboliques et de la forme $(\lambda,\lambda')\mapsto 2\mu\lambda\lambda'$ sur $F$. On fixe une telle forme. 
 
   Soit $\tilde{T}\in \underline{\cal T}$. On v\'erifie qu'il existe un unique couple $({\bf d},\mu)$, avec ${\bf d}\in  {\cal D}(d,d')$ et $\mu\in {\cal M}({\bf d})$,  tel que le triplet $(W_{T},U_{T},q_{M,T})$ soit conjugu\'e par un \'el\'ement de $M(F)$ \`a $(W_{{\bf d}},U_{{\bf d}},q_{{\bf d},\mu})$. On peut donc supposer que pour tout $\tilde{T}\in {\cal T}$, le triplet $(W_{T},U_{T},q_{M,T})$ est \'egal \`a $(W_{{\bf d}},U_{{\bf d}},q_{{\bf d},\mu})$ pour un tel couple $({\bf d},\mu)$. Cela d\'ecompose ${\cal T}$ en union disjointe de sous-ensembles ${\cal T}_{{\bf d},\mu}$. Fixons $({\bf d},\mu)$ et consid\'erons la contribution de ${\cal T}_{{\bf d},\mu}$ \`a $\epsilon_{geom,\nu_{1}}(\pi,\pi')$. Introduisons le groupe tordu $(M_{{\bf d}},\tilde{M}_{{\bf d}})$ associ\'e \`a $W_{{\bf d}}$. A chaque $\tilde{T}\in {\cal T}_{{\bf d},\mu}$ est associ\'e un sous-tore tordu maximal $(T_{{\bf d}},\tilde{T}_{{\bf d}})$ (que l'on notera simplement $\tilde{T}_{{\bf d}}$) de ce groupe tordu de sorte que $\tilde{T}$ soit l'ensemble des sommes directes d'un \'el\'ement de $\tilde{T}_{{\bf d}}$ sur $W_{{\bf d}}$ et de $q_{{\bf d},\mu}$ sur $U_{{\bf d}}$. Quand $\tilde{T}$ d\'ecrit ${\cal T}_{{\bf d},\mu}$, le tore $\tilde{T}_{{\bf d}}$ d\'ecrit un ensemble de repr\'esentants des classes de conjugaison par $M_{{\bf d}}(F)$ dans l'ensemble des sous-tores tordus maximaux  $\tilde{T}'$ de $\tilde{M}_{{\bf d}}$ tels que $\tilde{T}'_{\theta}$ soit anisotrope. Le groupe $Norm_{M}(\tilde{T})$ est le produit de $Norm_{M_{{\bf d}}}(\tilde{T}_{{\bf d}})$ et du groupe orthogonal $G^+_{q_{{\bf d},\mu}}$ de l'espace quadratique $(U_{{\bf d}},q_{{\bf d},\mu})$. Donc
 $$W(M,\tilde{T})=W(M_{{\bf d}},\tilde{T}_{{\bf d}})\times( G^+_{q_{{\bf d},\mu}}(F)/G_{q_{{\bf d},\mu}}(F)).$$
 Le dernier quotient a $2 $ \'el\'ements si ${\bf d}<d$, un seul si ${\bf d}=d$. Cela conduit \`a poser
 $$c({\bf d})=\left\lbrace\begin{array}{cc}1/2,&\text{ si }{\bf d}<inf(d,d'),\\ 1,&\text{ si }{\bf d}=inf(d,d').\\ \end{array}\right.$$
 D'autre part, les fonctions $c_{\tilde{\pi}}$ et $c_{\tilde{\pi}'}$ sur $\tilde{T}(F)$ s'identifient \`a des fonctions sur $\tilde{T}_{{\bf d}}(F)$. On peut en fait \'etendre leur d\'efinition et les d\'efinir sur tout $\tilde{M}_{{\bf d},reg}(F)$. Il faut prendre garde au fait qu'elles d\'ependent de $\mu$. Plus pr\'ecis\'ement, la donn\'ee de $\mu$ d\'efinit un plongement de $\tilde{M}_{{\bf d}}$ dans $\tilde{M}$: on identifie tout \'el\'ement $\tilde{m}_{{\bf d}}\in \tilde{M}_{{\bf d}}$ \`a l'\'el\'ement de $\tilde{M}$ qui est la somme directe de $\tilde{m}_{{\bf d}}$ sur $W_{{\bf d}}$ et de $q_{{\bf d},\mu}$ sur $U_{{\bf d}}$.   Posons
 $${\cal X}(d,d')=\sqcup_{{\bf d}\in {\cal D}(d,d')}( {\cal M}({\bf d})\times  \tilde{M}_{{\bf d}}(F)_{ani}/conj).$$
 On peut consid\'erer que les fonctions $c_{\tilde{\pi}}$ et $c_{\tilde{\pi}'}$ sont d\'efinies sur cette vari\'et\'e. Il en est de m\^eme des fonctions $D^{\tilde{M}}$ et $\Delta$. Il y a aussi une fonction \'evidente sur cette vari\'et\'e qui associe \`a tout point l'indice ${\bf d}$ de la composante qui contient ce point.    Avec ces d\'efinitions, on obtient l'\'egalit\'e
 $$(1) \qquad \epsilon_{geom,\nu_{1}}=\int_{{\cal X}(d,d')}c({\bf d})\vert 2\vert _{F}^{r^2+r+r(d-{\bf d})}c_{\tilde{\pi}}(\tilde{x})c_{\tilde{\pi}'}(\tilde{x})D^{\tilde{M}}(\tilde{x})\Delta(\tilde{x})^rd\tilde{x}.$$
  Posons
 $$\Xi^*(d,d')=\sqcup_{{\bf d}\in {\cal D}(d,d')}\Xi^*_{reg,{\bf d}}.$$
 D'apr\`es 1.4,  il y a une application naturelle de ${\cal X}(d,d')$ vers $\Xi^*(d,d')$ qui est un rev\^etement analytique. On prendra garde au fait que cette application ne pr\'eserve pas localement les mesures: le calcul de 1.4 montre que le jacobien est $\vert 2\vert_{F} ^{{\bf d}/2} $.

  \bigskip
  
  \section{Endoscopie}
  
  \bigskip
  
  \subsection{Rappels sur le calcul des fonctions $c_{\sigma,{\cal O}}$}
  
  Consid\'erons un espace quadratique $(V',q')$ de dimension $d'$ impaire. Soit $V'=W\oplus V'_{\sharp}$ une d\'ecomposition orthogonale telle que $dim(W)$ soit paire. On note $G'$, $G_{W}$ et $G'_{\sharp}$ les  groupes sp\'eciaux orthogonaux de $V'$, $W$, $V'_{\sharp}$. On suppose $G'_{\sharp}$ d\'eploy\'e. Soient $\sigma'$ une repr\'esentation admissible irr\'eductible de $G'(F)$ et $t$ un \'el\'ement semi-simple de $G_{W}(F)$ tel que $G'_{t}=T\times G'_{\sharp}$, o\`u $T$ est un sous-tore de $G_{W}$. Le caract\`ere $\Theta_{\sigma'}$ se d\'eveloppe au voisinage de $t$ selon 2.1(5), et l'ensemble $Nil(\mathfrak{g}'_{t}(F))=Nil(\mathfrak{g}'_{\sharp}(F))$ contient une unique orbite r\'eguli\`ere ${\cal O}_{reg}$. Nous allons rappeler une formule qui calcule le coefficient $c_{\sigma',{\cal O}_{reg}}(t)$. Fixons un sous-tore d\'eploy\'e maximal $T_{\sharp}$ de $G'_{\sharp}$ et un \'el\'ement r\'egulier $Y_{\sharp}\in \mathfrak{t}_{\sharp}(F)$. On d\'efinit la fonction $D^{G'_{\sharp}}$ sur $\mathfrak{g}'_{\sharp}(F)$ comme on l'a d\'efinie sur $G'_{\sharp}(F)$. On a alors l'\'egalit\'e
  $$(1) \qquad c_{\sigma',{\cal O}_{reg}}(t)=\vert W(G'_{\sharp},T_{\sharp})\vert ^{-1}lim_{\lambda\in F^{\times},\lambda\to 0}\Theta_{\sigma'}(texp(\lambda Y_{\sharp}))D^{G'_{\sharp}}(\lambda Y_{\sharp})^{1/2}$$
  cf. [W1] p.115.
  
  Consid\'erons maintenant un espace quadratique $(V,q)$ de dimension $d$ paire. Soit $V=W\oplus V_{\sharp}$ une d\'ecomposition orthogonale telle que $dim(W)$ soit paire. On note $G$, $G_{W}$ et $G_{\sharp}$ les  groupes sp\'eciaux orthogonaux de $V$, $W$, $V_{\sharp}$. On suppose $G_{\sharp}$ quasi-d\'eploy\'e. Soient $\sigma$ une repr\'esentation admissible irr\'eductible de $G(F)$ et $t$ un \'el\'ement semi-simple de $G_{W}(F)$ tel que $G_{t}=T\times G_{\sharp}$, o\`u $T$ est un sous-tore de $G_{W}$. Fixons un sous-tore maximal $T_{\sharp}$ de $G_{\sharp}$ contenu dans un sous-groupe de Borel, et un \'el\'ement r\'egulier $Y_{\sharp}\in \mathfrak{t}_{\sharp}(F)$. Si $dim(V_{\sharp})\leq2$, il n'y a encore qu'une orbite nilpotente r\'eguli\`ere ${\cal O}_{reg}$ dans $\mathfrak{g}_{\sharp}(F)$ et on a l'\'egalit\'e
$$(2) \qquad c_{\sigma,{\cal O}_{reg}}(t)=\vert W(G_{\sharp},T_{\sharp})\vert ^{-1}lim_{\lambda\in F^{\times},\lambda\to 0}\Theta_{\sigma}(texp(\lambda Y_{\sharp}))D^{G_{\sharp}}(\lambda Y_{\sharp})^{1/2}.$$
Supposons maintenant $dim(V_{\sharp})\geq 4$. Notons $q_{\sharp}$ la restriction de $q$ \`a $V_{\sharp}$ et  $\delta_{\sharp}$ son discriminant. Si $\delta_{\sharp}=1$, l'ensemble des orbites nilpotentes r\'eguli\`eres dans $\mathfrak{g}_{\sharp}(F)$ est param\'etr\'e par ${\cal N}_{\sharp}=F^{\times}/F^{\times,2}$, cf.  1.7.  On pose $\eta_{\sharp}=1$ et on note ${\cal F}_{\sharp}$ l'ensemble des couples $(F',F')$ o\`u $F'$ est une extension quadratique de $F$ (il s'agit de vraies extensions quadratiques, l'alg\`ebre $F'=F\oplus F$ est exclue; la m\^eme remarque vaut ci-dessous). Si $\delta_{\sharp}\not=1$, posons $E_{\sharp}=F(\sqrt{\delta_{\sharp}})$. Fixons $\eta_{\sharp}\in F^{\times}$ tel que $q_{\sharp}$ soit \'equivalente \`a une somme orthogonale de plans hyperboliques et de l'espace $E_{\sharp}$ muni de la forme $(v,v')\mapsto \eta_{\sharp}trace_{E_{\sharp}/F}(\tau_{E_{\sharp}/F}(v)v')$. L'ensemble des orbites nilpotentes r\'eguli\`eres dans $\mathfrak{g}_{\sharp}(F)$ est param\'etr\'e par ${\cal N}_{\sharp}=\eta_{\sharp} Norm_{E_{\sharp}/F}(E_{\sharp}^{\times})/F^{\times,2}\subset F^{\times}/F^{\times,2} $. Consid\'erons l'ensemble des couples d'extensions quadratiques $(F_{1},F_{2})=(F(\sqrt{\delta_{1}}),F(\sqrt{\delta_{2}})) $ de $F$ tels que $\delta_{1}\delta_{2}=\delta_{\sharp}$.On y introduit la relation d'\'equivalence: $(F_{1},F_{2})$ est \'equivalent \`a $(F_{2},F_{1})$. On fixe un ensemble de repr\'esentants ${\cal F}_{\sharp}$ des classes d'\'equivalence. On a param\'etr\'e en 1.4 les classes de conjugaison d'\'el\'ements semi-simples r\'eguliers de $G_{\sharp}(F)$. Un param\'etrage analogue vaut pour les \'el\'ements semi-simples r\'eguliers de $\mathfrak{g}_{\sharp}(F)$: il suffit de remplacer les donn\'ees $y_{i}$ telles que $y_{i}\tau_{i}(y_{i})=1$ par des donn\'ees  $Y_{i}$ telles que $Y_{i}+\tau_{i}(Y_{i})=0$. Posons $J=\{1,...,dim(V_{\sharp})/2\}$. Pour $j\in J$ posons $F_{\pm j}=F$ et, si $j\geq3$, $F_{j}=F\oplus F$. Consid\'erons un couple $(F_{1},F_{2})\in {\cal F}_{\sharp}$.   Fixons une famille $(Y_{j})_{j\in J}$, en "position g\'en\'erale", telle que, pour tout $j\in J$, $Y_{j}$ soit un \'el\'ement de $F_{j}$ tel que $Y_{j}+\tau_{j }(Y_{j})=0$. Fixons deux familles $c^+=(c^+_{j})_{j=1,2}$ et $c^-=(c^-_{j})_{j=1,2}$ d'\'el\'ements de $F^{\times}$ telles que  pour  $\epsilon=\pm$,
$$sgn_{F_{1}/F}(c_{1}^{\epsilon})=\epsilon sgn_{F_{1}/F}(\eta_{\sharp}(1-Norm_{F_{1}/F}(Y_{1})Norm_{F_{2}/F}(Y_{2})^{-1})),$$
$$sgn_{F_{2}/F}(c_{2}^{\epsilon})=\epsilon sgn_{F_{2}/F}(\eta_{\sharp}(1-Norm_{F_{2}/F}(Y_{2})Norm_{F_{1}/F}(Y_{1})^{-1})),$$
On v\'erifie que, pour $\epsilon=\pm$, les familles $\xi=(J,(F_{\pm j})_{j\in J},(F_{j})_{j\in J},(Y_{j})_{j\in J})$ et $c^{\epsilon}$ param\`etrent une classe de conjugaison d'\'el\'ements semi-simples r\'eguliers dans $\mathfrak{g}_{\sharp}(F)$. Plus exactement, elles en param\`etrent deux, qui sont conjugu\'ees par le groupe orthogonal. On fixe un \'el\'ement $Y_{F_{1},F_{2}}^{\epsilon}$ dans l'une de ces classes. On suppose que $Y_{F_{1},F_{2}}^+$ et $Y_{F_{1},F_{2}}^-$ sont stablement conjugu\'es. On note $T_{F_{1},F_{2}}^{\epsilon}$ le commutant de $Y_{F_{1},F_{2}}^{\epsilon}$ dans $G_{\sharp}$. Pour $\mu\in {\cal N}_{\sharp}$, on a alors l'\'egalit\'e
$$(3) \qquad c_{\sigma,{\cal O}_{\mu}}(t)=lim_{\lambda\in F^{\times,2}, \lambda\to 0}( \vert W(G_{\sharp},T_{\sharp})\vert ^{-1}\Theta_{\sigma}(texp(\lambda Y_{\sharp}))D^{G_{\sharp}}(\lambda Y_{\sharp})^{1/2}$$
$$+\frac{1}{2}\sum_{(F_{1},F_{2})\in {\cal F}_{\sharp}}\sum_{\epsilon=\pm}\epsilon \vert W(G_{\sharp},T_{F_{1},F_{2}}^{\epsilon})\vert ^{-1}sgn_{F_{1}/F}(\mu \eta_{\sharp})\Theta_{\sigma}(t exp(\lambda Y_{F_{1},F_{2}}^{\epsilon}))D^{G_{\sharp}}(\lambda Y_{F_{1},F_{2}}^{\epsilon})^{1/2})$$
cf. [W1], p.115 (les d\'efinitions de cette r\'ef\'erence sont un peu plus compliqu\'ees, Dieu sait pourquoi). Remarquons que l'on pourrait remplacer le terme $sgn_{F_{1}/F}(\mu\eta_{\sharp})$ par $sgn_{F_{2}/F}(\mu\eta_{\sharp})$: si $\delta_{\sharp}=1$, $F_{1}=F_{2}$; si $\delta_{\sharp}\not=1$, le produit de ces deux termes vaut $sgn_{E_{\sharp}/F}(\mu \eta_{\sharp})=1$. La formule suivante nous sera utile:
$$(4)\qquad \vert {\cal N}_{\sharp}\vert ^{-1}\sum_{\mu\in {\cal N}_{\sharp}}c_{\sigma,{\cal O}_{\mu}}(t)=lim_{\lambda\to 0}\vert W(G_{\sharp},T_{\sharp})\vert ^{-1}\Theta_{\sigma}(texp(\lambda Y_{\sharp}))D^{G_{\sharp}}(\lambda Y_{\sharp})^{1/2}.$$
En effet, cette formule est imm\'ediate quand $dim(V_{\sharp})\leq 2$. Supposons $dim(V_{\sharp})\geq4$. Quand $\mu$ d\'ecrit ${\cal N}_{\sharp}$, $\mu\eta_{\sharp}$ d\'ecrit $F^{\times}/F^{\times,2}$ si $\delta_{\sharp}=1$, le sous-groupe $Norm_{E_{\sharp}/F}(E_{\sharp}^{\times})/F^{\times,2}$ quand $\delta_{\sharp}\not=1$. Pour $(F_{1},F_{2})\in {\cal F}_{\sharp}$, les conditions $\delta_{1}\delta_{2}=\delta_{\sharp}$ et $\delta_{1}\not=1$, $\delta_{2}\not=1$ entra\^{\i}nent que $sgn_{F_{1}/F}$ est non trivial sur le groupe pr\'ec\'edent. Donc
$$\sum_{\mu\in {\cal N}_{\sharp}}sgn_{F_{1}/F}(\mu\eta_{\sharp})=0,$$
et (4) se d\'eduit de (3).

\bigskip

\subsection{D\'efinition d'une multiplicit\'e stable}

On consid\`ere deux couples $(V,q)$ et $(V',q')$ d'espaces quadratiques de dimensions $d$ et $d'$. On suppose $d$ pair et $d'$ impair. On note $G$ et $G'$ les groupes sp\'eciaux orthogonaux et on les suppose quasi-d\'eploy\'es. On consid\`ere, dans un groupe de Grothendieck convenable, une combinaison lin\'eaire finie
$$\Sigma=\sum_{k}a_{k}\sigma_{k} $$
o\`u les $a_{k}$ sont des nombres complexes et les $\sigma_{k}$ appartiennent \`a $Temp(G(F))$. On pose $\Theta_{\Sigma}=\sum_{k}a_{k}\Theta_{\sigma_{k}}$ et on suppose que $\Theta_{\Sigma}$ est  stable en tant que distribution ou, ce qui revient au m\^eme, qu'elle est constante, en tant que fonction, sur les classes de conjugaison stable. On consid\`ere de m\^eme une combinaison lin\'eaire finie 
$$\Sigma'=\sum_{k}a'_{k}\sigma'_{k},$$
o\`u les $\sigma'_{k}$ appartiennent \`a $Temp(G'(F))$. On d\'efinit $\Theta_{\Sigma'}$ et on suppose que $\Theta_{\Sigma'}$ est stable.

On d\'efinit $\Xi^*(q,q')$ comme en 2.2. Soit $\xi=(I,(F_{\pm i})_{i\in I},(F_{i})_{i\in I},(y_{i})_{i\in I})\in \Xi^*(q,q')$. On pose $I_{\xi}=I=I^*$. En imitant la preuve de 2.2(5), on voit qu'il existe \`a \'equivalence pr\`es une unique d\'ecomposition orthogonale $(V',q')=(W,q_{W})\oplus (V'_{\sharp},q'_{\sharp})$ telle que $dim(W)=d_{\xi}$, $\delta(q_{W})=\delta_{\xi}$ et le groupe sp\'ecial orthogonal  de $(V'_{\sharp},q'_{\sharp})$ soit d\'eploy\'e. L'\'el\'ement $\xi$ param\`etre deux classes de conjugaison stable dans le groupe sp\'ecial orthogonal de $(W,q_{W})$. Fixons un \'el\'ement $t'$ dans l'une d'elles. Par lin\'earit\'e, on peut d\'efinir le terme $c_{\Sigma',{\cal O}_{reg}}(t')$.  Il est facile de voir que, parce que $\Sigma'$ est stable, il ne d\'epend pas du choix de $t'$. On le note simplement $c_{\Sigma'}(\xi)$.  

Supposons $d_{\xi}<d$. Il existe de m\^eme une d\'ecomposition orthogonale $(V,q)=(W,q_{W})\oplus (V_{\sharp},q_{\sharp})$ telle que $dim(W)=d_{\xi}$, $\delta(q_{W})=\delta_{\xi}$ et le groupe sp\'ecial orthogonal $G_{\sharp}$ de $(V_{\sharp},q_{\sharp})$ soit quasi-d\'eploy\'e. L'espace $(W,q_{W})$ n'est pas en g\'en\'eral le m\^eme que ci-dessus et la d\'ecomposition n'est pas forc\'ement unique. Fixons-la. De nouveau, $\xi$ param\`etre deux classes de conjugaison stable dans le groupe sp\'ecial orthogonal de $(W,q_{W})$. Fixons un \'el\'ement $t$ dans l'une d'elles. Notons ${\cal N}_{\sharp}$ l'ensemble des orbites nilpotentes r\'eguli\`eres dans $\mathfrak{g}_{\sharp}(F)$. Par lin\'earit\'e, on d\'efinit le terme $c_{\Sigma,{\cal O}}(t)$ pour  ${\cal O}\in {\cal N}_{\sharp}$. Posons
$$c_{\Sigma}(\xi)=\vert {\cal N}_{\sharp}\vert ^{-1}\sum_{{\cal O}\in {\cal N}_{\sharp}}c_{\Sigma,{\cal O}}(t).$$
Ce terme ne d\'epend pas du choix de $t$. En effet, la formule 3.1(4) le calcule comme la limite d'une expression o\`u $t$ n'intervient que par une valeur $\Theta_{\Sigma}(texp(\lambda Y_{\sharp}))$ du caract\`ere $\Theta_{\Sigma}$. Quand on change le choix de $t$, disons qu'on le remplace par $t_{1}$, on  peut choisir $Y_{1,\sharp}$ qui v\'erifie relativement \`a $t_{1}$ les m\^emes conditions que $Y_{\sharp}$ v\'erifiait relativement \`a $t$ et qui est tel que $t_{1}exp(\lambda Y_{1,\sharp})$ soit stablement conjugu\'e \`a $texp(\lambda Y_{\sharp})$. L'ind\'ependance du choix de $t$ r\'esulte alors de la stabilit\'e de $\Sigma$.

Supposons maintenant $d_{\xi}=d$. L'\'el\'ement $\xi$ param\`etre deux classes de conjugaison stable dans $G(F)$, qui sont conjugu\'ees par le groupe orthogonal. On fixe des repr\'esentants $t_{1}$ et $t_{2}$ de chacune d'elles et on pose $c_{\Sigma}(\xi)=\Theta_{\Sigma}(t_{1})+\Theta_{\Sigma}(t_{2})$.

  On pose $r=(\vert d-d'\vert -1)/2$ et
$$(1) \qquad S(\Sigma,\Sigma')=S(\Sigma',\Sigma)=\int_{\Xi^*(q,q')}2^{d_{\xi}}c_{\Sigma}(\xi)c_{\Sigma'}(\xi)D^{inf(d,d')}(\xi)\Delta(\xi)^rd\xi,$$
cf. 1.5 pour les notations.

\bigskip

\subsection{ Transfert de multiplicit\'es}

On consid\`ere une paire d'espaces quadratiques $(V,q)$ et $(V',q')$ comme en 2.1. On consid\`ere de plus quatre espaces quadratiques $(V_{+},q_{+})$, $(V_{-},q_{-})$, $(V'_{+},q'_{+})$, $(V'_{-},q'_{-})$. On utilisera quelques notations que l'on esp\`ere maintenant \'evidentes: $d_{+}$ est la dimension de $V_{+}$, $G_{+}$ est le groupe sp\'ecial orthogonal de $(V_{+},q_{+})$ etc... On impose les hypoth\`eses suivantes:

(1) $d_{+}$ et $d_{-}$ sont pairs, $d_{+}+d_{-}=d$, $\delta(q_{+})\delta(q_{-})=\delta(q)$;

(2) $d'_{+}$ et $d'_{-}$ sont impairs, $d'_{+}+d'_{-}=d'+1$;

(3) les groupes $G'_{+}$, $G'_{-}$ sont d\'eploy\'es et les espaces $(V_{+},q_{+})$ et $(V_{-},q_{-}) $ v\'erifient la condition (QD) de 1.7.

Le groupe $G_{+}\times G_{-}$ est un groupe endoscopique de $G$ et $G'_{+}\times G'_{-}$ est un groupe endoscopique de $G'$. On utilise les facteurs de transfert d\'efinis en 1.7. Dans les formules (1) et (3) de ce paragraphe interviennent des termes $\Delta_{-2\delta(q')}(\xi_{+},\xi_{-},c)$ et $\Delta_{\nu_{0}\delta(q)}(\xi_{+},\xi_{-},c)$. On remarque que ce sont les m\^emes car $-2\delta(q')=\nu_{0}\delta(q)$. On note simplement ce terme $\Delta(\xi_{+},\xi_{-},c)$. Les hypoth\`eses de 2.1 entra\^{\i}nent que  $G'$ est d\'eploy\'e si et seulement si l'hypoth\`ese (QD) de ce paragraphe est v\'erifi\'ee. Donc, ou bien on a les \'egalit\'es
 
 $(\underline{V},\underline{q})=(V,q)$ et $(\underline{V}',\underline{q}')=(V',q')$;
 
 \noindent ou bien on a  
 
 $(\underline{V},\underline{q})\not=(V,q)$ et $(\underline{V}',\underline{q}')\not=(V',q')$.

  Consid\'erons  une combinaison lin\'eaire finie
 $$\Sigma'_{+}=\sum_{k}b'_{k,+}\sigma'_{k,+}$$
 o\`u les $\sigma'_{k,+}$ appartiennent \`a $Temp(G'_{+}(F))$. On consid\`ere  des combinaisons lin\'eaires analogues $\Sigma'_{-}$, $\Sigma'$, $\Sigma_{+}$, $\Sigma_{-}$ et $\Sigma$ relatives aux groupes $G'_{-}$, $G'$, $G_{+}$, $G_{-}$ et $G$. On suppose
 
 (1) les caract\`eres $\Theta_{\Sigma'_{+}}$, $\Theta_{\Sigma'_{-}}$, $\Theta_{\Sigma_{+}}$ et $\Theta_{\Sigma_{-}}$ sont stables;
 
 (2) $\Theta_{\Sigma'}$ est le transfert de $\Theta_{\Sigma'_{+}}\times \Theta_{\Sigma'_{-}}$ et $\Theta_{\Sigma}$ est le transfert de $\Theta_{\Sigma_{+}}\times \Theta_{\Sigma_{-}}$.
 
 On a d\'efini des termes $S(\Sigma_{+},\Sigma'_{+})$, $S(\Sigma_{+},\Sigma'_{-})$ etc... en 3.2. En prolongeant par bilin\'earit\'e l'application $(\sigma,\sigma')\mapsto m_{geom}(\sigma,\sigma')$ de 2.1, on d\'efinit $m_{geom}(\Sigma,\Sigma')$. On pose
 $$\mu(G')=\left\lbrace\begin{array}{cc}1,&\text{ si }G'\text{ est d\'eploy\'e},\\ -1&,\text{ sinon}.\\ \end{array}\right.$$
 
 \ass{Proposition}{Sous ces hypoth\`eses, on a l'\'egalit\'e
 $$m_{geom}(\Sigma,\Sigma')=\frac{1}{2}(S(\Sigma_{+},\Sigma'_{+})S(\Sigma_{-},\Sigma'_{-})+\mu(G')S(\Sigma_{+},\Sigma'_{-})S(\Sigma_{-},\Sigma'_{+})).$$}
 
 Preuve. On suppose pour fixer les id\'ees $d<d'$. Le membre de gauche de l'\'egalit\'e \`a prouver est une int\'egrale sur un rev\^etement de $\Xi^*(q,q')$. Le membre de droite est la somme de deux int\'egrales, l'une sur $\Xi^*(q_{+},q'_{+})\times \Xi^*(q_{-},q'_{-})$ et l'autre sur $\Xi^*(q_{+},q'_{-})\times \Xi^*(q_{-},q'_{+})$.  Si $G$ et $G'$ sont quasi-d\'eploy\'es, l'application "r\'eunion disjointe" envoie chacun de ces ensembles dans $\Xi^*(q,q')$ et pr\'eserve localement les mesures. Ce n'est plus tout-\`a-fait vrai si $G$ ou $G'$ n'est pas quasi-d\'eploy\'e, \`a cause de la condition (2) de 2.2 qui n'est pas v\'erifi\'ee sur l'image de l'application r\'eunion disjointe. Notons  $\Xi_{et}^*(q,q')$ la r\'eunion de $\Xi^*(q,q')$ et de l'ensemble r\'eduit au terme $\xi=\emptyset$ (c'est-\`a-dire l'ensemble d'indices $I$ est vide). On a $\Xi_{et}^*(q,q')=\Xi^*(q,q')$ si $G$ et $G'$ sont quasi-d\'eploy\'es.  L'application r\'eunion disjointe ci-dessus est toujours \`a valeurs dans cet ensemble. L'\'egalit\'e de l'\'enonc\'e est  de la forme
 $$\int_{\Xi_{et}^*(q,q')}f_{1}(\xi)d\xi=\int_{\Xi_{et}^*(q,q')}f_{2}(\xi)d\xi,$$
 o\`u $f_{1}(\xi)=0$ si $\xi=\emptyset\not\in \Xi^*(q,q')$. 
 Pour la prouver, il suffit de fixer $\xi\in \Xi_{et}^*(q,q')$ et de prouver l'\'egalit\'e $f_{1}(\xi)=f_{2}(\xi)$. Fixons donc $\xi=(I,(F_{\pm i})_{i\in I},(F_{i})_{i\in I},(y_{i})_{i\in I})$, en supposant d'abord $\xi\in \Xi^*(q,q')$. On pose ${\bf d}=d_{\xi}$ et $\boldsymbol{\delta}=\delta_{\xi}$. La fibre du rev\^etement ${\cal C}(q;q')$ au-dessus de $\xi$ est en bijection avec une classe $C^{\epsilon}(\xi)\subset C(\xi)$. Pour $c=(c_{i})_{i\in I^*}\in C^{\epsilon}(\xi)$, notons $x(\xi,c)$ le point de la fibre param\'etr\'e par $c$. D'apr\`es 2.2(8), on a
 $$f_{1}(\xi)=\sum_{c\in C^{\epsilon}(\xi)}c_{\Sigma}(x(\xi,c))c_{\Sigma'}(x(\xi,c))D^G(x(\xi,c))\Delta(x(\xi,c))^r,$$
o\`u $r=(d'-d-1)/2$. En utilisant les notations introduites en   1.5, c'est aussi
 $$(1) \qquad f_{1}(\xi)=D^{d}(\xi)\Delta(\xi)^{r}\sum_{c\in C^{\epsilon}(\xi)}c_{\Sigma}(x(\xi,c))c_{\Sigma'}(x(\xi,c)).$$ 
 Pour tout $I'\subset I$, posons $\xi(I')=(I',(F_{\pm i})_{i\in I'},(F_{i})_{i\in I'},(y_{i})_{i\in I'})$. L'application
 $(I_{1},I_{2})\mapsto (\xi(I_{1}),\xi(I_{2}))$ est une bijection de l'ensemble des couples $(I_{1},I_{2})$ tels que $I_{1}\sqcup I_{2}=I$ sur l'ensemble des $(\xi_{1},\xi_{2})\in \Xi_{reg}^*\times \Xi_{reg}^*$ tels que $\xi_{1}\sqcup \xi_{2}=\xi$. Notons ${\cal I}^+$, resp ${\cal I}^-$, l'ensemble des couples $(I_{1},I_{2})$ tels que $I_{1}\sqcup I_{2}=I$ et $(\xi(I_{1}),\xi(I_{2}))\in \Xi^*(q_{+},q'_{+})\times \Xi^*(q_{-},q'_{-})$, resp. $(\xi(I_{1}),\xi(I_{2}))\in \Xi^*(q_{+},q'_{-})\times \Xi^*(q_{-},q'_{+})$. En utilisant la d\'efinition 3.2(1), on a l'\'egalit\'e
 $$(2) \qquad f_{2}(\xi)=2^{{\bf d}-1}(\sum_{(I_{1},I_{2})\in {\cal I}^+}c_{\Sigma_{+}}(\xi(I_{1}))c_{\Sigma'_{+}}(\xi(I_{1}))c_{\Sigma_{-}}(\xi(I_{2}))c_{\Sigma'_{-}}(\xi(I_{2}))$$
 $$D^{inf(d_{+},d'_{+})}(\xi(I_{1}))D^{inf(d_{-},d'_{-})}(\xi(I_{2}))\Delta(\xi(I_{1}))^{r^+_{1}}\Delta(\xi(I_{2}))^{r^-_{2}})+2^{{\bf d}-1}\mu(G')(\sum_{(I_{1},I_{2})\in {\cal I}^-}c_{\Sigma_{+}}(\xi(I_{1}))$$
 $$c_{\Sigma'_{-}}(\xi(I_{1}))c_{\Sigma_{-}}(\xi(I_{2}))c_{\Sigma'_{+}}(\xi(I_{2}))D^{inf(d_{+},d'_{-})}(\xi(I_{1}))D^{inf(d_{-},d'_{+})}(\xi(I_{2}))\Delta(\xi(I_{1}))^{r^-_{1}}\Delta(\xi(I_{2}))^{r^+_{2}}),$$
 o\`u on a pos\'e $ \vert d_{+}-d'_{+}\vert =1+2r^+_{1}$, $\vert d_{-}-d'_{-}\vert =1+2r^-_{2}$, $\vert d_{+}-d'_{-}\vert =1+2r^-_{1}$, $\vert d_{-}-d'_{+}\vert =1+2r^+_{2}$.
 
 Introduisons le sous-espace quadratique $(W_{{\bf d},\boldsymbol{\delta}},q_{{\bf d},\boldsymbol{\delta}})$ de $(V,q)$, cf. 2.2, les d\'ecompositions orthogonales $(V,q)=(W_{{\bf d},\boldsymbol{\delta}},q_{{\bf d},\boldsymbol{\delta}})\oplus (V_{\sharp},q_{\sharp})$, $(V',q')=(W_{{\bf d},\boldsymbol{\delta}},q_{{\bf d},\boldsymbol{\delta}})\oplus (V'_{\sharp},q'_{\sharp})$ et les notations aff\'erentes. Les \'el\'ements $x(\xi,c)$ appartiennent au groupe sp\'ecial orthogonal $G_{{\bf d},\boldsymbol{\delta}}(F)$ et sont en position g\'en\'erale. Le terme $c_{\Sigma'}(x(\xi,c))$ est calcul\'e par la formule 3.1(1). Dans celle-ci, $Y_{\sharp}$ appartient \`a l'alg\`ebre de Lie d'un tore d\'eploy\'e. Sa classe de conjugaison stable est \'egale \`a sa classe de conjugaison.  Elle est param\'etr\'ee par $(J,(F_{\pm j})_{j\in J},(F_{j})_{j\in J},(Y_{j})_{j\in J})$, o\`u $J=\{1,...,(d'-{\bf d}-1)/2\}$ et, pour tout $j\in J$, $F_{\pm j}=F$, $F_{j}=F\oplus F$, $Y_{j}$ est un \'el\'ement de $F_{j}^{\times}$ tel que $Y_{j}+\tau_{j}(Y_{j})=0$   . La classe de conjugaison stable de l'\'el\'ement $x(\xi,c)exp(\lambda Y_{\sharp})$ est param\'etr\'e par la r\'eunion disjointe de $\xi$ et de $\zeta=(J,(F_{\pm j})_{j\in J},(F_{j})_{j\in J},(exp(\lambda Y_{j}))_{j\in J})$. De m\^eme que l'on a d\'efini $\xi(I')$ pour $I'\subset I$, on d\'efinit $\zeta(J')$ pour $J'\subset J$. Les \'el\'ements de $G'_{+}(F)/stconj\times G'_{-}(F)/stconj$ dont l'image par la correspondance endoscopique est la classe de conjugaison stable de $x(\xi,c)exp(\lambda Y_{\sharp})$ sont param\'etr\'es par les couples $(\xi(I_{1})\sqcup \zeta(J_{1}),\xi(I_{2})\sqcup \zeta(J_{2}))$, o\`u $(I_{1},I_{2},J_{1},J_{2})$ parcourent l'ensemble des quadruplets tels que $I_{1}\sqcup I_{2}=I$, $J_{1}\sqcup J_{2}=J$ et
 $$d_{\xi(I_{1})}+2\vert J_{1}\vert +1=d'_{+},\,\,d_{\xi(I_{2})}+2\vert J_{2}\vert +1=d'_{-}.$$
 Assimilons toute fonction d\'efinie sur un ensemble de classes de conjugaison stable \`a une fonction d\'efinie sur l'ensemble de param\`etres correspondant. Gr\^ace \`a 1.7(1), la formule 1.6(1) s'\'ecrit
 $$(3) \qquad \Theta_{\Sigma'}(x(\xi,c)exp(\lambda Y_{\sharp}))D^{d'}(\xi\sqcup \zeta)^{1/2}=\sum_{I_{1},I_{2},J_{1},J_{2}}\Theta_{\Sigma'_{+}}(\xi(I_{1})\sqcup \zeta(J_{1}))D^{d'_{+}}(\xi(I_{1})\sqcup \zeta(J_{1}))^{1/2}$$
 $$\Theta_{\Sigma'_{-}}(\xi(I_{2})\sqcup \zeta(J_{2})) D^{d'_{-}}(\xi(I_{2})\sqcup \zeta(J_{2}))^{1/2}\Delta(\xi(I_{1})\sqcup \zeta(J_{1}),\xi(I_{2})\sqcup \zeta(J_{2}),c),$$
 o\`u l'ensemble de sommation est celui que l'on vient de d\'ecrire. Montrons que
 
 (4) on a l'\'egalit\'e $\Delta(\xi(I_{1})\sqcup \zeta(J_{1}),\xi(I_{2})\sqcup \zeta(J_{2}),c)=\Delta(\xi(I_{1}),\xi(I_{2}),c)$.
 
 D'apr\`es les d\'efinitions de 1.5, il suffit de prouver que, pour tout $i\in I$, le terme
 $$A=(-1)^{d_{\zeta}/2}P'_{\xi\sqcup \zeta}(y_{i})P_{\xi\sqcup \zeta}(-1)y_{i}^{-d_{\zeta}/2}$$
 est le produit de $B=P'_{\xi}(y_{i})P_{\xi}(-1)$ et d'une norme de l'extension $F_{i}/F_{\pm i}$. Pour tout $j\in J$, \'ecrivons $exp(\lambda Y_{j})=(a_{j},a_{j}^{-1})\in F\oplus F=F_{j}$. Alors
 $$A=B\prod_{j\in J}(-y_{i})^{-1}(y_{i}-a_{j})(y_{i}-a_{j}^{-1})(-1-a_{j})(-1-a_{j}^{-1}),$$
 d'o\`u
 $$A=B\prod_{j\in J}(y_{i}-a_{j})(y_{i}^{-1}-a_{j})(-1-a_{j}^{-1})^2.$$
 Le terme index\'e par $j$ dans ce produit est la norme de $(y_{i}-a_{j})(-1-a_{j}^{-1})$, et cela prouve (4).
 
 Faisons maintenant tendre $\lambda$ vers $0$. Remarquons que 
 $$D^{d'}(\xi\sqcup \zeta)^{1/2}=D^{d'}(\xi)^{1/2}D^{G'_{\sharp}}(exp(\lambda Y_{\sharp}))^{1/2}.$$
  D'autre part, $W(G'_{\sharp},T_{\sharp})$ est le groupe de Weyl d'un syst\`eme de racines de type $B_{\vert J\vert }$, notons $w(B_{\vert J\vert} )$ son nombre d'\'el\'ements. Gr\^ace \`a 3.1(1), le membre de gauche de (3) tend vers 
 $$ w(B_{\vert J\vert} ) D^{G'}(\xi)^{1/2}c_{\Sigma'}(x(\xi,c)).$$
De la m\^eme fa\c{c}on, un terme comme 
$$\Theta_{\Sigma'_{+}}(\xi(I_{1})\sqcup \zeta(J_{1}))D^{d'_{+}}(\xi(I_{1})\sqcup \zeta(J_{1}))^{1/2}$$
tend vers
$$w(B_{\vert J_{1}\vert} )D^{d'_{+}}(\xi(I_{1}))^{1/2}c_{\Sigma'_{+}}(\xi(I_{1})).$$
On obtient
$$w(B_{\vert J\vert} ) D^{d'}(\xi)^{1/2}c_{\Sigma'}(x(\xi,c))=\sum_{I_{1},I_{2},J_{1},J_{2}}\Delta(\xi(I_{1}),\xi(I_{2}),c)w(B_{\vert J_{1}\vert} )w(B_{\vert J_{2}\vert} )D^{d'_{+}}(\xi(I_{1}))^{1/2}$$
$$D^{d'_{-}}(\xi(I_{2}))^{1/2}c_{\Sigma'_{+}}(\xi(I_{1}))c_{\Sigma'_{-}}(\xi(I_{2})).$$
On d\'ecompose la somme en une somme sur les couples $(I_{1},I_{2})$ et une somme int\'erieure sur les $(J_{1},J_{2})$. Les couples $(I_{1},I_{2})$ sont soumis aux seules conditions
$$(5) \qquad I_{1}\sqcup I_{2}=I,\,\,d_{\xi(I_{1})}<d'_{+},\,\,d_{\xi(I_{2})}<d'_{-}.$$
Fixons un tel couple et consid\'erons la somme int\'erieure. Posons ${\bf j}_{1}=(d'_{+}-d_{\xi(I_{1})}-1)/2$, ${\bf j}_{2}=(d'_{-}-d_{\xi(I_{2})}-1)/2$. Les couples $(J_{1},J_{2})$ sont soumis aux seules conditions $J_{1}\sqcup J_{2}=J$, $\vert J_{1}\vert ={\bf j}_{1}$ et $\vert J_{2}\vert ={\bf j}_{2}$. Le terme que l'on somme est simplement $w(B_{{\bf j}_{1}})w(B_{{\bf j}_{2}})$. La somme vaut donc
$$\frac{\vert J\vert !}{{\bf j}_{1}!{\bf j}_{2}!}w(B_{{\bf j}_{1}})w(B_{{\bf j}_{2}}),$$
et on montre que ceci n'est autre que $w(B_{\vert J\vert} )$. On obtient finalement l'\'egalit\'e
$$(6)\qquad c_{\Sigma'}(x(\xi,c))= D^{d'}(\xi)^{-1/2}\sum_{I_{1},I_{2}}\Delta(\xi(I_{1}),\xi(I_{2}),c) B'(I_{1},I_{2}),$$
o\`u l'on somme sur les $(I_{1},I_{2})$ v\'erifiant (5) et o\`u l'on a pos\'e
$$B'(I_{1},I_{2})= D^{d'_{+}}(\xi(I_{1}))^{1/2}D^{d'_{-}}(\xi(I_{2}))^{1/2}c_{\Sigma'_{+}}(\xi(I_{1}))c_{\Sigma'_{-}}(\xi(I_{2})).$$

Le terme $c_{\Sigma}(x(\xi,c))$ est calcul\'e par la formule 3.1(2) ou 3.1(3) selon que $dim(V_{\sharp})\leq 2$ ou $dim(V_{\sharp})\geq 4$. Supposons d'abord  $dim(V_{\sharp})\geq4$. Consid\'erons un couple $(F_{1},F_{2})\in {\cal F}_{\sharp}$. Pour $\epsilon=\pm$, la classe de conjugaison de  $exp(\lambda Y_{F_{1},F_{2}}^{\epsilon})$   est param\'etr\'ee par 
$$\zeta=(J,(F_{\pm j})_{j\in J},(F_{j})_{j\in J},(exp(\lambda Y_{j}))_{j\in J})$$
 et $c^{\epsilon}=(c_{j}^{\epsilon})_{j=1,2}$, avec les m\^emes d\'efinitions qu'en 3.1. Donc la classe de conjugaison de $x(\xi,c)exp(\lambda Y_{F_{1},F_{2}}^{\epsilon})$ est param\'etr\'ee par $\xi\sqcup \zeta$ et $c\sqcup c^{\epsilon}$.  Les \'el\'ements de $G_{+}(F)/stconj\times G_{-}(F)/stconj$ dont l'image par la correspondance endoscopique est la classe de conjugaison stable de $x(\xi,c)exp(\lambda Y_{F_{1},F_{2}}^{\epsilon})$ sont param\'etr\'es par les couples $(\xi(I_{1})\sqcup \zeta(J_{1}),\xi(I_{2})\sqcup \zeta(J_{2}))$, o\`u $(I_{1},I_{2},J_{1},J_{2})$ parcourent l'ensemble des quadruplets tels que $I_{1}\sqcup I_{2}=I$, $J_{1}\sqcup J_{2}=J$,
 $$(7)\qquad d_{\xi(I_{1})}+2\vert J_{1}\vert =d_{+},\,\,d_{\xi(I_{2})}+2\vert J_{2}\vert =d_{-}$$
et
$$(8) \qquad \delta_{\xi(I_{1})}\delta_{\zeta(J_{1})}=\delta(q_{+}),\,\,\delta_{\xi(I_{2})}\delta_{\zeta(J_{2})}=\delta(q_{-}).$$
Il faut prendre garde au fait qu'ici le param\'etrage n'est pas bijectif. D'apr\`es la description de 1.7, \`a un quadruplet sont associ\'ees deux classes de conjugaison stable dont l'image est celle requise, sauf dans le cas o\`u l'un des groupes $G_{+}$ ou $G_{-}$ est r\'eduit \`a $\{1\}$, o\`u il n'y a qu'une classe. Ignorons dans un premier temps cette difficult\'e, nous indiquerons plus tard comment la r\'esoudre. Alors, gr\^ace \`a 1.6(1) et 1.7(3), le terme index\'e par $(F_{1},F_{2})$ dans le membre de droite de 3.1(3) est \'egal \`a
$$(9) \qquad \frac{1}{2}D^d(\xi)^{-1/2}\sum_{I_{1},I_{2},J_{1},J_{2}}\sum_{\epsilon=\pm}\epsilon  \vert W(G_{\sharp},T_{F_{1},F_{2}}^{\epsilon}\vert^{-1} sgn_{F_{1}/F}(\nu_{0}\eta_{\sharp}) \Theta_{\Sigma_{+}}(\xi(I_{1})\sqcup \zeta(J_{1}))$$
$$D^{d_{+}}(\xi(I_{1})\sqcup \zeta(J_{1}))^{1/2}\Theta_{\Sigma_{-}}(\xi(I_{2})\sqcup \zeta(J_{2})) D^{d_{-}}(\xi(I_{2})\sqcup \zeta(J_{2}))^{1/2}\Delta(\xi(I_{1})\sqcup \zeta(J_{1}),\xi(I_{2})\sqcup \zeta(J_{2}),c\sqcup c^{\epsilon}).$$
Consid\'erons un couple $(J_{1},J_{2})$ intervenant ci-dessus et supposons que le sous-ensemble $\{1,2\}$ de $J$ est inclus dans  $J_{1}$ ou dans  $J_{2}$. On voit alors que $\Delta(\xi(I_{1})\sqcup \zeta(J_{1}),\xi(I_{2})\sqcup \zeta(J_{2}),c\sqcup c^{\epsilon})$ est ind\'ependant de $\epsilon$. Le nombre d'\'el\'ements de $W(G_{\sharp},T_{F_{1},F_{2}}^{\epsilon})$ est aussi ind\'ependant de $\epsilon$, notons-le simplement  $\vert W(G_{\sharp},T_{F_{1},F_{2}})\vert $.   La somme en $\epsilon$ est alors la somme de deux termes oppos\'es, qui est nulle. On peut donc se limiter aux couples $(J_{1},J_{2})$ tels que l'un des termes contienne $1$ et l'autre contienne $2$. Si par exemple $J_{1}$ contient $1$ et $J_{2}$ contient $2$, on a $\delta_{\zeta(J_{1})}=\delta_{1}$ et $\delta_{\zeta(J_{2})}=\delta_{2}$ (rappelons que $F_{j}=F(\sqrt{\delta_{j}})$ pour $j=1,2$). Mais alors (8) impose les valeurs de $\delta_{1}$ et $\delta_{2}$.  On obtient le r\'esultat suivant. Consid\'erons la condition

(10) les deux couples $(\delta_{\xi(I_{1})}\delta(q_{+}),\delta_{\xi(I_{2})}\delta(q_{-}))$ et $(\delta_{1},\delta_{2})$ sont \'egaux \`a permutation pr\`es.

 Si elle n'est pas v\'erifi\'ee, il n'y a pas de couple $(J_{1},J_{2})$ v\'erifiant les conditions requises. Poursuivons le calcul en supposant (10) v\'erifi\'ee. Pour fixer les id\'ees, on suppose que les deux couples de cette relation sont \'egaux, le calcul \'etant similaire si l'on doit en permuter un.  Alors les couples $(J_{1},J_{2})$ autoris\'es  sont ceux pour lesquels 
 
 (11) si $\delta_{1}\not=\delta_{2}$, c'est-\`a-dire si $\delta_{\sharp}\not=1$, alors $1\in J_{1}$ et $2\in J_{2}$;
 
 (12) si $\delta_{\sharp}=1$, alors $1\in J_{1}$ et $2\in J_{2}$ ou $1\in J_{2}$ et $2\in J_{1}$.
 
 Pour un tel couple, montrons que, si $\lambda$ est assez proche de $0$, on a
 $$(13)\qquad \epsilon sgn_{F_{1}/F}(\nu_{0}\eta_{\sharp})\Delta(\xi(I_{1})\sqcup \zeta(J_{1}),\xi(I_{2})\sqcup \zeta(J_{2}),c\sqcup c^{\epsilon})= \Delta(\xi(I_{1}),\xi(I_{2}),c).$$
 
On doit encore montrer que, pour $i\in I$, le terme $A$ d\'efini plus haut est \'egal au produit de $B$ et d'une norme de l'extension $F_{i}/F_{\pm i}$. Le terme $A$ est produit de $B$, de termes similaires \`a ceux trait\'es plus haut, et de deux termes 
 $$(-y_{i})^{-1}(y_{i}-exp(\lambda Y_{j}))(y_{i}-exp(-\lambda Y_{j}))(-1-exp(\lambda Y_{j}))(-1-exp(-\lambda Y_{j}))$$
 pour $j=1,2$. Ceci est un \'el\'ement de $F_{\pm i}$. Quand $\lambda$ est  proche de $0$, il est proche de 
 $$4(-y_{i})^{-1}(y_{i}-1)^2,$$
 qui est la norme de $2(y_{i}-1)$. D'o\`u l'assertion. Cela ne suffit pas car, dans le membre de gauche de (13), il y a deux facteurs suppl\'ementaires: le terme $\epsilon sgn_{F_{1}/F}(\nu_{0}\eta_{\sharp})$ et, en supposant par exemple $2\in J_{2}$, un terme $sgn_{F_{2}/F}(\nu_{0}\delta(q)C_{2})$ correspondant \`a l'\'el\'ement $2\in J_{2}^*$. On doit prouver que le produit de ces deux termes vaut $\epsilon$. Rappelons que
 $$\nu_{0}\delta(q)C_{2}=(-1)^{d/2}\nu_{0}\delta(q)c_{2}^{\epsilon}P'_{\xi\sqcup \zeta}(exp(\lambda Y_{2}))P_{\xi\sqcup\zeta}(-1)exp(\lambda Y_{2})^{1-d/2}(exp(\lambda Y_{2})-1)^{-1}(exp(\lambda Y_{2})+1).$$
  Pour $j\in J\setminus \{1,2\}$, posons $exp(\lambda Y_{j})= (a_{j},a_{j}^{-1})\in F^{\times}\oplus F^{\times}=F_{j}^{\times}$. On a les \'egalit\'es 
 $$P_{\xi\sqcup\zeta}(-1)=P_{\xi}(-1)P_{\zeta}(-1), $$
 $$P'_{\xi\sqcup\zeta}(exp(\lambda Y_{2}))=P_{\xi}(exp(\lambda Y_{2}))P'_{\zeta}(exp(\lambda Y_{2})),$$
 $$P'_{\zeta}(exp(\lambda Y_{2}))=(exp(\lambda Y_{2})-exp(-\lambda Y_{2})(exp(\lambda Y_{2})-exp(\lambda Y_{1}))(exp(\lambda Y_{2})-exp(-\lambda Y_{1}))$$
 $$\prod_{j\in J\setminus \{1,2\}}(exp(\lambda Y_{2})-a_{j})(exp(\lambda Y_{2})-a_{j}^{-1}).$$
 On a donc l'\'egalit\'e
 $$\nu_{0}\delta(q)C_{2}=D_{1}D_{2}D_{3}D_{4}D_{5}$$
 o\`u
 $$D_{1}=\nu_{0}\delta(q)c_{2}^{\epsilon},$$
 $$D_{2}=(exp(\lambda Y_{2})-exp(-\lambda Y_{2}))(exp(\lambda Y_{2})-1)^{-1}(exp(\lambda Y_{2})+1),$$
 $$D_{3}=exp(-\lambda Y_{2})(exp(\lambda Y_{2})-exp(\lambda Y_{1}))(exp(\lambda Y_{2})-exp(-\lambda Y_{1})),$$
 $$D_{4}=(-exp(\lambda Y_{2}))^{-d_{\xi}/2}P_{\xi}(exp(\lambda Y_{2}))P_{\xi}(-1),$$
 $$D_{5}=P_{\zeta}(-1)\prod_{j\in J\setminus\{1,2\}}(-exp(\lambda Y_{2}))^{-1}(exp(\lambda Y_{2})-a_{j})(exp(\lambda Y_{2})-a_{j}^{-1}).$$
     Introduisons  la relation d'\'equivalence dans $\bar{F}^{\times}$: $a\equiv b$ si et seulement si $ab^{-1}$ est une norme de l'extension $F_{2}/F$. Chacun des termes ci-dessus appartient \`a $F^{\times}$. A \'equivalence pr\`es, on peut les remplacer par des termes qui leur sont assez proches quand $\lambda$ est assez proche de $0$. Les racines de $P_{\zeta}$ sont proches de $1$, donc $P_{\zeta}(-1)$ est  \'equivalent \`a $(-2)^{2\vert  J\vert }\equiv 1$. Pour $j\in J\setminus \{1,2\}$, le facteur index\'e par $j$ dans $D_{5}$ est \'egal \`a $a_{j}^{-1}Norm_{F_{2}/F}(exp(\lambda Y_{2})-a_{j})$, qui est \'equivalent \`a $1$ puisque $a_{j}$ est proche de $1$. Donc $D_{5}\equiv 1$. Le terme $D_{4}$ est  \'equivalent \`a $(-1)^{d_{\xi}/2}P_{\xi}(1)P_{\xi}(-1)$, ou encore \`a $(-1)^{d_{\xi}/2}P_{\xi}(1)P_{\xi}(-1)^{-1}$. Ceci est le produit sur les $i\in I$ des termes
 $$(-1)^{[F_{\pm i}:F]}Norm_{F_{i}/F}( Y_{i}),$$
 o\`u l'on a pos\'e $Y_{i}=\frac{1-y_{i}}{1+y_{i}}$. On a $\tau_{i}(Y_{i})=-Y_{i}$. Donc $Norm_{F_{i}/F_{\pm i}}(Y_{i})\in -\delta_{i}F_{\pm i}^{\times,2}$, o\`u l'on rappelle que $F_{i}=F_{\pm i}(\sqrt{\delta_{i}})$. Le terme ci-dessus est donc \'equivalent \`a $Norm_{F_{\pm i}/F}(\delta_{i})$, et $D_{4}$ est \'equivalent \`a $\delta_{\xi}$. Le terme $D_{3}$ est    \'equivalent \`a 
 $$\lambda^2(Y_{2}-Y_{1})(Y_{2}+Y_{1})\equiv -Y_{1}^2(1-Norm_{F_{2}/F}(Y_{2})Norm_{F_{1}/F}(Y_{1})^{-1})$$
 $$\equiv -\delta_{1}(1-Norm_{F_{2}/F}(Y_{2})Norm_{F_{1}/F}(Y_{1})^{-1}).$$
 Le terme $D_{2}$ est   \'equivalent \`a $1$. Par d\'efinition de $\delta_{\sharp}$, on a les \'egalit\'es
  $$\delta(q)=\delta_{\sharp}\delta_{\xi} =\delta_{1}\delta_{2}\delta_{\xi}\equiv -\delta_{1}\delta_{\xi}.$$
  En rassemblant ces calculs, on obtient
  $$\nu_{0}\delta(q)C_{2}\equiv \nu_{0}c_{2}^{\epsilon}(1-Norm_{F_{2}/F}(Y_{2}-Norm_{F_{1}/F}(Y_{1})^{-1}).$$
  En utilisant la d\'efinition de $c_{2}^{\epsilon}$, on obtient $sgn_{F_{2}/F}(\nu_{0}\delta(q)C_{2})=\epsilon sgn_{F_{2}/F}(\nu_{0}\eta_{\sharp})$, qui est \'egal \`a $\epsilon sgn_{F_{1}/F}(\nu_{0}\eta_{\sharp})$, ainsi qu'on l'a remarqu\'e en 3.1. Cela ach\`eve la preuve de (13).
  
 Maintenant, le terme que l'on somme dans (9) ne d\'epend plus de $\epsilon$. La somme en $\epsilon$ revient \`a une multiplication par $2$, qui compense le premier facteur $\frac{1}{2}$. Faisons tendre $\lambda\in F^{\times,2}$ vers $0$. Comme plus haut, un terme comme $$\Theta_{\Sigma_{+}}(\xi(I_{1})\sqcup \zeta(J_{1}))D^{d_{+}}(\xi(I_{1})\sqcup \zeta(J_{1}))^{1/2}$$
tend vers
$$w( J_{1} )D^{d_{+}}(\xi(I_{1}))^{1/2}c_{\Sigma_{+}}(\xi(I_{1})),$$
o\`u $w(J_{1})$ est le nombre d'\'el\'ements du groupe de Weyl d'un groupe sp\'ecial orthogonal "pair" contenant un \'el\'ement param\'etr\'e par $\zeta(J_{1})$. On v\'erifie l'\'egalit\'e $w(J_{1})=w(B_{\vert J_{1}\vert -1})$. Alors la limite quand $\lambda$ tend vers $0$ de l'expression (9) vaut
$$\vert W(G_{\sharp},T_{F_{1},F_{2}})\vert ^{-1}\sum_{I_{1},I_{2},J_{1},J_{2}}\Delta(\xi(I_{1}),\xi(I_{2}),c)w(B_{\vert J_{1}\vert-1} )w(B_{\vert J_{2}\vert-1} )D^{d_{+}}(\xi(I_{1}))^{1/2}$$
$$D^{d_{-}}(\xi(I_{2}))^{1/2}c_{\Sigma_{+}}(\xi(I_{1}))c_{\Sigma_{-}}(\xi(I_{2})).$$
Il convient maintenant de tenir compte de la difficult\'e que l'on a signal\'e plus haut et de corriger cette formule en cons\'equence. A chaque quadruplet $(I_{1},I_{2},J_{1},J_{2})$ sont en fait associ\'ees deux classes de conjugaison stable dans $G_{+,reg}(F)\times G_{-,reg}(F)$ (remarquons que les conditions impos\'ees aux quadruplets assurent que  $J_{1}$ et $J_{2}$ ont au moins un \'el\'ement, donc que $G_{+}$ et $G_{-}$ sont non triviaux si l'ensemble de sommation n'est pas vide). Or les calculs ci-dessus ne d\'ependent pas de la classe choisie. Donc le terme qui nous int\'eresse, \`a savoir la limite quand $\lambda$ tend vers $0$  du terme index\'e par $(F_{1},F_{2})$ dans le membre de droite de 3.1(3), est \'egale au terme ci-dessus multipli\'e par $2$. On d\'ecompose la somme ci-dessus en une somme sur les couples $(I_{1},I_{2})$ et une somme int\'erieure sur les couples $(J_{1},J_{2})$. Fixons $(I_{1},I_{2})$.  La somme int\'erieure est
$$\sum_{J_{1},J_{2}}w(B_{\vert J_{1}\vert-1} )w(B_{\vert J_{2}\vert-1} ).$$
Les couples $(J_{1},J_{2})$ sont soumis aux conditions (7), (11) et (12) (toujours en supposant les deux couples de (10) \'egaux). Posons ${\bf j}_{1}=(d_{+}-d_{\xi(I_{1})})/2$, ${\bf j}_{2}=(d_{-}-d_{\xi(I_{1})})/2$. Si $\delta_{\sharp}\not=1$,  l'application $(J_{1},J_{2})\mapsto ( J_{1}\setminus\{1\}, J_{2}\setminus\{2\})$ est une bijection de notre ensemble de couples sur celui des couples $(J'_{1},J'_{2})$ tels que $\vert J'_{1}\vert ={\bf j}_{1}-1$, $\vert J'_{2}\vert ={\bf j}_{2}-1$ et $J'_{1}\sqcup J'_{2}=J\setminus\{1,2\}$. Dans ce cas, la somme int\'erieure vaut
$$\frac{(\vert J\vert-2) !}{({\bf j}_{1}-1)!({\bf j}_{2}-1)!}w(B_{{\bf j}_{1}-1})w(B_{{\bf j}_{2}-1}),$$
et ceci est \'egal \`a $w(B_{\vert J\vert -2})$. Or on v\'erifie que ce nombre est \'egal \`a $\vert W(G_{\sharp},T_{F_{1},F_{2}})\vert /2$. Si $\delta_{\sharp}=1$, il y a deux fois plus de couples $(J_{1},J_{2})$: la condition (12) nous permet d'intervertir les places de $1$ et $2$. Mais $W(G_{\sharp},T_{F_{1},F_{2}})$ est lui-aussi deux fois plus gros: puisque $F_{1}=F_{2}$, il y a des \'el\'ements qui permutent les deux premi\`eres composantes du tore. Le r\'esultat est donc le m\^eme. En tenant compte de la multiplication par $2$ que l'on a r\'etabli ci-dessus, on obtient que la limite quand $\lambda$ tend vers $0$  du terme index\'e par $(F_{1},F_{2})$ dans le membre de droite de 3.1(3) est \'egale \`a
$$ D^d(\xi)^{-1/2}\sum_{I_{1},I_{2}}\Delta(\xi(I_{1}),\xi(I_{2}),c) D^{d_{+}}(\xi(I_{1}))^{1/2}D^{d_{-}}(\xi(I_{2}))^{1/2}c_{\Sigma_{+}}(\xi(I_{1}))c_{\Sigma_{-}}(\xi(I_{2})).$$
Rappelons \`a quelles conditions est soumis le couple $(I_{1},I_{2})$. Il y a la condition de d\'epart
 $$(14)\qquad I_{1}\sqcup I_{2}=I,$$
  la condition (10) et enfin il doit exister un couple $(J_{1},J_{2})$ intervenant dans le calcul ci-dessus tel que (7) soit v\'erifi\'e. Cette derni\`ere condition est \'equivalente \`a
$$d_{\xi(I_{1})}\leq d_{+}-2,\,\, d_{\xi(I_{2})}\leq d_{-}-2.$$
On peut la remplacer par la r\'eunion des conditions
$$d_{\xi(I_{1})}\leq d_{+},\,\, d_{\xi(I_{2})}\leq d_{-}$$
 et
 
 (15) si $d_{\xi(I_{1})}=d_{+}$, $\delta_{\xi(I_{1})}=\delta(q_{+})$; si $d_{\xi(I_{2})}=d_{-}$, $\delta_{\xi(I_{2})}=\delta(q_{-}).$

En effet, la r\'eunion de cette condition et de (10) interdit les \'egalit\'es  $d_{\xi(I_{1})}=d_{+}$ ou $d_{\xi(I_{2})}=d_{-}$ puisque $\delta_{1}$ et $\delta_{2}$ sont tous deux diff\'erents de $1$.
 
 On doit maintenant traiter le premier terme du membre de droite de 3.1(3). Le calcul est essentiellement le m\^eme que ci-dessus, en rempla\c{c}ant $(\delta_{1},\delta_{2})$ par $(\delta_{\sharp},1)$. La seule diff\'erence notable est que l'on peut avoir des couples $(J_{1},J_{2})$ dont l'un des termes est vide.   Cela se produit pour un couple $(I_{1},I_{2})$ tel que, par exemple $d_{\xi(I_{2})}=d_{-}$. D\'etaillons ce cas. Il n'y a plus qu'un choix pour $(J_{1},J_{2})$, \`a savoir $J_{1}=J$ et $J_{2}=\emptyset$. En supposant que $G_{-}\not=\{1\}$, il n'est plus vrai que les deux classes de conjugaison stable de $G_{+,reg}(F)\times G_{-,reg}(F)$ param\'etr\'ees par $(I_{1},I_{2},J_{1},J_{2})$ donnent la m\^eme contribution. L'analogue de (13) reste vrai et la contribution du quadruplet est
 $$\vert W(G_{\sharp},T_{\sharp})\vert ^{-1}D^d(\xi)^{-1/2}D^{d_{+}}(\xi(I_{1})\sqcup \zeta(J))^{1/2}D^{d_{-}}(\xi(I_{2}))^{1/2}\Delta(\xi(I_{1}),\xi(I_{2}),c)$$
 $$(\Theta_{\Sigma_{+}}(y_{+})\Theta_{\Sigma_{-}}(y_{-})+\Theta_{\Sigma_{+}}(y'_{+})\Theta_{\Sigma_{-}}(y'_{-})),$$
 o\`u $(y_{+},y_{-})$ et $(y'_{+},y'_{-})$ sont  les deux classes de conjugaison stable en question. Comme pr\'ec\'edemment, quand $\lambda$ tend vers $0$, les termes $D^{d_{+}}(\xi(I_{1})\sqcup \zeta(J))^{1/2}\Theta_{\Sigma_{+}}(y_{+})$ et $D^{d_{+}}(\xi(I_{1})\sqcup \zeta(J))^{1/2}\Theta_{\Sigma_{+}}(y'_{+})$ tendent tous deux vers $w(J)D^{d_{+}}(\xi(I_{1}))^{1/2}c_{\Sigma_{+}}(\xi(I_{1}))$. Ici $w(J)=\vert W(G_{\sharp},T_{\sharp})\vert $. La limite de l'expression pr\'ec\'edente est donc
 $$D^d(\xi)^{-1/2}D^{d_{+}}(\xi(I_{1}))^{1/2}D^{d_{-}}(\xi(I_{2}))^{1/2}\Delta(\xi(I_{1}),\xi(I_{2}),c)c_{\Sigma_{+}}(\xi(I_{1}))(\Theta_{\Sigma_{-}}(y_{-})+\Theta_{\Sigma_{-}}(y'_{-})).$$
 Les deux \'el\'ements $y_{-}$ et $y'_{-}$ sont les deux classes de conjugaison stable dans $G_{-,reg}(F)$ param\'etr\'ees par $\xi(I_{2})$. Par d\'efinition, la derni\`ere somme ci-dessus est \'egale \`a $c_{\Sigma_{-}}(\xi(I_{2}))$, cf. 3.2. La contribution du couple $(I_{1},I_{2})$ a donc la m\^eme forme que pr\'ec\'edemment. On laisse un calcul plus d\'etaill\'e au lecteur. Le r\'esultat est le m\^eme que ci-dessus, en rempla\c{c}ant dans la condition (10) le couple $(\delta_{1},\delta_{2})$ par le couple $(\delta_{\sharp},1)$. 
 
 Faisons maintenant la somme des contributions de tous les termes du membre de droite de 3.1(3). Notons que pour tout couple $(I_{1},I_{2})$, il y a un unique couple $(\delta_{1},\delta_{2})$ (en incluant ce dernier couple $(\delta_{\sharp},1)$) tel que (10) soit v\'erifi\'e. On obtient alors la formule
 $$(16)\qquad  c_{\Sigma}(x(\xi,c))=D^d(\xi)^{-1/2}\sum_{I_{1},I_{2}}\Delta(\xi(I_{1}),\xi(I_{2}),c) B(I_{1},I_{2}),$$
  o\`u on somme sur les $(I_{1},I_{2})$ v\'erifiant (14) et (15) et o\`u l'on a pos\'e
  $$B(I_{1},I_{2})=D^{d_{+}}(\xi(I_{1}))^{1/2}D^{d_{-}}(\xi(I_{2}))^{1/2}c_{\Sigma_{+}}(\xi(I_{1}))c_{\Sigma_{-}}(\xi(I_{2}))$$
 
 On a suppos\'e $dim(V_{\sharp})\geq4$. Le cas $dim(V_{\sharp})\leq 4$ est similaire, avec encore des subtilit\'es dues au d\'edoublement des classes. Le r\'esultat est le m\^eme.
 
 Gr\^ace \`a (1), (6) et (16), on a l'\'egalit\'e
 $$f_{1}(\xi)=A(\xi) \sum_{I_{1},I_{2}}\sum_{I'_{1},I'_{2}}\sum_{c\in C(\xi)^{\epsilon}}B(I_{1},I_{2})B'(I'_{1},I'_{2})\Delta(\xi(I_{1}),\xi(I_{2}),c)\Delta(\xi(I'_{1}),\xi(I'_{2}),c),$$
 o\`u l'on somme sur les $(I_{1},I_{2})$ v\'erifiant (14) et (15) et les $(I'_{1},I'_{2})$ v\'erifiant (5) et o\`u l'on a pos\'e
  $$A(\xi)=\Delta(\xi)^{r}D^d(\xi)^{1/2}D^{d'}(\xi)^{-1/2}.$$
  Supposons ${\bf d}\not=0$.  On a une \'egalit\'e
  $$\Delta(\xi(I_{1}),\xi(I_{2}),c)\Delta(\xi(I'_{1}),\xi(I'_{2}),c)=\alpha (\prod_{i\in I_{2}}sgn_{F_{i}/F_{\pm i}}(c_{i}))(\prod_{i\in I'_{2}}sgn_{F_{i}/F_{\pm i}}(c_{i})),$$
  o\`u $\alpha$ est un certain signe ind\'ependant de $c$. Rappelons que $C(\xi)^{\epsilon}$ est une classe dans $C(\xi)$ modulo le sous-groupe $C(\xi)^1$  form\'e des $c=(c_{i})_{i\in I}$ tels que $\prod_{i\in I}sgn_{F_{i}/F_{\pm i}}(c_{i})=1$. La somme sur $c\in C(\xi)^{\epsilon}$ des termes ci-dessus n'est non nulle que si $I_{2}=I'_{2}$ ou $I_{2}\sqcup I'_{2}=I$, autrement dit que si $(I_{1},I_{2})$ et $(I'_{1},I'_{2})$ sont \'egaux \`a l'ordre pr\`es. Si ces deux couples sont \'egaux, on a 
  $$\Delta(\xi(I_{1}),\xi(I_{2}),c)\Delta(\xi(I'_{1}),\xi(I'_{2}),c)=1$$
  pour tout $c$ et la somme vaut $\vert C(\xi)^{\epsilon}\vert =2^{{\bf d}-1}$. Supposons $(I_{1},I_{2})=(I'_{2},I'_{1})$. La formule 1.7(1) relie le terme $\Delta(\xi(I'_{1}),\xi(I'_{2}),c)$ au facteur de transfert  relatif au groupe endoscopique $(G'_{+},G'_{-})$ du groupe $G'$, . On sait que, quand on \'echange les deux termes $G'_{+}$ et $G'_{-}$, le facteur de transfert ne change pas si $G'$ est d\'eploy\'e tandis qu'il est multipli\'e par $-1$ si $G'$ n'est pas d\'eploy\'e. On en d\'eduit
  $$\Delta(\xi(I'_{1}),\xi(I'_{2}),c)=\mu(G')\Delta(\xi(I'_{2}),\xi(I'_{1}),c)$$
  pour tout $c\in C_{\xi}^{\epsilon}$. La somme qui nous int\'eresse vaut alors $\mu(G')2^{{\bf d}-1}$. Donc
  $$(17) \qquad f_{1}(\xi)=A(\xi)2^{{\bf d}-1}(A_{1}+\mu(G')A_{2}),$$
  o\`u $A_{1}$ est la somme des $B(I_{1},I_{2})B'(I_{1},I_{2})$ sur les couples $(I_{1},I_{2})$ v\'erifiant (5), (14) et (15), tandis que $A_{2}$ est la somme des $B(I_{1},I_{2})B'(I_{2},I_{1})$ sur les couples $(I_{1},I_{2})$ v\'erifiant (14) et (15) et tels que $(I_{2},I_{1})$ v\'erifie (5). On a suppos\'e ${\bf d}\not=0$. Supposons maintenant ${\bf d}=0$, c'est-\`a-dire $\xi=\emptyset$. Sous l'hypoth\`ese $\xi\in \Xi^*(q,q')$, cela implique que $G$ et $G'$ sont quasi-d\'eploy\'es, donc $\mu(G')=1$. Il n'y a plus qu'un couple $(I_{1},I_{2})$ qui intervient: le couple $(\emptyset,\emptyset)$ et on obtient
  $$f_{1}(\xi)=A(\xi)A_{1}=A_{\xi}A_{2}.$$
  La formule (17) reste vraie. Consid\'erons enfin le cas $\xi=\emptyset$ et $G$ ou $G'$ n'est pas quasi-d\'eploy\'e. Par d\'efinition, $f_{1}(\xi)=0$ et $\mu(G')=-1$. Avec les d\'efinitions ci-dessus, on a $A_{1}=A_{2}$ et la formule (17) est encore v\'erifi\'ee.
  
  Comparons les formules (17) et (2). On voit que les ensembles de sommation qui d\'efinissent $A_{1}$ et $A_{2}$ ne sont autres que ${\cal I}^+$ et ${\cal I}^-$. Pour prouver l'\'egalit\'e $f_{1}(\xi)=f_{2}(\xi)$ et
  la proposition, il reste \`a prouver que les termes que l'on somme sont \'egaux. Faisons-le pour les premi\`eres sommes. On fixe donc $(I_{1},I_{2})\in {\cal I}^+$ et on compare les termes index\'es par ce couple dans  les formules (17) et (2). Dans les deux appara\^{\i}t en facteur le produit
  $$2^{d_{\xi}-1}c_{\Sigma_{+}}(\xi(I_{1}))c_{\Sigma'_{+}}(\xi(I_{1}))c_{\Sigma_{-}}(\xi(I_{2}))c_{\Sigma'_{-}}(\xi(I_{2})).$$
  Les termes restants sont 
 $$(18) \qquad \Delta(\xi)^{r}D^d(\xi)^{1/2}D^{d'}(\xi)^{-1/2}D^{d_{+}}(\xi(I_{1}))^{1/2}D^{d'_{+}}(\xi(I_{1}))^{1/2}D^{d_{-}}(\xi(I_{2}))^{1/2}D^{d'_{-}}(\xi(I_{2}))^{1/2}$$
 pour la formule (17) et
 $$(19)\qquad D^{inf(d_{+},d'_{+})}(\xi(I_{1}))D^{inf(d_{-},d'_{-})}(\xi(I_{2}))\Delta(\xi(I_{1}))^{r^+_{1}}\Delta(\xi(I_{2}))^{r^-_{2}}$$
 pour la formule (2). En utilisant 1.5(1), on a les \'egalit\'es
$$\Delta(\xi)^{r}D^d(\xi)^{1/2}D^{d'}(\xi)^{-1/2}=\Delta(\xi)^{-1/2},$$
$$D^{d_{+}}(\xi(I_{1}))^{1/2}D^{d'_{+}}(\xi(I_{1}))^{1/2}=D^{inf(d_{+},d'_{+})}(\xi(I_{1}))\Delta(\xi(I_{1}))^{(max(d_{+},d'_{+})-inf(d_{+},d'_{+}))/2}$$
$$=D^{inf(d_{+},d'_{+})}(\xi(I_{1}))\Delta(\xi(I_{1}))^{r^+_{1}+1/2},$$
et de m\^eme
 $$D^{d_{-}}(\xi(I_{2}))^{1/2}D^{d'_{-}}(\xi(I_{2}))^{1/2}=D^{inf(d_{-},d'_{-})}(\xi(I_{2}))\Delta(\xi(I_{2}))^{r^-_{2}+1/2}.$$
 Donc (18) est le produit de (19) et de $\Delta(\xi)^{-1/2}\Delta(\xi(I_{1}))^{1/2}\Delta(\xi(I_{2}))^{1/2}$. Ce dernier terme vaut $1$. Cela ach\`eve la d\'emonstration. $\square$

\bigskip

\subsection{Transfert de valeurs de facteurs $\epsilon$}

 On consid\`ere deux espaces $U$ et $U'$ comme en 2.3 et deux espaces quadratiques $(V,q)$ et $(V',q')$, dont on note $G$ et $G'$ les groupes sp\'eciaux orthogonaux. On suppose $dim(V)=dim(U)=d$, $dim(V')=dim(U')=d'$ et $G$ et $G'$ quasi-d\'eploy\'es. Le groupe $G$ est un groupe endoscopique de $\tilde{M}$ et $G'$ est un groupe endoscopique de $\tilde{M}'$. On utilise les facteurs de transfert d\'efinis en 1.8.

  On consid\`ere une combinaison lin\'eaire finie
$$\Sigma=\sum_{k}a_{k}\sigma_{k},$$
o\`u les $\sigma_{k}$ appartiennent \`a $Temp(G(F))$. On consid\`ere une combinaison lin\'eaire analogue $\Sigma'$ relative au groupe $G'$. On consid\`ere une combinaison lin\'eaire finie
$$\Pi=\sum_{k}b_{k}\pi_{k},$$
o\`u les $\pi_{k}$ appartiennent \`a $Temp(M(F))$ et sont autoduales. Chaque $\pi_{k}$ se prolonge en une repr\'esentation $\tilde{\pi}_{k}$ de $\tilde{M}(F)$, cf. 2.3. On pose
$$\tilde{\Pi}=\sum_{k}b_{k}\tilde{\pi}_{k}.$$
On consid\`ere une combinaison lin\'eaire analogue $\Pi'$ relative au groupe $M'$, dont on d\'eduit une combinaison lin\'eaire $\tilde{\Pi}'$. Par lin\'earit\'e, on d\'efinit les caract\`eres $\Theta_{\tilde{\Pi}}$ et $\Theta_{\tilde{\Pi}'}$. On suppose:

$\bullet$ les caract\`eres $\Theta_{\Sigma}$ et $\Theta_{\Sigma'}$ sont stables;

$\bullet$ $\Theta_{\tilde{\Pi}}$ est le transfert de $\Theta_{\Sigma}$ et $\Theta_{\tilde{\Pi}'}$ est le transfert de $\Theta_{\Sigma'}$.

Notons $-1$ l'\'el\'ement central de $G(F)$ qui agit sur $V$ par multiplication par $-1$.  Toute repr\'esentation irr\'eductible $\sigma$ de $G(F)$ poss\`ede un caract\`ere central $\omega_{\sigma}$. Posons
$$\Sigma(-1)=\sum_{k}a_{k}\omega_{\sigma_{k}}(-1)\sigma_{k}.$$
Si une distribution localement int\'egrable $D$ sur $G(F)$ est stablement invariante, la distribution $g\mapsto D(-g)$ l'est aussi. Donc $\Theta_{\Sigma(-1)}$ est stable et on peut d\'efinir $S(\Sigma(-1),\Sigma')$.

 En prolongeant par bilin\'earit\'e l'application $(\pi,\pi')\mapsto \epsilon_{geom,\nu_{1}}(\pi,\pi')$ de 2.3, on d\'efinit $\epsilon_{geom,\nu_{1}}(\Pi,\Pi')$. 
 
 \ass{Proposition}{Sous ces hypoth\`eses, on a l'\'egalit\'e
 $$\epsilon_{geom,\nu_{1}}(\Pi,\Pi')= (\delta(q),2\nu_{1})_{F}S(\Sigma(-1),\Sigma').$$}
 
 Preuve. On suppose comme toujours $d<d'$. D'apr\`es 2.4(1), le membre de gauche de l'\'egalit\'e de l'\'enonc\'e est une int\'egrale sur ${\cal X}(d,d')$, qui est un rev\^etement de $\Xi^*(d,d')$. Le membre de droite est une int\'egrale sur $\Xi^*(q,q')$, qui est inclus dans $\Xi^*(d,d')$. On peut donc \'ecrire l'\'egalit\'e de l'\'enonc\'e sous la forme
 $$\int_{\Xi^*(d,d')}f_{1}(\xi)d\xi=\int_{\Xi^*(d,d')}f_{2}(\xi)d\xi.$$
 Il suffit de fixer $\xi\in \Xi^*(d,d')$ et de prouver l'\'egalit\'e $f_{1}(\xi)=f_{2}(\xi)$.
 
 Fixons donc $\xi=(I,(F_{\pm i})_{i\in I},(F_{i})_{i\in I},(y_{i})_{i\in I})\in \Xi^*(d,d')$. Posons ${\bf d}=d_{\xi}$. La fibre de ${\cal X}(d,d')$ au-dessus de $\xi$ est en bijection avec ${\cal M}({\bf d})\times \Gamma_{pair}(\xi)$. Pour un \'el\'ement $(\mu,\gamma)$ de cet ensemble, notons $\tilde{x}(\xi,\mu,\gamma)$ l'\'el\'ement de la fibre param\'etr\'e par $(\mu,\gamma)$. En se rappelant que l'on doit tenir compte du jacobien de l'application de ${\cal X}(d,d')$ vers $\Xi^*(d,d')$, de la   formule  2.4(1) se d\'eduit l'\'egalit\'e 
 $$f_{1}(\xi)=c({\bf d})\vert 2\vert _{F}^{r^2+r+r(d-{\bf d})-{\bf d}/2}\sum_{(\mu,\gamma)\in {\cal M}({\bf d})\times\Gamma_{pair}(\xi)}c_{\tilde{\Pi}}(\tilde{x}(\xi,\mu,\gamma))c_{\tilde{\Pi}'}(\tilde{x}(\xi,\mu,\gamma))$$
 $$D^{\tilde{M}}(\tilde{x}(\xi,\mu,\gamma))\Delta(\tilde{x}(\xi,\mu,\gamma))^r.$$ 
 On a l'\'egalit\'e $\Delta(\tilde{x}(\xi,\mu,\gamma))=\Delta(\xi)$.  
 Il est clair que $D^{\tilde{M}}(\tilde{x}(\xi,\mu,\gamma))$ ne d\'epend que de $d$ et $\xi$. Notons-le $\tilde{D}^d(\xi)$. On obtient
  $$(1) \qquad f_{1}(\xi)=c({\bf d})\vert 2\vert _{F}^{r^2+r+r(d-{\bf d})-{\bf d}/2}\tilde{D}^d(\xi)\Delta(\xi)^r$$
  $$\sum_{(\mu,\gamma)\in {\cal M}({\bf d})\times\Gamma_{pair}(\xi)}c_{\tilde{\Pi}}(\tilde{x}(\xi,\mu,\gamma))c_{\tilde{\Pi}'}(\tilde{x}(\xi,\mu,\gamma)).$$
 La seconde fonction est donn\'ee par
 $$(2) \qquad f_{2}(\xi)=\left\lbrace\begin{array}{cc}2^{{\bf d}}(\delta(q),2\nu_{1})_{F}c_{\Sigma(-1)}(\xi)c_{\Sigma'}(\xi)D^d(\xi)\Delta(\xi)^r,&\text{ si }\xi\in \Xi^*(q,q'),\\ 0,&\text{ sinon. }\\ \end{array}\right.$$

 Soit $(\mu,\gamma)\in {\cal M}({\bf d})\times\Gamma_{pair}(\xi)$. Introduisons l'espace quadratique $(U_{{\bf d}},q_{{\bf d},\mu})$ de 2.4. Notons-le plut\^ot $(V_{\sharp},q_{\sharp})$ et notons $(V'_{\sharp},q'_{\sharp})$ la somme orthogonale  $(V_{\sharp},q_{\sharp})\oplus (Z_{2r+1},q_{2r+1,-\nu_{1}})$. Notons $G_{\sharp}$ et $G'_{\sharp}$ les groupes sp\'eciaux orthogonaux de ces espaces. L'\'el\'ement $\tilde{x}(\xi,\mu,\gamma)$ de $\tilde{M}'(F)$ (ou plus exactement un repr\'esentant de cette classe) a pour commutant connexe le produit d'un tore et de $G'_{\sharp}$. Ce dernier groupe est d\'eploy\'e et $c_{\tilde{\Pi}'}(\tilde{x}(\xi,\mu,\gamma))$ est donn\'e par une formule analogue \`a 3.1(1), c'est-\`a-dire
 $$(3) \qquad c_{\tilde{\Pi}'}(\tilde{x}(\xi,\mu,\gamma))=\vert W(G'_{\sharp},T_{\sharp})\vert ^{-1}lim_{\lambda\to 0}\Theta_{\tilde{\Pi}'}(\tilde{x}(\xi,\mu,\gamma)exp(\lambda Y_{\sharp}))D^{G'_{\sharp}}(\lambda Y_{\sharp})^{1/2}.$$
  En tant qu'\'el\'ement de $\mathfrak{g}'_{\sharp}(F)$, la classe de conjugaison de $Y_{\sharp}$, qui est \'egale \`a sa classe de conjugaison stable, est param\'etr\'ee par $(J,(F_{\pm j})_{j\in J},(F_{j})_{j\in J},(Y_{j})_{j\in J})$, o\`u $J=\{1,...,(d'-{\bf d}-1)/2\}$ et, pour tout $j\in J$, $F_{\pm j}=F$, $F_{j}=F\oplus F$ et $Y_{j}$ est un \'el\'ement de $F_{j}^{\times}$ tel que $Y_{j}+\tau_{j}(Y_{j})=0$. Cela signifie que l'on peut d\'ecomposer $(V'_{\sharp},q'_{\sharp})$ en somme orthogonale 
  $$D\oplus (\oplus_{j\in J} F_{j}),$$
  o\`u $D$ est une droite et chaque $F_{j}$ est un plan hyperbolique, et $Y$ agit par multiplication par $0$ sur $D$ et par $Y_{j}$ sur $F_{j}$. La restriction de $q'_{\sharp}$ \`a $D$ est forc\'ement \'equivalente \`a la forme $(\alpha,\alpha')\mapsto 2\mu \alpha\alpha'$ sur $F$. On en d\'eduit le param\'etrage de $\tilde{x}(\xi,\mu,\gamma)exp(\lambda Y_{\xi})$ dans $\tilde{M}'(F)$. Sa classe de conjugaison stable est param\'etr\'ee par $\xi\sqcup \zeta$, o\`u $\zeta=(J,(F_{\pm j})_{j\in J},(F_{j})_{j\in J},(exp(2\lambda Y_{j}))_{j\in J})$. Sa classe de conjugaison est param\'etr\'ee par l'\'el\'ement suppl\'ementaire $(2\mu,\gamma)$ de $\Gamma_{imp}(\xi\sqcup \zeta)$. D'apr\`es 1.6(2) et 1.8(1),  l'\'egalit\'e (3) se transforme en
  $$c_{\tilde{\Pi}'}(\tilde{x}(\xi,\mu,\gamma))=\vert W(G'_{\sharp},T_{\sharp})\vert ^{-1}lim_{\lambda\to 0}\Delta(\xi\sqcup \zeta,(2\mu,\gamma))d'(\lambda)\Theta_{\Sigma'}(\xi\sqcup \zeta),$$
  o\`u on a pos\'e
  $$d'(\lambda)=D_{0}^{\tilde{M}'}(\tilde{x}(\xi,\mu,\gamma))^{-1/2}D^{G'_{\sharp}}(\lambda Y_{\sharp})^{1/2}D^{d'}(\xi\sqcup \zeta)^{1/2}.$$
  On v\'erifie comme en 3.3(4) l'\'egalit\'e
  $$\Delta(\xi\sqcup \zeta,(2\mu,\gamma))=\Delta(\xi,(2\mu,\gamma)).$$
  En se reportant \`a la d\'efinition de la fonction $D_{0}^{\tilde{M}'}$, cf. 1.6, on voit que
  $$D_{0}^{\tilde{M}'}(\tilde{x})=\vert 2\vert_{F} ^{-(d'+1)/2}D^{\tilde{M}'}(\tilde{x})$$
  pour tout $\tilde{x}\in \tilde{M}'_{reg}(F)$. L'exposant $-(d'+1)/2$ est obtenu ainsi: notons $T$ le commutant de $M'_{\tilde{x}}$ dans $M'$; alors l'exposant est $dim(M'_{\tilde{x}})-dim(T)$. On obtiendrait $-d/2$ si l'on rempla\c{c}ait $\tilde{M}'$ par $\tilde{M}$. Donc
  $$d'(\lambda)=\vert 2\vert _{F}^{(d'+1)/4}\tilde{D}^{d'}(\xi)^{-1/2}D^{d'}(\xi\sqcup\zeta)^{1/2}.$$
  La formule 3.1(1) montre que
  $$\vert W(G'_{\sharp},T_{\sharp})\vert ^{-1}lim_{\lambda\to 0}D^{d'}(\xi\sqcup\zeta)^{1/2}\Theta_{\Sigma'}(\xi\sqcup \zeta)=D^{d'}(\xi)^{1/2}c_{\Sigma'}(\xi).$$ 
  D'o\`u
  $$(4)\qquad c_{\tilde{\Pi}'}(\tilde{x}(\xi,\mu,\gamma))=\vert 2\vert _{F}^{(d'+1)/4}\tilde{D}^{d'}(\xi)^{-1/2}D^{d'}(\xi)^{1/2}\Delta(\xi,(2\mu,\gamma))c_{\Sigma'}(\xi).$$
  
 Calculons maintenant $c_{\tilde{\Pi}}(\tilde{x}(\xi,\mu,\gamma))$ en supposant d'abord $dim(V_{\sharp})\geq4$. Alors ce terme est donn\'e par une formule analogue \`a 3.1(3). Consid\'erons la contribution d'un couple $(F_{1},F_{2})\in {\cal F}_{\sharp}$. Il intervient un terme $\Theta_{\tilde{\Pi}}(\tilde{x}(\xi,\mu,\gamma)exp(\lambda Y_{F_{1},F_{2}}^{\epsilon}))$. Comme ci-dessus, on peut le calculer comme une somme sur les \'el\'ements $y\in G_{reg}(F)/stconj$ qui correspondent \`a la classe de conjugaison stable de $\tilde{x}(\xi,\mu,\gamma)exp(\lambda Y_{F_{1},F_{2}}^{\epsilon})$. Ces $y$ ne d\'ependent pas de $\epsilon$. Le terme que l'on somme est le produit de
 $$\Delta_{G,\tilde{M}}(y,\tilde{x}(\xi,\mu,\gamma)exp(\lambda Y_{F_{1},F_{2}}^{\epsilon}))$$
et d'un terme qui ne d\'epend pas de $\epsilon$. Le point est qu'en se reportant \`a la formule 1.8(2) qui calcule  ce facteur de transfert, on voit que le facteur ci-dessus ne d\'epend pas non plus de $\epsilon$. Puisque $\epsilon$ intervient en facteur dans la formule 3.1(3), la somme sur $\epsilon$ de ces termes est nulle. Il ne reste que la contribution du premier terme de 3.1(3).    En tant qu'\'el\'ement de $\mathfrak{g}_{\sharp}(F)$, la classe de conjugaison stable de $Y_{\sharp}$  est param\'etr\'ee par $(J,(F_{\pm j})_{j\in J},(F_{j})_{j\in J},(Y_{j})_{j\in J})$, o\`u :
 
$\bullet$  $J=\{1,...,(d-{\bf d})/2\}$;
 
  $\bullet$ si $\delta(q_{\sharp})=1$,  $F_{\pm j}=F$ et $F_{j}=F\oplus F$ pour tout $j\in J$; si $\delta(q_{\sharp})\not=1$, il en est ainsi pour tout $j\in J\setminus \{1\}$, tandis que $F_{\pm 1}=F$ et $F_{1}=F(\sqrt{\delta(q_{\sharp})})$;
  
  $\bullet$ pour tout $j\in J$, $Y_{j}$ est un \'el\'ement de $F_{j}^{\times}$ tel que $Y_{j}+\tau_{j}(Y_{j})=0$.
  
    La classe de conjugaison stable de $\tilde{x}(\xi,\mu,\gamma)exp(\lambda Y_{\sharp})$ est param\'etr\'ee par $\xi\sqcup \zeta$, o\`u $\zeta=(J,(F_{\pm j})_{j\in J},(F_{j})_{j\in J},(exp(2\lambda Y_{j}))_{j\in J})$. Si $\delta_{\sharp}=1$, sa classe de conjugaison est param\'etr\'ee par l'\'el\'ement suppl\'ementaire $\gamma\in \Gamma_{pair}(\xi\sqcup \zeta)=\Gamma_{pair}(\xi)$. Poursuivons le calcul en supposant $\delta(q_{\sharp})\not=1$, le cas o\`u ce terme vaut $1$ ne diff\'erant que par les notations. L'espace quadratique $V_{\sharp}$ est somme orthogonale des $F_{j}$ pour $j\in J$. Pour $j\not=1$, $F_{j}$ est un plan hyperbolique. Pour $j=1$, $F_{1}$ est muni d'une forme $(v,v')\mapsto trace_{F_{1}/F}(\tau_{1}(v)v'\eta)$. Mais, par d\'efinition de l'espace quadratique $(U_{{\bf d}},q_{{\bf d},\mu})$, on sait que $V_{\sharp}$ est somme orthogonale de plans hyperboliques et d'un plan $F^2$ muni de la forme
$$((\alpha,\beta),(\alpha',\beta'))\mapsto 2\mu \alpha\alpha'+2\nu_{1} \beta\beta'.$$
On en d\'eduit les \'egalit\'es $\delta(q_{\sharp})=-\nu_{1}\mu$ et $\eta\equiv \mu\,\,mod\,\,Norm_{F_{1}/F}(F_{1}^{\times})$. Notons $\gamma_{\sharp}$ la famille index\'ee par l'unique \'el\'ement $1\in J$, dont l'unique \'el\'ement est $\mu exp(\lambda Y_{1})$. On voit alors que  la classe de conjugaison de $\tilde{x}(\xi,\mu,\gamma)exp(\lambda Y_{\sharp})$ est param\'etr\'ee par le terme suppl\'ementaire $\gamma\sqcup \gamma_{\sharp}\in \Gamma_{pair}(\xi\sqcup\zeta)$.  .

Si $\delta_{\xi\sqcup \zeta}\not=\delta(q)$, il ne correspond \`a $\tilde{x}(\xi,\mu,\gamma)exp(\lambda Y_{\sharp})$ aucun \'el\'ement de $G_{reg}(F)/stconj$. Alors $\Theta_{\tilde{\Pi}}(\tilde{x}(\xi,\mu,\gamma)exp(\lambda Y_{\sharp}))=0$. Supposons $\delta_{\xi\sqcup \zeta}=\delta(q)$. En vertu du calcul de $\delta_{\zeta}=\delta(q_{\sharp})$ ci-dessus, cette condition \'equivaut \`a $-\nu_{1}\mu \delta_{\xi}=\delta(q)$. Dans ce cas, il y a deux \'el\'ements de $G_{reg}( F)/stconj$ qui correspondent \`a $\tilde{x}(\xi,\mu,\gamma)exp(\lambda Y_{\sharp})$. Comme dans le calcul de 3.3, ces deux \'el\'ements ont en fait la m\^eme contribution. Le fait qu'il y en ait deux va simplement multiplier le r\'esultat par $2$. Ces \'el\'ements sont param\'etr\'es par $-(\xi\sqcup \zeta)$. On obtient
$$c_{\tilde{\Pi}}(\tilde{x}(\xi,\mu,\gamma))=2\vert W(G_{\sharp},T_{\sharp})\vert ^{-1}lim_{\lambda\to 0}\Delta(\xi\sqcup \zeta,\gamma\sqcup\gamma_{\sharp})d(\lambda)\Theta_{\Sigma}(-(\xi\sqcup\zeta)),$$
o\`u l'on a pos\'e
$$d(\lambda)=D_{0}^{\tilde{M}}(\tilde{x}(\xi,\mu,\gamma))^{-1/2}D^{G_{\sharp}}(\lambda Y_{\sharp})^{1/2}D^d(-(\xi\sqcup\zeta)).$$
On calcule comme ci-dessus
$$d(\lambda)=\vert 2\vert _{F}^{d/4}\tilde{D}^d(\xi)^{-1/2}D^d(-(\xi\sqcup \zeta))^{1/2}.$$

     Montrons que

(5) supposons $-\nu\mu\delta_{\xi}=\delta(q)$; alors on a l'\'egalit\'e 
$$\Delta(\xi\sqcup\zeta,\gamma\sqcup \gamma_{\sharp})=(\delta(q),2\nu_{1})_{F}\Delta(\xi,\gamma)\prod_{i\in I}sgn_{F_{i}/F_{\pm i}}(-2\mu).$$

Le rapport $\Delta(\xi\sqcup\zeta,\gamma\sqcup \gamma_{\sharp})\Delta(\xi,\gamma)^{-1}$ est le produit de termes index\'es par $I$ et, dans le cas o\`u $\delta(q_{\sharp})\not=1$, d'un terme suppl\'ementaire index\'e par $j=1\in J$.  Faisons le calcul dans le cas o\`u $\delta(q_{\sharp})\not=1$. Soit $i\in I$. Le terme correspondant est $sgn_{F_{i}/F_{\pm i}}(B_{i})$, o\`u
$$B_{i}=P_{\xi}(1)^{-1}P'_{\xi}(y_{i})^{-1}P_{\xi\sqcup \zeta}(1)P'_{\xi\sqcup}(y_{i})y_{i}^{-\vert  J\vert }.$$
Ecrivons $exp(2\lambda Y_{j})=(a_{j},a_{j}^{-1})\in F\oplus F$ pour $j\in J\setminus \{1\}$. Alors $B_{i}$ est le produit des $y_{i}^{-1}(y_{i}-a_{j})(y_{i}-a_{j}^{-1})(1-a_{j})(1-a_{j}^{-1})$ sur ces $j$ et du terme
$$y_{i}^{-1}(y_{i}-exp(2\lambda Y_{1}))(y_{i}-exp(-2\lambda Y_{1}))(1-exp(2\lambda Y_{1}))(1-exp(-2\lambda Y_{1})).$$
Les premiers termes sont \'egaux \`a $Norm_{F_{i}/F_{\pm i}}((y_{i}-a_{j})(1-a_{j}^{-1}))$.  Pour $\lambda$ proche de $0$, le terme restant est proche de
$$-4\lambda^2Y_{1}^2y_{i}^{-1}(y_{i}-1)^2.$$
Puisque $Y_{1}+\tau_{1}(Y_{1})=0$ et $F_{1}=F(\sqrt{\delta(q_{\sharp})})$,  $Y_{1}^2$ est le produit de $\delta(q_{\sharp})$ et d'un \'el\'ement de $F^{\times,2}$. Le terme ci-dessus est alors le produit de $\delta(q_{\sharp})Norm_{F_{i}/F_{\pm i}}(y_{i}-1))$ et d'un \'el\'ement de $F^{\times,2}$. D'o\`u
$$sgn_{F_{i}/F_{\pm i}}(B_{i})=sgn_{F_{i}/F_{\pm i}}(\delta(q_{\sharp})).$$
On se rappelle que $\delta(q_{\sharp})=-\nu_{1}\mu$ et que l'on a fix\'e  un \'el\'ement $\delta_{i}$ de $F_{\pm i}^{\times}$ tel que $F_{i}=F_{\pm i}(\sqrt{\delta_{i}})$. Pour tout $\alpha\in F^{\times}$, on a l'\'egalit\'e 
$$sgn_{F_{i}/F_{\pm i}}(\alpha)=(Norm_{F_{\pm i}/F}(\delta_{i}),\alpha)_{F}$$
cf. [S] p.216, exercice.
Donc
$$sgn_{F_{i}/F_{\pm i}}(B_{i})=(Norm_{F_{\pm i}/F}(\delta_{i}),2\nu_{1})_{F}sgn_{F_{i}/F_{\pm i}}(-2\mu),$$
puis, en se rappelant la d\'efinition $\delta_{\xi}=\prod_{i\in I}Norm_{F_{\pm i}/F}(\delta_{i})$,
$$(6) \qquad \prod_{i\in I}sgn_{F_{i}/F_{\pm i}}(B_{i})=(\delta_{\xi},2\nu_{1})_{F}\prod_{i\in I}sgn_{F_{i}/F_{\pm i}}(-2\mu).$$
Le terme suppl\'ementaire est $sgn_{F_{1}/F}(B_{1})$, o\`u
$$B_{1}=-\mu exp(-\lambda Y_{1})P_{\xi\sqcup\zeta}(1)P'_{\xi\sqcup\zeta}(exp(2\lambda Y_{1}))exp(2\lambda Y_{1})^{1-d_{\xi\sqcup \zeta}/2 }(exp(2\lambda Y_{1})-1).$$
 On peut d\'ecomposer $B_{1}=D_{1}D_{2}D_{3}$, o\`u
 $$D_{1}=-\mu exp(-\lambda Y_{1})(1-exp(2\lambda Y_{1}))(1-exp(-2\lambda Y_{1}))(exp(2\lambda Y_{1})-exp(-2\lambda Y_{1}))(exp(2\lambda Y_{1})-1),$$
 $$D_{2}=P_{\xi}(1)P_{\xi}(exp(2\lambda Y_{1}))exp(2\lambda Y_{1})^{-d_{\xi}/2},$$
 $$D_{3}=\prod_{j\in J\setminus\{1\}}(1-a_{j})(1-a_{j})^{-1}(exp(2\lambda Y_{1})-a_{j})(exp(2\lambda Y_{1}-a_{j}^{-1})exp(2\lambda Y_{1})^{-1}.$$
 On voit comme ci-dessus que $D_{3}$ est une norme de l'extension $F_{1}/F$. Pour $\lambda$ proche de $0$, $D_{2}$ est proche de $P_{\xi}(1)^2$, qui est aussi une norme. Enfin, $D_{1}$ est proche de $2^5\mu\lambda^4Y_{1}^4$, qui est le produit de $2\mu$ et d'une norme. On obtient
 $$sgn_{F_{1}/F}(B_{1})=sgn_{F_{1}/F}(2\mu)=(\delta(q_{\sharp}),2\mu)_{F}.$$
D'apr\`es cette relation et l'\'egalit\'e (6), il suffit,  pour prouver (5), de prouver l'\'egalit\'e
$$(\delta_{\xi},2\nu_{1})_{F}(\delta(q_{\sharp}),2\mu)_{F}=(\delta(q),2\nu_{1})_{F}.$$
Cela r\'esulte des \'egalit\'es $-2\mu=2\nu_{1}\delta(q_{\sharp})$ et $\delta_{\xi}\delta(q_{\sharp})=\delta(q)$. $\square$
 
On a les \'egalit\'es $\Theta_{\Sigma}(-(\xi\sqcup\zeta))=\Theta_{\Sigma(-1)}(\xi\sqcup \zeta)$ et $D^d(-(\xi\sqcup\zeta))=D^d(\xi\sqcup\zeta)$. La formule 3.1(4) montre que
$$\vert W(G_{\sharp},T_{\sharp})\vert ^{-1}lim_{\lambda\to 0}D^d(\xi\sqcup\zeta)^{1/2}\Theta_{\Sigma(-1)}(\xi\sqcup\zeta)=D^d(\xi)^{1/2}c_{\Sigma(-1)}(\xi).$$
On obtient finalement

$\bullet$ si $\mu=-\nu_{1} \delta_{\xi}\delta(q)$,
$$(7) \qquad c_{\tilde{\Pi}}(\tilde{x}(\xi,\mu,\gamma))= 2\vert 2\vert _{F}^{d/4}\tilde{D}^d(\xi)^{-1/2}D^d(\xi)^{1/2} c_{\Sigma(-1)}(\xi)(\delta(q),2\nu_{1})_{F}\Delta(\xi,\gamma)\prod_{i\in I}sgn_{F_{i}/F_{\pm i}}(-2\mu);$$
 
 $\bullet$ si $\mu\not=-\nu_{1} \delta_{\xi}\delta(q)$, $c_{\tilde{\Pi}}(\tilde{x}(\xi,\mu,\gamma))= 0$.

 On a suppos\'e $dim(V_{\sharp})\geq4$. Un calcul similaire vaut si $dim(V_{\sharp})=2$. Si $dim(V_{\sharp})=0$, il y a une diff\'erence. Supposons $d\geq2$. Il y a encore deux \'el\'ements de $G_{reg}(F)/stconj$ param\'etr\'es par $-\xi$ (ici, $\zeta$ dispara\^{\i}t). Ils n'ont plus de raison de donner la m\^eme contribution. Mais, dans ce cas, $c_{\Sigma(-1)}(\xi)$ est justement d\'efini comme la somme des valeurs de $\Theta_{\Sigma(-1)}$ sur ces deux \'el\'ements. On obtient la m\^eme formule, priv\'ee du premier facteur $2$. Si $d=0$, il n'y a plus qu'une classe de conjugaison stable, et obtient le m\^eme r\'esultat. En se reportant \`a la d\'efinition de $c({\bf d})$, on voit que le r\'esultat ci-dessus est g\'en\'eral, \`a condition de remplacer le premier facteur $2$ par $c({\bf d})^{-1}$.
 
 Revenons \`a la formule (1). D'apr\`es le r\'esultat ci-dessus, la somme en $(\mu,\gamma)$ est limit\'ee au sous-ensemble de ces couples pour lesquels $\mu=-\nu_{1}\delta_{\xi}\delta(q)$. Mais $\mu$ appartient \`a ${\cal M}({\bf d})$. Si ${\bf d}<d$, on a ${\cal M}({\bf d})=F^{\times}/F^{\times,2}$ et il n'y a pas de probl\`eme. Par contre, si ${\bf d}=d$, on a ${\cal M}({\bf d})=\{-\nu_{1}\}$, et l'appartenance de $\mu$ \`a cet ensemble impose $\delta_{\xi}=\delta(q)$, ce qui \'equivaut \`a $\xi\in \Xi^*(q,q')$. Si $\xi\not\in \Xi^*(q,q')$, la somme est donc vide et $f_{1}(\xi)=0$. On a aussi $f_{2}(\xi)=0$, d'o\`u l'\'egalit\'e cherch\'ee dans ce cas. Supposons d\'esormais $\xi\in \Xi^*(q,q')$. On v\'erifie sur les d\'efinitions que le facteur $\Delta(\xi,(2\mu,\gamma))$ qui intervient dans l'\'egalit\'e (4) est \'egal \`a
$$ \Delta(\xi,\gamma)\prod_{i\in I}sgn_{F_{i}/F_{\pm i}}(-2\mu).$$
Alors le produit de (4) et de (7) devient ind\'ependant de $(\mu,\gamma)$. Sommer sur ce couple revient \`a multiplier par le nombre d'\'el\'ements de l'ensemble de sommation. On a d\'ej\`a dit que $\mu$ \'etait en fait fix\'e. Ce nombre d'\'el\'ements est donc celui de $\Gamma_{pair}(\xi)$, c'est-\`a-dire $2^{{\bf d}}$. En utilisant (4) et (7) (o\`u l'on se rappelle que le premier facteur $2$ doit \^etre remplac\'e par $c({\bf d})^{-1}$), on obtient la formule suivante. Posons
$$r(\xi)=r^2+r+r(d-{\bf d})-{\bf d}/2+(d'+1)/4+d/4,$$
$$E(\xi)=\vert 2\vert_{F} ^{r(\xi)}\tilde{D}^d(\xi)^{1/2}\tilde{D}^{d'}(\xi)^{-1/2}D^d(\xi)^{-1/2}D^{d'}(\xi)^{1/2}.$$
Alors
$$f_{1}(\xi)=2^{{\bf d}}(\delta(q),2\nu_{1})_{F}E(\xi)c_{\Sigma(-1)}(\xi)c_{\Sigma'}(\xi)D^d(\xi)\Delta(\xi)^r.$$
D'apr\`es (2), pour d\'emontrer l'\'egalit\'e $f_{1}(\xi)=f_{2}(\xi)$, il reste \`a prouver l'\'egalit\'e $E(\xi)=1$. On doit calculer $\tilde{D}^d(\xi)$ et $\tilde{D}^{d'}(\xi)$. Pour cela, on \'ecrit $U=W_{{\bf d}}\oplus U_{{\bf d}}$ comme en 2.4. On repr\'esente un \'el\'ement $\tilde{x}$ de la classe de conjugaison stable de $\tilde{G}_{reg}(F)$ param\'etr\'e par $\xi$ comme la somme d'une forme bilin\'eaire sur $W_{{\bf d}}$ elle-aussi param\'etr\'ee par $\xi$ et d'une forme bilin\'eaire sym\'etrique sur $U_{{\bf d}}$. Notons $M_{1}$ et $M_{2}$ les groupes lin\'eaires $GL(W_{{\bf d}})$ et $GL(U_{{\bf d}})$. On a la d\'ecomposition
$$\mathfrak{m}=\mathfrak{m}_{1}\oplus \mathfrak{r}\oplus\mathfrak{m}_{2},$$
o\`u
$$\mathfrak{r}=(W_{{\bf d}}\otimes_{F}U_{{\bf d}}^*)\oplus (W_{{\bf d}}^*\otimes_{F}U_{{\bf d}}).$$
L'automorphisme $1-\theta_{\tilde{x}}$ respecte cette d\'ecomposition, donc $\tilde{D}^d(\xi)$ est produit de trois termes. Chacun d'eux est la valeur absolue du d\'eterminant de $1-\theta_{\tilde{x}}$ agissant sur un facteur de la d\'ecomposition ci-dessus, quotient\'e par le noyau de cet automorphisme. Le terme correspondant au facteur $\mathfrak{m}_{1}$ n'est autre que $\tilde{D}^{{\bf d}}(\xi)$. On v\'erifie que $\theta_{\tilde{x}}$ agit sans point fixe sur $\mathfrak{r}$ et que le terme correspondant \`a ce facteur est $\Delta(\xi)^{dim(U_{{\bf d}})}=\Delta(\xi)^{d-{\bf d}}$. Sur $\mathfrak{m}_{2}$, $\theta_{\tilde{x}}$ n'a que deux valeurs propres, $1$ et $-1$. L'espace propre pour $1$ est l'alg\`ebre de Lie d'un groupe sp\'ecial orthogonal. Il est de dimension $dim(U_{{\bf d}})(dim(U_{{\bf d}})-1)/2$. L'espace propre pour $-1$ est donc de dimension $dim(U_{{\bf d}})^2-dim(U_{{\bf d}})(dim(U_{{\bf d}})-1)/2=(d-{\bf d})(d-{\bf d}+1)/2$. Chaque valeur propre $-1$ contribue par $\vert 2\vert _{F}$. Le terme correspondant \`a $\mathfrak{m}_{2}$ est donc $\vert 2\vert _{F}^{(d-{\bf d})(d-{\bf d}+1)/2}$, et finalement
$$\tilde{D}^d(\xi)=\vert 2\vert _{F}^{(d-{\bf d})(d-{\bf d}+1)/2}\Delta(\xi)^{d-{\bf d}}\tilde{D}^{{\bf d}}(\xi).$$
Une m\^eme formule vaut pour $\tilde{D}^{d'}(\xi)$. En utilisant ces formules ainsi que 1.5(1), on obtient
$E(\xi)=\vert 2\vert _{F}^{s(\xi)}$, o\`u
$$s(\xi)=r(\xi)+(d-{\bf d})(d-{\bf d}+1)/4-(d'-{\bf d})(d'-{\bf d}+1)/4.$$
Un simple calcul montre que $s(\xi)=0$, donc $E(\xi)=1$, ce qui ach\`eve la d\'emonstration. $\square$

\bigskip
\subsection{Une premi\`ere cons\'equence}

Conservons les espaces et les groupes  du paragraphe pr\'ec\'edent.  

\ass{Corollaire}{ (i) Soient $\Sigma'$ comme en 3.4, $\pi'$ un \'el\'ement autodual de $Temp(M'(F))$ et $b'\in {\mathbb C}^{\times}$. Supposons que $\Theta_{\Sigma'}$ soit stable et que $b'\Theta_{\tilde{\pi}'}$ soit le transfert de $\Theta_{\Sigma'}$. Alors le caract\`ere central $\omega_{\pi'}$ est trivial.

(ii) Soient $\Sigma$ comme en 3.4, $\pi$ un \'el\'ement autodual de $Temp(M(F))$ et $b\in {\mathbb C}^{\times}$. Supposons que $\Theta_{\Sigma}$ soit stable et que $b\Theta_{\tilde{\pi}}$ soit le transfert de $\Theta_{\Sigma}$. Alors le caract\`ere central $\omega_{\pi}$ est  \'egal au caract\`ere $\alpha\mapsto (\delta(q),\alpha)_{F}$.

(iii)  Soient $\Sigma'$, $\pi'$ et $b'$ v\'erifiant les hypoth\`eses de (i) et $\Sigma$, $\pi$ et $b$ v\'erifiant les hypoth\`eses de (ii). On a l'\'egalit\'e 
$$bb'(\delta(q),-1)_{F}^{(d'-1)/2}\epsilon(1/2,\pi\times \pi',\psi_{F})= S(\Sigma(-1),\Sigma').$$}

{\bf Remarque.} La d\'ependance de $\psi_{F}$ du premier terme ci-dessus n'est qu'apparente, car $b$ et $b'$ d\'ependent aussi de $\psi_{F}$. En effet, les hypoth\`eses sont que $b\Theta_{\tilde{\pi}}$ et $b'\Theta_{\tilde{\pi}'}$ sont des transferts de $\Theta_{\Sigma}$ et $\Theta_{\Sigma'}$. Mais les normalisations  que l'on a choisies de $\tilde{\pi}$ et $\tilde{\pi}'$ d\'ependent de $\psi_{F}$.

Preuve. Dans la situation de (iii), la proposition pr\'ec\'edente et le r\'esultat rappel\'e en 2.3 entra\^{\i}nent
$$(1) \qquad bb'\omega_{\pi}((-1)^{(d'-1)/2}2\nu_{1})\omega_{\pi'}((-1)^{1+d/2}2\nu_{1})\epsilon(1/2,\pi\times \pi',\psi_{F})=(\delta(q),2\nu_{1})_{F}S(\Sigma(-1),\Sigma').$$
L'\'el\'ement $\nu_{1}$ est un ingr\'edient de la preuve, mais on peut le choisir quelconque et ni le facteur $\epsilon$, ni le terme $S(\Sigma(-1),\Sigma')$ n'en d\'ependent. L'\'egalit\'e ci-dessus \'etant vraie pour tout $\nu_{1}$, le caract\`ere $\omega_{\pi}\omega_{\pi'}$ est forc\'ement \'egal \`a $\alpha\mapsto (\delta(q),\alpha)_{F}$. 

Dans la situation de (i), on remplace l'espace $V$ par $0$, on compl\`ete les donn\'ees $\Sigma'$, $\pi'$ et $b'$ par $\Sigma$ r\'eduit \`a l'unique repr\'esentation irr\'eductible de $G(F)=\{1\}$, $\pi$ l'unique repr\'esentation irr\'eductible de $M(F)=\{1\}$ et $b=1$. Alors $\omega_{\pi}=1$ et $\delta(q)=1$. La relation que l'on vient de prouver entra\^{\i}ne que $\omega_{\pi'}=1$.

Dans la situation de (ii), on remplace l'espace $V'$ par une droite, on compl\`ete les donn\'ees $\Sigma$, $\pi$ et $b$ par $\Sigma'$ r\'eduit \`a l'unique repr\'esentation irr\'eductible de $G'(F)=\{1\}$, $\pi'$ la repr\'esentation triviale de $M'(F)=F^{\times}$ et $b'=1$. On v\'erifie ais\'ement que ces donn\'ees satisfont les hypoth\`eses requises. De nouveau, la relation ci-dessus entra\^{\i}ne la conclusion de (ii).

  En revenant \`a la situation de (iii), on remplace dans l'\'egalit\'e (1) les caract\`eres par leurs valeurs que l'on vient de calculer et on obtient la relation cherch\'ee. $\square$

\bigskip

\section{Preuve du th\'eor\`eme principal}

\bigskip

\subsection{Repr\'esentations du groupe de Weil-Deligne}

On note $W_{F}$ le groupe de Weil de $\bar{F}/F$ et $W_{DF}$ le groupe de Weil-Deligne, c'est-\`a-dire $W_{DF}=W_{F}\times SL(2,{\mathbb C})$. Pour tout entier $N \geq1$, notons $\Phi_{temp}(GL(N))$ l'ensemble des classes de conjugaison par $GL(N,{\mathbb C})$ d'homomorphismes continus $ \varphi:W_{DF}\to GL(N,{\mathbb C})$ qui v\'erifient les conditions suivantes:

$\bullet$ $\varphi$ est semi-simple;

$\bullet$ la restriction de $\varphi$ \`a $SL(2,{\mathbb C})$ est alg\'ebrique;

$\bullet$ $\varphi$ est temp\'er\'e, c'est-\`a-dire que l'image de $W_{F}$ par $\varphi$ est relativement compacte.

Notons $\Phi_{temp,irr}(GL(N))$ le sous-ensemble des \'el\'ements irr\'eductibles de $\Phi_{temp}(GL(N))$. D'apr\`es la correspondance de Langlands (th\'eor\`eme de Harris-Taylor et Henniart), tout $\varphi\in \Phi_{temp,irr}(GL(N))$ d\'etermine une repr\'esentation admissible $\pi(\varphi)$ de $GL(N,F)$, qui est unitaire et de la s\'erie discr\`ete.

D\'efinissons une involution $\varphi\mapsto \varphi^{\theta}$ dans $\Phi_{temp}(GL(N))$ par $\varphi^{\theta}(w)={^t\varphi(w)}^{-1}$ pour tout $w\in W_{DF}$. On note $\Phi_{temp,N}$ le sous-ensemble des $\varphi\in \Phi_{temp}(GL(N))$ tels que $\varphi$ est conjugu\'e \`a $\varphi^{\theta}$ et $\Phi_{temp,N,irr}$ le sous-ensemble des \'el\'ements irr\'eductibles de $\Phi_{temp,N}$.  On note $\Phi_{temp}(GL)$, $\Phi_{temp}$ etc... la r\'eunion des $\Phi_{temp}(GL(N))$, $\Phi_{temp,N}$ etc... pour $N\geq1$. Pour $\varphi\in \Phi_{temp}(GL)$,on note $N(\varphi)$ l'entier tel que $\varphi$ appartienne \`a $\Phi_{temp}(GL(N))$.

Tout \'el\'ement $\varphi\in \Phi_{temp,N}$ admet une d\'ecomposition 
$$(1) \qquad \varphi=(\oplus_{i\in I}l_{i}\varphi_{i})\oplus (\oplus_{j\in J}l_{j}(\varphi_{j}\oplus \varphi_{j}^{\theta}))$$
v\'erifiant les conditions suivantes. Les ensembles $I$ et $J$ sont finis et disjoints. Pour tout $i\in I$, resp. $j\in J$, $l_{i}$, resp. $l_{j}$, est un entier strictement positif. L'application $i\mapsto \varphi_{i}$ est injective et prend ses valeurs dans $\Phi_{temp,irr}$. L'application $j\mapsto \{\varphi_{j},\varphi_{j}^{\theta}\}$ est injective et prend ses valeurs dans l'ensemble des ensembles de la forme $\{\varphi',(\varphi')^{\theta}\}$, o\`u $\varphi'$ est un \'el\'ement de $\Phi_{temp}(GL)$ qui n'est pas autodual. On a l'\'egalit\'e
$$N=(\sum_{i\in I}l_{i}N(\varphi_{i}))+(\sum_{j\in J}2l_{j}N(\varphi_{j})).$$
 De (1) se d\'eduit une repr\'esentation 
 $$\pi(\varphi)^L=(\otimes_{i\in I}(\pi(\varphi_{i})\otimes...\otimes \pi(\varphi_{i})))\otimes(\otimes_{j\in J}(\pi(\varphi_{j})\otimes...\otimes \pi(\varphi_{j}))\otimes(\pi(\varphi_{j})^{\vee}\otimes...\otimes\pi(\varphi_{j})^{\vee}))$$
 d'un groupe de L\'evi $L(F)$ de $GL(N,F)$. Chaque $\pi(\varphi_{i})$, resp. $\pi(\varphi_{j})$, $\pi(\varphi_{j})^{\vee}$, est r\'ep\'et\'ee $l_{i}$ fois, resp. $l_{j}$ fois. Choisissons un sous-groupe parabolique $P$ de composante de L\'evi $L$ et posons $\pi(\varphi)=Ind_{P}^{GL(N)}(\pi(\varphi)^L)$. C'est une repr\'esentation admissible, irr\'eductible, temp\'er\'ee et autoduale de $GL(N,F)$.
 
Posons $U_{N}={\mathbb C}^N$. Plus concr\`etement, la d\'ecomposition (1) provient d'une d\'ecomposition
$$(2) \qquad U_{N}=(\oplus_{i\in I}U_{l_{i}}\otimes_{{\mathbb C}}U_{N(\varphi_{i})})\oplus\oplus_{j\in J}((U_{l_{j}}\otimes_{{\mathbb C}} U_{N(\varphi_{j})})\oplus (U_{l_{j}}\otimes_{{\mathbb C}}U_{N(\varphi_{j})})).$$
Pour $i\in I$, le groupe $W_{DF}$ agit sur $U_{N(\varphi_{i})}$ par $\varphi_{i}$. Pour $j\in J$, il agit sur la premi\`ere copie de $U_{N(\varphi_{j})}$ par $\varphi_{j}$ et sur la seconde copie par $\varphi_{j}^{\theta}$. Il agit trivialement sur les autres espaces. 

Fixons une forme quadratique non d\'eg\'en\'er\'ee sur ${\mathbb C}^N$. On note $O(N,{\mathbb C})$ et $SO(N,{\mathbb C})$ ses groupes orthogonaux et sp\'eciaux orthogonaux. Si $N$ est pair, fixons une forme symplectique sur ${\mathbb C}^N$ et notons $Sp(N,{\mathbb C})$ son groupe symplectique.  Notons $\Phi_{temp,N}^{orth}$, resp. $\Phi_{temp,N}^{symp}$ si $N$ est pair, l'ensemble des $\varphi\in \Phi_{temp,N}$  dont l'image est incluse dans $O(N,{\mathbb C})$, resp. $Sp(N,{\mathbb C})$ (plus exactement des $\varphi$ qui sont conjugu\'es \`a un \'el\'ement v\'erifiant cette propri\'et\'e). Pour simplifier l'\'ecriture, on pose $\Phi_{temp,N}^{symp}=\emptyset$ si $N$ est impair. On pose $\Phi_{temp,N,irr}^{orth}=\Phi_{temp,N,irr}\cap \Phi_{temp,N}^{orth}$, $\Phi_{temp,N,irr}^{symp}=\Phi_{temp,N,irr}\cap \Phi_{temp,N}^{symp}$.  On v\'erifie que $\Phi_{temp,N,irr}$ est la r\'eunion disjointe de $\Phi_{temp,N,irr}^{orth}$ et de $\Phi_{temp,N,irr}^{symp}$. Soit $\varphi\in \Phi_{temp,N}$, que l'on \'ecrit sous la forme (1). On note $I^{orth}$, resp. $I^{symp}$, le sous-ensemble des $i\in I$ tels que $\varphi_{i}$ appartienne \`a $\Phi_{temp,irr}^{orth}$, resp. $\Phi_{temp,irr}^{symp}$. L'\'el\'ement $\varphi$ appartient \`a $\Phi_{temp,N}^{orth}$, resp. $ \Phi_{temp,N}^{symp}$, si et seulement si les coefficients $l_{i}$ sont pairs pour tout $i\in I^{symp}$, resp. $i\in I^{orth}$. Soit $\varphi\in \Phi_{temp,N}^{orth}$. Par composition avec le d\'eterminant, on obtient un caract\`ere de $W_{DF}$ \`a valeurs dans $\{\pm 1\}$. Il est forc\'ement trivial sur $SL(2,{\mathbb C})$.  
Sa restriction \`a $W_{F}$ s'identifie par la th\'eorie du corps de classes \`a un caract\`ere quadratique de $F^{\times}$. On note $\delta(\varphi)$ l'\'el\'ement de $F^{\times}/F^{\times,2}$ tel que ce caract\`ere soit $a\mapsto (\delta(\varphi),a)_{F}$. 

Il est utile de calculer le commutant $S_{\varphi}$ dans $SO(N,{\mathbb C})$ de l'image d'un \'el\'ement $\varphi\in \Phi_{temp,N}^{orth}$. Consid\'erons la d\'ecomposition (2). L'espace $U_{N}$ est muni d'une forme quadratique. Pour $i\in I$, les espaces $U_{N(\varphi_{i})}$ sont munis d'une forme soit quadratique, soit symplectique. Pour $j\in J$, la premi\`ere copie de $U_{N(\varphi_{j})}$ est en dualit\'e avec la seconde.  De ces donn\'ees s'en d\'eduisent d'autres: pour $i\in I$,  l'espace $U_{l_{i}}$ est muni d'une forme du m\^eme type que celle sur $U_{N(\varphi_{i})}$; pour $j\in J$, la premi\`ere copie de $U_{l_{j}}$ est en dualit\'e avec la seconde. Le commutant dans $O(N,{\mathbb C})$ de l'image de $\varphi$ est alors
$$(\prod_{i\in I^{orth}}O(l_{i},{\mathbb C}))\times(\prod_{i\in I^{symp}}Sp(l_{i},{\mathbb C}))\times (\prod_{j\in J}GL(l_{j},{\mathbb C})).$$
Consid\'erons un \'el\'ement $x=((x_{i})_{i\in I^{orth}},(x_{i})_{i\in I^{symp}},(x_{j})_{j\in J})$ de ce produit. En tenant compte de la fa\c{c}on dont ce produit est plong\'e dans $O(N,{\mathbb C})$, on voit que le d\'eterminant de $x$ agissant dans $U_{N}$ est \'egal \`a
$$\prod_{i\in I^{orth}}det(x_{i})^{N(\varphi_{i})}.$$
Le commutant dans $SO(N,{\mathbb C})$ est donc le sous-groupe des $x$ tels que ce produit vaille $1$. Notons $I^{orth,pair}$, resp.  $I^{orth,imp}$, le sous-ensemble des $i\in I^{orth}$ tels que $N(\varphi_{i})$ est pair, resp. impair. Notons $S_{\varphi}^0$ la composante neutre de $S_{\varphi}$. On obtient que $S_{\varphi}/S_{\varphi}^0$ est le sous-groupe des \'el\'ements $(e_{i})_{i\in I^{orth}}\in ({\mathbb Z}/2{\mathbb Z})^{I^{orth}}$ tels que $\sum_{i\in I^{orth,imp}}e_{i}=0$. Supposons $N$ pair. Le centre du groupe $O(N,{\mathbb C})$  est \'egal \`a $\{\pm 1\} $. Notons $z_{\varphi}$ l'image de $-1$ dans $S_{\varphi}/S_{\varphi}^0$. On a $z_{\varphi}=(l_{i})_{i\in I^{orth}}$.

Supposons $N$ pair. On calcule de m\^eme le commutant  dans $Sp(N,{\mathbb C})$ de l'image d'un \'el\'ement $\varphi\in \Phi_{temp,N}^{symp}$. Avec des notations similaires \`a celles du cas orthogonal, on obtient que $S_{\varphi}/S_{\varphi}^0=({\mathbb Z}/2{\mathbb Z})^{I^{symp}}$ et que $z_{\varphi}=(l_{i})_{i\in I^{symp}}$.
 
\bigskip

\subsection{Conjectures pour les groupes sp\'eciaux orthogonaux impairs}

Consid\'erons  un espace quadratique $(V',q')$ de dimension $d'$ impaire, pour lequel on utilise les notations maintenant habituelles.   On utilise les constructions de 1.7 et 1.8. En particulier, les transferts apparaissant ci-dessous sont relatifs aux facteurs de transfert d\'efinis dans ces paragraphes. 

Consid\'erons l'ensemble des homomorphismes $\varphi':W_{DF}\to Sp(d'-1,{\mathbb C})$ tels que,  par composition avec  l'inclusion de ce dernier groupe dans $GL(d'-1,{\mathbb C})$, on obtienne un \'el\'ement de $\Phi_{temp}(GL(d'-1))$. Notons $\Phi_{temp}(G')$ l'ensemble des classes de conjugaison par $Sp(d'-1,{\mathbb C})$ dans  cet ensemble. L'application qui, \`a $\varphi'\in \Phi_{temp}(G')$, associe sa composition avec l'inclusion dans $GL(d'-1,{\mathbb C})$, est une bijection de $\Phi_{temp}(G')$ sur $\Phi_{temp,d'-1}^{symp}$. On peut identifier ces deux ensembles. On conjecture qu'il existe une partition
$$Temp(G'(F))=\sqcup_{\varphi'\in \Phi_{temp}(G')}\Pi^{G'}(\varphi')$$
v\'erifiant les propri\'et\'es ci-dessous.

Fixons $\varphi'\in \Phi_{temp}(G')$. Notons ${\cal E}^{G'}(\varphi')$ l'ensemble des caract\`eres $\epsilon'$ de $S_{\varphi'}/S_{\varphi'}^0$ tels que $\epsilon'(z_{\varphi'})=1$ si $G'$ est d\'eploy\'e, $\epsilon'(z_{\varphi'})=-1$ sinon. Alors il existe une bijection $\epsilon'\mapsto \sigma'(\varphi',\epsilon')$ de ${\cal E}^{G'}(\varphi')$ sur $\Pi^{G'}(\varphi')$. Pour $s'\in S_{\varphi'}/S_{\varphi'}^0$, on pose
$$(1) \qquad \Theta_{s'}^{G'}(\varphi')=\sum_{ \epsilon'\in {\cal E}^{G'}(\varphi')}\epsilon'(s')\Theta_{\sigma'(\varphi',\epsilon')}.$$

Supposons $G'$ d\'eploy\'e.   Alors $\Theta_{0}^{G'}(\varphi')$ est une distribution stable (on note $0$ l'\'el\'ement neutre de $S_{\varphi}/S_{\varphi}^0$ puisqu'on a adopt\'e une notation additive). Introduisons un espace $U'$ de dimension $d'$ sur $F$ et le groupe tordu $(M',\tilde{M}')$ associ\'e. La repr\'esentation $\varphi$ de $W_{DF}$ est \`a valeurs dans $GL(d'-1,{\mathbb C})$. Notons ${\bf 1}$ la repr\'esentation triviale de $W_{DF}$ de dimension $1$ et posons $\varphi'_{>}=\varphi'\oplus {\bf 1}$. Alors $\varphi'_{>}$ est un \'el\'ement de $\Phi_{temp,d'}$ et on en d\'eduit une repr\'esentation $\pi(\varphi'_{>})$ de $M'(F)$. Cette repr\'esentation est autoduale et se prolonge comme en 2.3 en une repr\'esentation $\tilde{\pi}(\varphi'_{>})$ de $\tilde{M}(F)$. Rappelons que $G'$ est un groupe endoscopique de $(M',\tilde{M}')$, cf. 1.8. Alors il existe $c^{G'}(\varphi')\in {\mathbb C}^{\times}$ tel que $\vert c^{G'}(\varphi')\vert =1$ et que $c^{G'}(\varphi')\Theta_{\tilde{\pi}(\varphi'_{>})}$ soit le transfert de $\Theta^{G'}_{0}(\varphi')$.  

{\bf Remarque.} Introduisons un espace $U$ de dimension $d'-1$ sur $F$ et le groupe tordu $(M,\tilde{M})$ associ\'e. Le groupe $G'$ est aussi un groupe endoscopique de $(M,\tilde{M})$. Le facteur de transfert est trivial dans ce cas. Une propri\'et\'e beaucoup plus caract\'eristique de $\Pi^{G'}(\varphi')$ est que le transfert de $\Theta^{G'}(\varphi')$ \`a $\tilde{M}(F)$ est un multiple de $\Theta_{\tilde{\pi}(\varphi')}$. Mais nous n'utiliserons pas ce cas d'endoscopie tordue.

Revenons au cas o\`u $G'$ est quelconque. Introduisons des espaces quadratiques $(V'_{+},q'_{+})$ et $(V'_{-},q'_{-})$ v\'erifiant les m\^emes conditions qu'en 1.7. Soient $\varphi'_{+}\in \Phi_{temp}(G'_{+})$ et $\varphi'_{-}\in \Phi_{temp}(G'_{-})$. Supposons $\varphi'=\varphi'_{+}\oplus \varphi'_{-}$.    L'espace ${\mathbb C}^{d'-1}$ de $\varphi'$ se d\'ecompose  conform\'ement en somme directe ${\mathbb C}^{d'_{+}-1}\oplus {\mathbb C}^{d'_{-}-1}$ de sous-espaces stables par la repr\'esentation $\varphi'$. L'automorphisme qui agit par l'identit\'e sur le premier espace et par multiplication par $-1$ sur le second est un \'el\'ement de $S_{\varphi'}$. Notons $s'$ son image dans $S_{\varphi'}/S_{\varphi'}^0$. Alors il existe $\gamma^{G'}(\varphi'_{+},\varphi'_{-})\in {\mathbb C}^{\times}$ tel que $\vert \gamma^{G'}(\varphi'_{+},\varphi'_{-})\vert =1$ et que la distribution $\gamma^{G'}(\varphi'_{+},\varphi'_{-})\Theta_{s'}^{G'}(\varphi')$ soit le transfert de $\Theta^{G'_{+}}_{0}(\varphi'_{+})\times \Theta^{G'_{-}}_{0}(\varphi'_{-})$.  

D'apr\`es la premi\`ere remarque de 1.7, les conjectures ci-dessus sont insensibles au remplacement de $q'$ par $\alpha q'$, pour $\alpha\in F^{\times}$.

\bigskip

\subsection{Conjectures pour les groupes sp\'eciaux orthogonaux pairs}

Consid\'erons  un espace quadratique $(V,q)$ de dimension $d$ paire.  On utilise les constructions de 1.7 et 1.8. Les transferts intervenant ci-dessous sont relatifs aux facteurs de transfert d\'efinis dans ces paragraphes.

  Consid\'erons l'ensemble des homomorphismes $\varphi:W_{DF}\to O(d,{\mathbb C})$ tels que $\delta(\varphi)=\delta(q)$ et que,  par composition avec  l'inclusion de $O(d,{\mathbb C})$  dans $GL(d,{\mathbb C})$, on obtienne un \'el\'ement de $\Phi_{temp}(GL(d))$. Notons $\Phi_{temp}(G)$ l'ensemble des classes de conjugaison par $SO(d,{\mathbb C})$ dans  cet ensemble. La composition avec l'inclusion de $O(d,{\mathbb C})$  dans $GL(d,{\mathbb C})$ d\'efinit une application de $\Phi_{temp}(G)$ dans $\Phi_{temp,d}^{orth}$. Son image est bien s\^ur l'ensemble des $\varphi\in \Phi_{temp,d}^{orth}$ tels que $\delta(\varphi)=\delta(q)$. Le point f\^acheux est que l'application n'est pas injective en g\'en\'eral: deux homormophismes $\varphi:W_{DF}\to O(d,{\mathbb C})$ qui ont m\^eme image dans $\Phi_{temp,d}^{orth}$ sont conjugu\'es par un \'el\'ement de $O(d,{\mathbb C})$, mais pas forc\'ement par un \'el\'ement de $SO(d,{\mathbb C})$.  Pour poser des conjectures raisonnables, on doit  consid\'erer $O(d,{\mathbb C})$ comme le $L$-groupe de $G$, ce qui sous-entend des donn\'ees suppl\'ementaires. En particulier, ces donn\'ees permettent d'identifier un sous-tore maximal de $SO(d,{\mathbb C})$ au groupe dual d'un sous-tore maximal de $G$, cette identification \'etant bien d\'efinie modulo l'action du groupe de Weyl de $G$.
  
  On conjecture qu'il existe une partition
$$(1) \qquad Temp(G(F))=\sqcup_{\varphi\in \Phi_{temp}(G)}\Pi^{G}(\varphi)$$
v\'erifiant les propri\'et\'es ci-dessous.

Fixons $\varphi\in \Phi_{temp}(G)$. Notons ${\cal E}^{G}(\varphi)$ l'ensemble des caract\`eres $\epsilon$ de $S_{\varphi}/S_{\varphi}^0$ tels que $\epsilon(z_{\varphi})=1$ si  $(V,q)$ v\'erifie la condition (QD) de 1.7, $\epsilon(z_{\varphi})=-1$ sinon. Alors il existe une bijection $\epsilon\mapsto \sigma(\varphi,\epsilon)$ de ${\cal E}^{G}(\varphi)$ sur $\Pi^{G}(\varphi)$. Pour $s\in S_{\varphi}/S_{\varphi}^0$, on d\'efinit $\Theta^G_{s}(\varphi)$ par une formule similaire \`a 4.2(1).

Supposons   que $(V,q)$ v\'erifie la condition (QD) de 1.7. Alors $\Theta^{G}_{0}(\varphi)$ est une distribution stable. Introduisons un espace $U$ de dimension $d$ sur $F$ et le groupe tordu $(M,\tilde{M})$ associ\'e.    Rappelons que $G$ est un groupe endoscopique de $\tilde{M}$. On dispose de la repr\'esentation $\pi(\varphi)$ de $M(F)$, que l'on prolonge comme en 2.3 en une repr\'esentation $\tilde{\pi}(\varphi)$ de $\tilde{M}(F)$. Alors il existe $c^G(\varphi)\in {\mathbb C}^{\times}$ tel que $\vert c^G(\varphi)\vert =1$ et que $c^G(\varphi)\Theta_{\tilde{\pi}(\varphi)}$ soit le transfert de $\Theta^G_{0}(\varphi)$.

Revenons au cas o\`u $G$ est quelconque. Introduisons des espaces quadratiques $(V_{+},q_{+})$ et $(V_{-},q_{-})$ v\'erifiant les m\^emes conditions qu'en 1.7.  On suppose plus pr\'ecis\'ement que ces espaces v\'erifient la condition (QD) de 1.7. Soient $\varphi_{+}\in \Phi_{temp}(G_{+})$ et $\varphi_{-}\in \Phi_{temp}(G_{-})$. Posons $\varphi=\varphi_{+}\oplus \varphi_{-}$. C'est un \'el\'ement de $\Phi_{temp}(G)$.  Comme dans le paragraphe pr\'ec\'edent, la d\'efinition de $\varphi$ permet de d\'efinir un \'el\'ement $s\in S_{\varphi}/S_{\varphi}^0$. Alors il existe $\gamma^{G}(\varphi_{+},\varphi_{-})\in {\mathbb C}^{\times}$ tel que $\vert \gamma^{G}(\varphi_{+},\varphi_{-})\vert =1$ et que la distribution $\gamma^{G}(\varphi_{+},\varphi_{-})\Theta^G_{s}(\varphi)$
 soit le transfert de $\Theta^{G_{+}}_{0}(\varphi_{+})\times \Theta^{G_{-}}_{0}(\varphi_{-})$.  Comme plus haut, pour d\'efinir ces derni\`eres distributions, on doit consid\'erer $O(d_{+},{\mathbb C})$ et $O(d_{-},{\mathbb C})$ comme les $L$-groupes de $G_{+}$ et $G_{-}$, c'est-\`a-dire fixer des donn\'ees suppl\'ementaires.
 
 {\bf Remarque.} On pourrait rendre les conjectures plus canoniques de la fa\c{c}on suivante. Consid\'erons l'ensemble des couples $(\sigma,L)$, o\`u $\sigma\in Temp(G(F))$ et $L\in \Lambda(V)$, cf. 1.7. Le groupe orthogonal $G^+(F)$ agit diagonalement sur cet ensemble. Notons $\Lambda Temp(G(F))$ l'ensemble des orbites. Si l'on fixe $L\in \Lambda(V)$, l'application qui, \`a $\sigma\in Temp(G(F))$, associe l'orbite de $(\sigma,L)$ est une bijection de $Temp(G(F))$ sur $\Lambda Temp(G(F))$. Notons $\hat{\Lambda}$ l'ensemble des orbites de lagrangiens dans ${\mathbb C}^d$, pour l'action du groupe sp\'ecial orthogonal $SO(d,{\mathbb C})$. Il a deux \'el\'ements. Consid\'erons l'ensemble des couples $(\varphi,\hat{L})$, o\`u $\varphi\in \Phi_{temp}(G)$ et $\hat{L}\in \hat{\Lambda}$. Le groupe $O(d,{\mathbb C})$ agit diagonalement sur cet ensemble. Notons $\Lambda \Phi_{temp}(G)$ l'ensemble des orbites. De nouveau, si l'on fixe $\hat{L}\in \hat{\Lambda}$, l'application qui, \`a $\varphi\in \Phi_{temp}(G)$, associe l'orbite $\{(\varphi,\hat{L})\}$ de $(\varphi,\hat{L})$ est une bijection de $\Phi_{temp}(G)$ sur $\Lambda \Phi_{temp}(G)$.  Il doit exister une partition canonique
 $$\Lambda Temp(G(F))=\sqcup_{\{(\varphi,\hat{L})\}\in \Lambda\Phi_{temp}(G)}\Pi^{G}(\{(\varphi,\hat{L})\})$$
 dont (1) se d\'eduise de la fa\c{c}on suivante. Consid\'erons $O(d,{\mathbb C})$ comme le $L$-groupe de $G$, ce qui sous-entend que l'on fixe des donn\'ees suppl\'ementaires occultes. Celles-ci d\'efinissent une bijection entre $\hat{\Lambda}$ et $\Lambda(V)$. Soit $\hat{L}\in \hat{\Lambda}$ et $L$ son image dans $\Lambda(V)$. En utilisant ces  \'el\'ements, on identifie $\Lambda\Phi_{temp}(G)$ \`a $\Phi_{temp}(G)$ et $\Lambda Temp(G(F))$ \`a $Temp(G(F))$. Alors la partition ci-dessus devient la partition (1). Remarquons que cela ne d\'epend pas du choix de $\hat{L}$. On pourrait traduire de la m\^eme fa\c{c}on le reste des conjectures. On ne d\'eveloppera pas davantage cette voie un peu trop sophistiqu\'ee.

\bigskip

\subsection{Version faible des conjectures pour les groupes sp\'eciaux orthogonaux pairs}
 
 Soit $(V,q)$ comme dans le paragraphe pr\'ec\'edent. Notons $G^+$ son groupe orthogonal. Le groupe $G^+(F)$ agit naturellement dans $Temp(G(F))$. On note $\overline{Temp}(G(F))$ l'ensemble des orbites. Chacune d'elles a au plus deux \'el\'ements.  Pour $\bar{\sigma}\in \overline{Temp}(G(F))$, on d\'efinit le caract\`ere $\Theta_{\bar{\sigma}}$: si $\bar{\sigma}$ est r\'eduit \`a un \'el\'ement $\sigma$, on pose $\Theta_{\bar{\sigma}}=\Theta_{\sigma}$; si $\bar{\sigma}$ est form\'e de deux \'el\'ements $\sigma_{1}$ et $\sigma_{2}$, on pose $\Theta_{\bar{\sigma}}=\frac{1}{2}(\Theta_{\sigma_{1}}+\Theta_{\sigma_{2}})$.
 
 Consid\'erons le m\^eme ensemble d'homomorphismes $\varphi:W_{DF}\to O(d,{\mathbb C})$ que dans le paragraphe pr\'ec\'edent. Notons $\bar{\Phi}_{temp}(G)$ l'ensemble des classes de conjugaison par $O(d,{\mathbb C})$ dans cet ensemble. Cette fois, la composition avec l'inclusion de $O(d,{\mathbb C})$ dans $GL(d,{\mathbb C})$ d\'efinit une injection de $\bar{\Phi}_{temp}(G)$ dans $\Phi_{temp,d}^{orth}$. On conjecture qu'il existe une partition
$$\overline{Temp}(G(F))=\sqcup_{\varphi\in \bar{\Phi}_{temp}(G)}\bar{\Pi}^{G}(\varphi)$$
v\'erifiant des propri\'et\'es similaires \`a celles d\'ecrites au paragraphe pr\'ec\'edent. On ne r\'ecrit pas ces propri\'et\'es, il suffit d'ajouter judicieusement des $\bar{}$ un peu partout.

{\bf Remarque.} Il est parfois commode de consid\'erer les sous-ensembles de $\overline{Temp}(G(F))$ comme des sous-ensembles de $Temp(G(F))$ qui sont invariants (globalement) par l'action de $G^+(F)$. C'est ce que nous ferons si besoin est.

\bigskip

\subsection{Remarques sur les conjectures}

Il est probable que les travaux en cours de Arthur d\'emontreront la conjecture 4.2, du moins si l'on se restreint aux groupes  d\'eploy\'es (cf. [A1] th\'eor\`eme 30.1).  De m\^eme pour les conjectures de 4.4, du moins si l'on se restreint aux  groupes quasi-d\'eploy\'es. Arthur n'a pas encore publi\'e la preuve de son th\'eor\`eme. Nous ignorons bien entendu quel sera son r\'esultat final. Il est possible qu'il inclue le cas des groupes non quasi-d\'eploy\'es. En tout cas, il est clair que les m\'ethodes d'Arthur permettront de traiter ce cas \`a court terme. Il est moins clair qu'elles permettent de prouver les conjectures plus fines de 4.3. On verra.
En fait, les r\'esultats d'Arthur, outre qu'ils ne se limiteront pas au cas temp\'er\'e, seront certainement plus pr\'ecis, c'est-\`a-dire que les constantes $c^{G'}(\varphi)$ etc... que l'on a introduites seront  explicites, pour des normalisations convenables. Ces normalisations n'\'etant  peut-\^etre pas les m\^emes que les n\^otres, on a pr\'ef\'er\'e formuler les conjectures sous une forme plus vague. 

Dans les conjectures 4.3 et 4.4, on a consid\'er\'e que le $L$-groupe d'un groupe sp\'ecial orthogonal pair \'etait un groupe orthogonal plut\^ot que le produit semi-direct d'un groupe sp\'ecial orthogonal et de $W_{F}$. Cette pr\'esentation des conjectures est due \`a Moeglin ([M]). De m\^eme, le fait que la valeur centrale $\epsilon(z_{\varphi})$ des caract\`eres servant aux param\'etrages d\'epend de la forme du groupe se trouve dans [M], ainsi que dans [V] et [A2].

Dans la situation de 4.2, soit $\varphi'\in \Phi_{temp}(G')$, les propri\'et\'es de la correspondance de Langlands pour les groupes lin\'eaires impliquent que les caract\`eres centraux de $\pi(\varphi')$ et $\pi(\varphi'_{>})$ sont  triviaux. Pour le second caract\`ere, cela r\'esulte aussi du corollaire 3.5(i). De m\^eme, dans la situation de 4.3, resp. 4.4, pour $\varphi\in \Phi_{temp}(G)$, resp. $\varphi\in \bar{\Phi}_{temp}(G)$, le caract\`ere central de $\pi(\varphi)$ est le caract\`ere $\lambda\mapsto (\delta(q),\lambda)_{F}$.

\bigskip

\subsection{Remarques sur les constantes}

Consid\'erons la situation de 4.2, soit $\varphi'\in \Phi_{temp}(G')$. Pour tout $\epsilon'\in {\cal E}^{G'}(\varphi')$, on a la formule d'inversion
$$(1) \qquad \Theta_{\sigma'(\varphi',\epsilon')}=\vert S_{\varphi'}/S_{\varphi'}^0\vert ^{-1}\sum_{s'\in S_{\varphi'}/S_{\varphi'}^0}\epsilon'(s')\Theta^{G'}_{s'}(\varphi').$$
Si   les constantes  $\gamma^{G'}(\varphi'_{+},\varphi'_{-})$ qui figurent dans les conjectures sont connues, cela d\'etermine enti\`erement le param\'etrage du paquet $\Pi^{G'}(\varphi')$. Inversement, il est loisible de changer le param\'etrage de la fa\c{c}on suivante. Fixons un caract\`ere $\epsilon'_{0}$ de $S_{\varphi'}/S_{\varphi'}^0$ tel que $\epsilon'_{0}(z_{\varphi'})=1$. D\'efinissons un nouveau param\'etrage en notant $\sigma'(\varphi',\epsilon')$ la repr\'esentation pr\'ec\'edemment not\'ee $\sigma'(\varphi',\epsilon'\epsilon'_{0})$. Les conjectures sont encore v\'erifi\'ees, les constantes $\gamma^{G'}(\varphi'_{+},\varphi'_{-})$ \'etant multipli\'ees par $\epsilon'_{0}(s')$.

 Regardons ce qui se passe quand, dans la derni\`ere partie de la conjecture, on \'echange les r\^oles des couples $(G'_{+},\varphi'_{+})$ et $(G'_{-},\varphi'_{-})$. L'\'el\'ement $s'$ est remplac\'e par $s'z_{\varphi'}$. D'apr\`es la d\'efinition de ${\cal E}^{G'}(\varphi')$, on a $\Theta^{G'}_{s'z_{\varphi'}}(\varphi')=\Theta^{G'}_{s'}(\varphi')$ si $G'$ est d\'eploy\'e, tandis que $\Theta^{G'}_{s'z_{\varphi'}}(\varphi')=-\Theta^{G'}_{s'}(\varphi')$ si $G'$ n'est pas d\'eploy\'e. Les groupes endoscopiques $G'_{+}\times G'_{-}$ et $G'_{-}\times G'_{+}$ sont \'equivalents. Mais les facteurs de transfert, tels qu'on les a normalis\'es, ne sont pas forc\'ement les m\^emes. Permuter les deux groupes ne change pas ce facteur si $G'$ est d\'eploy\'e et le multiplie par $-1$ si $G'$ n'est pas d\'eploy\'e. Cela entra\^{\i}ne que l'on peut imposer aux constantes d'\^etre sym\'etriques, c'est-\`a-dire de v\'erifier l'\'egalit\'e $\gamma^{G'}(\varphi'_{-},\varphi'_{+})=\gamma^{G'}(\varphi'_{+},\varphi'_{-})$.  Des remarques similaires valent pour les conjectures de 4.3 et 4.4.
 
 Il y a quelques cas o\`u on peut d\'eterminer les constantes. Il y a d'abord un cas formel, celui de 4.3 avec $V=0$. On peut poser formellement $c^G(0)=1$. Le cas de 4.2 avec $dim(V')=1$ est moins formel. On a $G'=\{1\}$ mais $\pi(0_{>})$ est la repr\'esentation triviale de $F^{\times}$. On voit n\'eanmoins que dans ce cas, $ c^{G'}(0)=1$. Dans le cas de 4.2, avec $G'$ d\'eploy\'e et $G'_{+}=G'$, $G'_{-}=\{1\}$, le transfert est l'identit\'e, donc $\gamma^{G'}(\varphi',0)=1$. De m\^eme, dans le cas de 4.3,  si $(V,q)$ v\'erifie la condition (QD) de 1.7, on a $\gamma^G(\varphi,0)=1$. 
 
Consid\'erons la situation de 4.3, supposons que $\delta(q)$ n'est pas un  carr\'e et que $(V,q)$ ne v\'erifie pas la condition (QD) de 1.7. L'espace $(\underline{V},\underline{q})$ est \'equivalent \`a $(V,\alpha q)$, pour un \'el\'ement $\alpha\in F^{\times}$. Fixons un tel $\alpha$ et une racine carr\'ee $\sqrt{\alpha}$ dans $\bar{F}$. Identifions  $(\underline{V},\underline{q})$ \`a $(V,\alpha q)$ de sorte que l'isomorphisme $\beta:V\otimes_{F}\bar{F}\to  \underline{V}\otimes_{F}\bar{F}$ fix\'e en 1.7 soit $v\mapsto \sqrt{\alpha}^{-1}v$.   Le groupe $\underline{G}$ s'identifie \`a $G$ et le torseur int\'erieur $\psi_{G}$ devient l'identit\'e.  Les conjectures s'appliquent \`a $G$ comme \`a $\underline{G}$, mais ne disent pas la m\^eme chose. En effet, pour $\varphi\in \Phi_{temp}(G)=\Phi_{temp}(\underline{G})$, le paquet $\Pi^{\underline{G}}(\varphi)$ est param\'etr\'e par les caract\`eres $\epsilon$ tels que $\epsilon(z_{\varphi})=1$, tandis que le paquet $\Pi^G(\varphi)$ est param\'etr\'e par les $\epsilon$ tels que $\epsilon(z_{\varphi})=-1$.  Introduisons le caract\`ere $\epsilon_{\alpha}$ de $S_{\varphi}/S_{\varphi}^0$ d\'efini par
$$\epsilon_{\alpha}((e_{i})_{i\in I})=\prod_{i\in I}(\delta(\varphi_{i}),\alpha)_{F}^{e_{i}}$$
dans les notations de 4.1.  On v\'erifie que  $\epsilon(z_{\varphi})=-1$. Supposons les conjectures v\'erifi\'ees pour le groupe $\underline{G}$, ou plus pr\'ecis\'ement pour l'espace $(\underline{V},\underline{q})$. Soit $\varphi\in \Phi_{temp}(G)$. On dispose du paquet $\Pi^{\underline{G}}(\varphi)$ et d'un param\'etrage, que l'on note  ici $\epsilon\mapsto \sigma^{\underline{G}}(\varphi,\epsilon)$, de ce paquet par le groupe ${\cal E}^{\underline{G}}(\varphi)$. Posons $\Pi^G(\varphi)=\Pi^{\underline{G}}(\varphi)$ et, pour $\epsilon\in {\cal E}^G(\varphi)$, posons $\sigma^G(\varphi,\epsilon)=\sigma^{\underline{G}}(\varphi,\epsilon_{\alpha}\epsilon)$. En utilisant la seconde remarque de 1.7, on voit qu'avec ces d\'efinitions, les conjectures 4.3 sont encore v\'erifi\'ees pour $G$, avec les m\^emes constantes que pour $\underline{G}$. C'est \`a dire $\gamma^G(\varphi_{+},\varphi_{-})=\gamma^{\underline{G}}(\varphi_{+},\varphi_{-})$.  En particulier, d'apr\`es ce que l'on a vu plus haut, on a $\gamma^G(\varphi,0)=1$.

\bigskip

\subsection{Caract\`ere central}

Soient $(V,q)$ comme en 4.4. Le groupe $G(F)$ a pour centre $\{\pm 1\}$. Pour toute repr\'esentation admissible irr\'eductible $\sigma$ de $G(F)$,  $\sigma(-1)$ est une homoth\'etie de rapport $\pm 1$. On note  $\zeta(\sigma)$ ce rapport. Si $\sigma'$ est conjugu\'e de $\sigma$ par un \'el\'ement du groupe orthogonal, on a $\zeta(\sigma')=\zeta(\sigma)$. Cela permet de d\'efinir $\zeta(\bar{\sigma})$ pour $\bar{\sigma}\in \overline{Temp}(G(F))$.

\ass{Lemme}{On admet les conjectures 4.4. Soit $\varphi\in \bar{\Phi}_{temp}(G)$. Alors il existe $\zeta(\varphi)\in \{\pm 1\}$ tel que $ \zeta(\bar{\sigma})=\zeta(\varphi)$ pour tout $\bar{\sigma}\in\bar{\Pi}(\varphi)$. Ce terme $\zeta(\varphi)$ est le m\^eme pour $G$ et $\underline{G}$.}

Preuve. Supposons la premi\`ere assertion du lemme prouv\'ee pour le groupe  $\underline{G}$, ce qui nous fournit un terme $\zeta(\varphi)$ relatif \`a ce groupe. Consid\'erons les distributions 
$$\Theta^G(\varphi)=\sum_{\bar{\sigma}\in \bar{\Pi}^G(\varphi)}\Theta_{\bar{\sigma}}$$
et
$$ \Gamma^G(\varphi)=\sum_{\bar{\sigma}\in \bar{\Pi}^G(\varphi)}\zeta(\bar{\sigma})\Theta_{\bar{\sigma}}$$
sur $G(F)$. On doit prouver que $\Gamma^G(\varphi)=\zeta(\varphi)\Theta^G(\varphi)$.  En consid\'erant ces distributions comme des fonctions localement int\'egrables, on a l'\'egalit\'e $\Gamma^G(\varphi)(g)=\Theta^G(\varphi)(-g)$ pour tout $g\in G(F)$. On introduit les distributions analogues $\Theta^{\underline{G}}(\varphi)$ et $\Gamma^{\underline{G}}(\varphi)$ sur $\underline{G}(F)$, qui v\'erifient une relation analogue. Les conjectures impliquent que $\Theta^G(\varphi)$ est le transfert de $\Theta^{\underline{G}}(\varphi)$. Mais la multiplication par $-1$ commute au transfert. Donc $\Gamma^G(\varphi)$ est le transfert de $\Gamma^{\underline{G}}(\varphi)$. Or, d'apr\`es l'hypoth\`ese que l'on a faite, on a $\Gamma^{\underline{G}}(\varphi)=\zeta(\varphi)\Theta^{\underline{G}}(\varphi)$. Cela entra\^{\i}ne l'\'egalit\'e cherch\'ee.

  Cela nous ram\`ene au cas o\`u $(V,q)$ v\'erifie (QD). Ecrivons $\varphi$ comme en 4.1(1). Posons
$$\varphi_{0}=\oplus_{i\in I; l_{i} \text{ impair}}\varphi_{i},$$
 et $d_{0}=N(\varphi_{0})$. Le nombre $d_{0}$ est pair. Posons $N=(d-d_{0})/2$.  Il y a un \'el\'ement $\varphi_{1}\in \Phi(GL( N)$ tel que $\varphi=\varphi_{1}\oplus \varphi_{0}\oplus \varphi_{1}^{\theta}$. L'espace $V$ poss\`ede un sous-espace isotrope de dimension $N$. En effet, si $V$ est somme de plans hyperboliques, c'est \'evident. Sinon, la condition $\delta(\varphi)=\delta(q)$ impose $N(\varphi_{0})\geq2$ et l'assertion s'ensuit. On peut donc d\'ecomposer $V$ en somme directe $V=X\oplus V_{0}\oplus Y$, o\`u $X$ et $Y$ sont des espaces isotropes de dimension $N$ et $V_{0}$ est l'orthogonal de $X\oplus Y$. Cette d\'ecomposition donne naissance \`a un groupe de L\'evi $L=GL(N)\times G_{0}$ de $G$, o\`u $G_{0}$ est le groupe sp\'ecial orthogonal de $V_{0}$. On fixe un sous-groupe parabolique de $G$ de composante de L\'evi $L$. On a $\varphi_{0}\in \bar{\Phi}_{temp}(G_{0})$ et $\varphi_{1}$ d\'etermine une repr\'esentation $\pi(\varphi_{1})$ de $GL(N,F)$. Les applications de transfert entre groupes sp\'eciaux orthogonaux pairs et groupes lin\'eaires tordus commutent \`a l'induction, pour peu que l'on ait effectu\'e des choix coh\'erents de facteurs de transfert, ce qui est le cas. On peut alors d\'eduire des conjectures, d'une part que le paquet $\bar{\Pi}^{G_{0}}(\varphi_{0})$ est form\'e de repr\'esentations de la s\'erie discr\`ete, d'autre part que le paquet $\bar{\Pi}^G(\varphi)$ est form\'e des sous-repr\'esentations irr\'eductibles  des induites $Ind_{P}^G(\pi(\varphi_{1})\times \sigma_{0})$ pour $\sigma_{0}\in \bar{\Pi}^{G_{0}}(\varphi_{0})$ (et de leurs conjugu\'es par le groupe orthogonal dans le cas o\`u $V_{0}=\{0\}$). Dans une telle sous-repr\'esentation irr\'eductible, l'\'el\'ement central $-1$ agit par $\omega_{\pi(\varphi_{1})}(-1)\zeta(\sigma_{0})$. Donc l'assertion du lemme r\'esulte de la m\^eme assertion pour $G_{0}$ et $\varphi_{0}$.
 
 Cela nous ram\`ene au cas o\`u $\bar{\Pi}^G(\varphi)$ est form\'e de repr\'esentations de la s\'erie discr\`ete.
 Avec les notations ci-dessus, on sait que $\Theta^G(\varphi)$ est une distribution stable. Puisque $\Gamma^G(\varphi)(g)=\Theta_{0}^G(\varphi)(-g)$ pour tout $g\in G(F)$, $\Gamma^G(\varphi)$ est stable elle-aussi. Il nous suffit de prouver l'assertion plus g\'en\'erale suivante: 
 
 (1) on consid\`ere une combinaison lin\'eaire
 $$\Delta=\sum_{\bar{\sigma}\in \bar{\Pi}^G(\varphi)}a_{\bar{\sigma}}\Theta_{\bar{\sigma}};$$
 supposons que $\Delta$ est stable; alors $\Delta$ est proportionnelle \`a $\Theta^G(\varphi)$.
 
 D'apr\`es 4.6(1), $\Delta$ est combinaison lin\'eaire des $\Theta^G_{s}(\varphi)$. Ecrivons
 $$\Delta=\sum_{s\in S(\varphi)/S(\varphi)^0\{0,z_{\varphi}\}}b_{s}\Theta^G_{s}(\varphi).$$
 Fixons un \'el\'ement r\'egulier elliptique $g\in G(F)$. Soient $(g_{k})_{k=1,...,l}$ des repr\'esentants des classes de conjugaison par $G(F)$ dans la classe de conjugaison stable de $g$. Puisque $\Delta$ est stable, on a l'\'egalit\'e
 $$\Delta(g)=l^{-1}\sum_{k=1,...,l}\Delta(g_{k}).$$
 Soit $s\in S(\varphi)/S(\varphi)^0\{0,z_{\varphi}\}$, $s\not=0$. Alors $\Theta^G_{s}(\varphi)$ est le transfert d'une distribution sur un groupe endoscopique diff\'erent de $G$. Ainsi qu'il est bien connu, l'ensemble $\{g_{k}; k=1,...,l\}$ peut \^etre muni d'une structure de groupe ab\'elien. La restriction  \`a ce groupe de toute distribution localement constante sur les \'el\'ements r\'eguliers et  transfert d'une distribution sur un groupe endoscopique diff\'erent de $G$ est combinaison lin\'eaire de caract\`eres non triviaux de ce groupe. Donc
 $$\sum_{k=1,...,l}\Theta^G_{s}(\varphi)(g_{k})=0.$$
 Par contre, on a $\Theta^G_{0}(\varphi)=\Theta^G(\varphi)$ qui est stable, donc
 $$\sum_{k=1,...,l}\Theta^G_{0}(\varphi)(g_{k})=l \Theta^G_{0}(g).$$
 Il en r\'esulte l'\'egalit\'e $\Delta(g)= b_{0}\Theta^G(\varphi)(g)$. Cela est vrai pour tout $g$ r\'egulier elliptique. Mais on sait qu'une combinaison lin\'eaire de caract\`eres de repr\'esentations de la s\'erie discr\`ete qui est nulle sur les \'el\'ements r\'eguliers elliptiques est nulle partout. Cela prouve (1) et le lemme. $\square$
 
 \bigskip
 
 \subsection{D\'etermination des constantes}

\ass{Lemme}{(i) Supposons la conjecture 4.2 v\'erifi\'ee. Alors les constantes $c^{G'}(\varphi')$ sont \'egales \`a $1$. Quitte \`a changer les param\'etrages, on peut supposer que les constantes $\gamma^{G'}(\varphi'_{+},\varphi'_{-})$ sont \'egales \`a $1$ si $G'$ est d\'eploy\'e, \`a $-1$ sinon.

(ii) Supposons la conjecture 4.3, resp. 4.4, v\'erifi\'ee. Alors les constantes $c^G(\varphi)$ sont \'egales \`a $\zeta(\varphi)\epsilon(1/2,\pi(\varphi),\psi_{F})^{-1}$. Quitte \`a changer les param\'etrages, on peut supposer que les constantes $\gamma^{G}(\varphi_{+},\varphi_{-})$ sont \'egales \`a $1$.}

Ce lemme sera d\'emontr\'e en 4.11 et 4.12.   Dans le (ii), on a not\'e $\epsilon(1/2,\pi(\varphi),\psi_{F})$ le facteur $\epsilon$ usuel de Godement-Jacquet. Il d\'epend de $\psi_{F}$, donc $c^G(\varphi)$ \'egalement. Ce n'est pas surprenant puisque la normalisation de la repr\'esentation $\tilde{\pi}(\varphi)$ d\'epend elle-aussi de $\psi_{F}$.

\bigskip

\subsection{Le th\'eor\`eme}

Soient $N$ et $N'$ deux entiers naturels, avec $N'$ pair, et soient $\varphi\in \Phi_{temp,N}$ et $\varphi'\in \Phi_{temp,N'}^{symp}$. On pose
$$E(\varphi,\varphi')=(\delta(\varphi),-1)_{F}^{N'/2}\epsilon(1/2,\pi(\varphi)\times \pi(\varphi'),\psi_{F}).$$
Ce terme est bilin\'eaire en $\varphi$ et $\varphi'$. On a

(1)  $E(\varphi,\varphi')$ appartient \`a $\{\pm 1\}$ et est ind\'ependant de $\psi_{F}$.

Cela r\'esulte de la preuve de [GP1], proposition 10.5. Celle-ci s'appuie  sur  les deux relations bien connues suivantes. Soient $\rho$ et $\rho'$ deux repr\'esentations admissibles irr\'eductibles de groupes lin\'eaires $GL(N,F)$ et $GL(N',F)$. Soit $\alpha\in F^{\times}$, notons $\psi_{F}^{\alpha}$ le caract\`ere $\lambda\mapsto\psi_{F}(\alpha\lambda)$ de $F$. Alors
$$ \epsilon(1/2,\rho\times \rho',\psi_{F}^{\alpha})=\omega_{\rho}(\alpha)^{N'}\omega_{\rho'}(\alpha)^N\epsilon(1/2,\rho\times \rho',\psi_{F});$$
$$(2) \qquad \epsilon(1/2,\rho\times \rho',\psi_{F})\epsilon(1/2,\rho^{\vee}\times(\rho')^{\vee},\psi_{F})=\omega_{\rho}(-1)^{N'}\omega_{\rho'}(-1)^N.$$
Ici $N'$ est pair et $\omega_{\pi(\varphi')}=1$.

On consid\`ere deux espaces quadratiques $(V',q')$ et $(V,q)$ v\'erifiant les conditions de 2.1. On fixe $\varphi'\in \Phi_{temp}(G')$ et $\varphi\in \Phi_{temp}(G)$. On \'ecrit
$$\varphi'=(\oplus_{i'\in I'}l_{i'}\varphi_{i'})\oplus (\oplus_{j'\in J'}l_{j'}(\varphi_{j'}\oplus \varphi_{j'}^{\theta})),$$
$$\varphi=(\oplus_{i\in I}l_{i}\varphi_{i})\oplus (\oplus_{j\in J}l_{j}(\varphi_{j}\oplus \varphi_{j}^{\theta})),$$
avec des ensembles d'indices disjoints. 

Pour $i'\in (I')^{symp}$, on pose
$$\boldsymbol{\epsilon}_{i'}= E(\varphi,\varphi'_{i'}).$$
 Pour $i\in I^{orth}$, on pose
$$\boldsymbol{\epsilon}_{i}= E(\varphi_{i},\varphi').$$
  On d\'efinit un caract\`ere $\boldsymbol{\epsilon}'$ de $S_{\varphi'}/S_{\varphi'}^0=({\mathbb Z}/2{\mathbb Z})^{(I')^{symp}}$ par
$$\boldsymbol{\epsilon'}((e_{i'})_{i'\in (I')^{symp}})=\prod_{i'\in (I')^{symp}}\boldsymbol{\epsilon}_{i'}^{e_{i'}}.$$
On d\'efinit un caract\`ere  $\boldsymbol{\epsilon}$ de $S_{\varphi}/S_{\varphi}^0\subset ({\mathbb Z}/2{\mathbb Z})^{I^{orth}}$ par
$$\boldsymbol{\epsilon}((e_{i})_{i\in I^{orth}})=\prod_{i\in I^{orth}}\boldsymbol{\epsilon}_{i}^{e_{i}}.$$
En utilisant la relation (2), on d\'emontre l'\'egalit\'e
$$\boldsymbol{\epsilon}(z_{\varphi})=\boldsymbol{\epsilon}'(z_{\varphi'})=E(\varphi,\varphi').$$
On a d\'efini $\mu(G')$ en 3.3. Si $\mu(G')= E(\varphi,\varphi')$, le caract\`ere $\boldsymbol{\epsilon}$, resp.$\boldsymbol{\epsilon}'$,  appartient \`a ${\cal E}^G( \varphi)$, resp. ${\cal E}^{G'}(\varphi')$. Sinon, aucun des deux caract\`eres n'appartient \`a l'ensemble  en question.

Pour $\sigma\in Temp(G(F))$ et $\sigma'\in Temp(G'(F))$, on a d\'efini en 2.1 la multiplicit\'e $m(\sigma,\sigma')$. On v\'erifie que si $\sigma_{1}$ est conjugu\'ee \`a $\sigma$ par un \'el\'ement du groupe orthogonal, on a $m(\sigma_{1},\sigma')=m(\sigma,\sigma')$. Cela permet de d\'efinir $m(\bar{\sigma},\sigma')$ pour $\bar{\sigma}\in \overline{Temp}(G(F))$.

Rappelons que l'on a admis en 2.3 quelques r\'esultats issus d'une hypoth\'etique formule des traces locale tordue. Dans l'\'enonc\'e ci-dessous, on suppose que les constantes des conjectures satisfont aux conclusions du lemme 4.8.

\ass{Th\'eor\`eme}{(i) On admet les conjectures de 4.2 et 4.3. Si $\mu(G')\not= E(\varphi,\varphi')$, on a $m(\sigma,\sigma')=0$ pour tout $(\sigma,\sigma')\in \Pi^G(\varphi)\times \Pi^{G'}(\varphi')$. Supposons $\mu(G')= E(\varphi,\varphi')$. Soit $(\epsilon,\epsilon')\in {\cal E}^G(\varphi)\times {\cal E}^{G'}(\varphi')$. Alors on a les \'egalit\'es
$$m(\sigma(\varphi,\epsilon),\sigma(\varphi',\epsilon'))=\left\lbrace\begin{array}{cc}1,&\text{ si }(\epsilon,\epsilon')=(\boldsymbol{\epsilon},\boldsymbol{\epsilon}'),\\ 0,&\text{ sinon.}\\ \end{array}\right.$$

(ii) On admet les conjectures de 4.2 et 4.4. On a le m\^eme r\'esultat qu'en (i), en rempla\c{c}ant $\Pi^G(\varphi)$ par $\bar{\Pi}^G(\varphi)$ et les repr\'esentations $\sigma(\varphi,\epsilon)$ par $\bar{\sigma}(\varphi,\epsilon)$.}

  C'est la conjecture 6.9 de [GP2] restreinte aux repr\'esentations temp\'er\'ees, \`a ceci pr\`es que Gross et Prasad utilisent les facteurs $\epsilon$ de repr\'esentations galoisiennes et les facteurs associ\'es aux paires de repr\'esentations du groupe lin\'eaire. Mais l'\'egalit\'e de ces deux types de facteurs fait partie des r\'esultats de Harris-Taylor et Henniart.

Le th\'eor\`eme sera d\'emontr\'e en 4.11 et 4.12. Les preuves de (i) et (ii) \'etant similaires, on se contentera de prouver la premi\`ere assertion. Pour la fin de l'article, on admet les conjectures 4.2 et 4.3.  

\bigskip

\subsection{Utilisation des r\'esultats ant\'erieurs}

Consid\'erons les donn\'ees du paragraphe pr\'ec\'edent. Pour $s\in S_{\varphi}/S_{\varphi}^0$ et $s'\in S_{\varphi'}/S_{\varphi'}^0$, on pose
$$m(\varphi,s;\varphi',s')=\sum_{(\epsilon,\epsilon')\in {\cal E}^G(\varphi)\times {\cal E}^{G'}(\varphi')}\epsilon(s)\epsilon'(s')m(\sigma(\varphi,\epsilon),\sigma(\varphi',\epsilon')).$$
 Consid\'erons des espaces $(V_{+},q_{+})$, $(V_{-},q_{-})$, $(V'_{+},q'_{+})$ et $(V'_{-},q'_{-})$ v\'erifiant les hypoth\`eses de 3.3. On impose que $(V_{+},q_{+}) $ et $(V_{-},q_{-}) $ v\'erifient la condition (QD) de 1.7. Soient $\varphi_{+}\in \Phi_{temp}(G_{+})$, $\varphi_{-}\in \Phi_{temp}(G_{-})$, $\varphi'_{+}\in \Phi_{temp}(G'_{+})$ et $\varphi'_{-}\in \Phi_{temp}(G'_{-})$. Supposons que $\varphi=\varphi_{+}\oplus \varphi_{-}$ et $\varphi'=\varphi'_{+}\oplus \varphi'_{-}$. Ces donn\'ees endoscopiques d\'eterminent des \'el\'ements $s\in S_{\varphi}/S_{\varphi}^0$ et $s'\in S_{\varphi'}/S_{\varphi'}^0$. On  d\'efinit les combinaisons lin\'eaires suivantes:
  $$\Sigma'_{+}=\sum_{\sigma'\in\Pi^{G'_{+}}(\varphi'_{+})}\sigma',$$
  $$\Sigma'=\gamma^{G'}(\varphi'_{+},\varphi'_{-})\sum_{\epsilon'\in {\cal E}^{G'}(\varphi')}\epsilon'(s')\sigma(\varphi',\epsilon'),$$
  et on d\'efinit de m\^eme $\Sigma'_{-}$, $\Sigma_{+}$, $\Sigma_{-}$ et $\Sigma$. Les conjectures nous disent que les hypoth\`eses de 3.3 sont satisfaites. Appliquons la proposition de ce paragraphe. D'apr\`es 2.1 et les d\'efinitions, le membre de gauche de l'\'egalit\'e de cette proposition vaut 
  $$(1) \qquad \gamma^G(\varphi_{+},\varphi_{-})\gamma^{G'}(\varphi'_{+},\varphi'_{-})m(\varphi,s;\varphi',s').$$
Le membre de gauche contient des termes $S(\Sigma_{+},\Sigma'_{+})$ etc... Consid\'erons la situation de 3.4, o\`u l'on remplace les couples $(V,q)$ et $(V',q')$ par $(V_{+},q_{+})$ et $(V'_{+},q'_{+})$, les donn\'ees $\Sigma$ et $\Sigma'$ par $\Sigma_{+}$ et $\Sigma'_{+}$ et o\`u l'on pose
$$\Pi=c^{G_{+}}(\varphi_{+})\pi(\varphi_{+}),\,\,\Pi'=c^{G'_{+}}(\varphi'_{+})\pi((\varphi'_{+})_{>}).$$
Les conjectures nous disent que les hypoth\`eses de 3.4 sont v\'erifi\'ees. En appliquant le corollaire 3.5(iii), on obtient
$$S(\Sigma_{+}(-1),\Sigma'_{+})=c^{G_{+}}(\varphi_{+})c^{G'_{+}}(\varphi'_{+})(\delta(q_{+}),-1)_{F}^{(d'_{+}-1)/2}\epsilon(1/2,\pi(\varphi_{+})\times \pi((\varphi'_{+})_{>}),\psi_{F}).$$
Rappelons que $(\varphi'_{+})_{>}=\varphi'_{+}\oplus {\bf 1}$. Donc
$$\epsilon(1/2,\pi(\varphi_{+})\times \pi((\varphi'_{+})_{>}),\psi_{F})=\epsilon(1/2,\pi(\varphi_{+}),\psi_{F})\epsilon(1/2,\pi(\varphi_{+})\times \pi(\varphi'_{+}),\psi_{F}).$$
D'autre part, le lemme 4.6 entra\^{\i}ne que $\Sigma_{+}(-1)=\zeta(\varphi_{+})\Sigma_{+}$.  Alors
$$S(\Sigma_{+},\Sigma'_{+})=\zeta(\varphi_{+})c^{G_{+}}(\varphi_{+})c^{G'_{+}}(\varphi'_{+})\epsilon(1/2,\pi(\varphi_{+}),\psi_{F})E(\varphi_{+},\varphi'_{+}).$$
On calcule de m\^eme les autres termes du membre de droite de l'\'egalit\'e de la proposition 3.3. Celui-ci est donc \'egal au produit de
$$(2) \qquad c^{G_{+}}(\varphi_{+})c^{G'_{+}}(\varphi'_{+})c^{G_{-}}(\varphi_{-})c^{G'_{-}}(\varphi'_{-})\zeta(\varphi_{+})\zeta(\varphi_{-})\epsilon(1/2,\pi(\varphi_{+}),\psi_{F})\epsilon(1/2,\pi(\varphi_{-}),\psi_{F})$$
et de
$$ (3) \qquad \frac{1}{2}(E(\varphi_{+},\varphi'_{+})E(\varphi_{-},\varphi'_{-})+\mu(G')E(\varphi_{+},\varphi'_{-})E(\varphi_{-},\varphi'_{+})).$$
Etudions l'expression (3). Ecrivons
$$\varphi'_{+}=(\oplus_{i'\in I'}l_{i',+}\varphi_{i'})\oplus (\oplus_{j'\in J'}l_{j',+}(\varphi_{j'}\oplus \varphi_{j'}^{\theta})),$$
et \'ecrivons de fa\c{c}on similaire $\varphi'_{-}$, $\varphi_{+}$ et $\varphi_{-}$. On a les \'egalit\'es
$$s=(l_{i,-})_{i\in I^{orth}},\,\,s'=(l'_{i',-})_{i'\in (I')^{symp}}.$$
Parce que l'application $E$ est bilin\'eaire et de carr\'e $1$, on a
$$E(\varphi_{+},\varphi'_{-})E(\varphi_{-},\varphi'_{+})= E(\varphi,\varphi'_{-})E(\varphi_{-},\varphi')$$
et 
$$E(\varphi_{+},\varphi'_{+})E(\varphi_{-},\varphi'_{-})= E(\varphi,\varphi')E(\varphi,\varphi'_{-})E(\varphi_{-},\varphi').$$
En utilisant 4.9(2), on v\'erifie que
$$E(\varphi,\varphi'_{-})=\prod_{i'\in (I')^{symp}}E(\varphi,\varphi'_{i'})^{l_{i',-}}=\boldsymbol{\epsilon}'(s'),$$
$$E(\varphi_{-},\varphi')=\prod_{i\in I^{orth}}E(\varphi_{i},\varphi')^{l_{i,-}}=\boldsymbol{\epsilon}(s).$$
On obtient que (3) est \'egal \`a
$$(4)\qquad  \boldsymbol{\epsilon}(s)\boldsymbol{\epsilon}'(s')(E(\varphi,\varphi')+\mu(G'))/2.$$
La proposition 3.3 dit que (1) est \'egal au produit de (2) et de (4).

	\bigskip

\subsection{La preuve dans le cas $G'$ d\'eploy\'e}

Modifions les hypoth\`eses de d\'epart. On consid\`ere seulement $(V',q')$, $\varphi'\in \Phi_{temp}(G')$ et on suppose $G'$ d\'eploy\'e. Posons $V=\{0\}$, $q=0$.  Les conjectures de 4.2 ne d\'ependent d'aucun param\`etre. On peut remplacer la forme  $q'$ par un de ses multiples, cela ne change rien.  Quitte \`a effectuer un tel changement, on peut supposer que les espaces quadratiques $(V',q')$ et $(V,q)$ v\'erifient les hypoth\`eses de 2.1. Posons $\varphi=0$ et choisissons $(V'_{+},q'_{+})=(V',q')$, $\varphi'_{+}=\varphi'$ dans les constructions du paragraphe pr\'ec\'edent. On a $E(0,\varphi')=\mu(G')=1$ et le terme 4.10(4) est \'egal \`a $1$. Donc les termes 4.10(1) et 4.10(2) sont \'egaux. D'apr\`es les remarques de 4.6, le terme 4.10(1) est \'egal \`a $m(0,0;\varphi',0)$. C'est le nombre d'\'el\'ements du paquet $\Pi^{G'}(\varphi')$ qui admettent un mod\`ele de Whittaker. C'est donc un entier naturel.  Le terme 4.10(2) est \'egal \`a $c^{G'}(\varphi')$. C'est un nombre complexe de module $1$. L'\'egalit\'e des deux termes entra\^{\i}ne $c^{G'}(\varphi')=1$. 

La m\^eme \'egalit\'e entra\^{\i}ne qu'il y a un unique \'el\'ement du paquet $\Pi^{G'}(\varphi')$ qui admet un mod\`ele de Whittaker. Comme on l'a dit en 4.6, on peut modifier le param\'etrage de sorte que cette repr\'esentation soit param\'etr\'ee par le caract\`ere trivial. Sous les m\^emes hypoth\`eses concernant $(V',q')$ et $(V,q)$, prenons maintenant $(V'_{+},q'_{+})$, $\varphi'_{+}$, $(V'_{-},q'_{-})$ et $\varphi'_{-}$ quelconques. Avec le param\'etrage que l'on vient de fixer, on a $m(0,0;\varphi',s')=1$ et le terme 4.10(1) vaut $\gamma^{G'}(\varphi'_{+},\varphi'_{-})$. Le terme 4.10(2) vaut $1$ d'apr\`es ce que l'on vient de prouver, appliqu\'e \`a $G'_{+}$ et $G'_{-}$. Le caract\`ere $\boldsymbol{\epsilon}'$ est trivial par construction donc le terme 4.10(4) vaut $1$. On obtient $\gamma^{G'}(\varphi'_{+},\varphi'_{-})=1$.

Consid\'erons maintenant un espace $(V,q)$ et $\varphi\in \Phi_{temp}(G)$. On suppose que $(V,q )$ v\'erifie la condition (QD) de 1.7. On peut alors trouver une droite quadratique $(V',q')$ de sorte que les hypoth\`eses de 2.1 soient satisfaites. Le groupe $G'=\{1\}$ est d\'eploy\'e. On choisit $(V_{+},q_{+})=(V,q)$ et $\varphi_{+}=\varphi$  dans les constructions du paragraphe pr\'ec\'edent. De nouveau, le terme 4.10(4) vaut $1$ et le terme 4.10(1) est le nombre d'\'el\'ements du paquet $\Pi^G(\varphi)$ qui admettent un mod\`ele de Whittaker relativement \`a l'orbite unipotente r\'eguli\`ere param\'etr\'ee par $\nu_{0}$. Le terme 4.10(2) se r\'eduit \`a $c^{G}(\varphi)\zeta(\varphi)\epsilon(1/2,\pi(\varphi),\psi_{F})$, qui est un nombre complexe de module $1$. L'\'egalit\'e des deux termes entra\^{\i}ne $c^G(\varphi)=\zeta(\varphi)\epsilon(1/2,\pi(\varphi),\psi_{F})^{-1}$.  

Il y a encore un unique \'el\'ement de $\Pi^G(\varphi))$ qui admet un mod\`ele de Whittaker du type ci-dessus. On modifie le param\'etrage de sorte que cet \'el\'ement soit param\'etr\'e par le caract\`ere trivial. On prend maintenant $(V_{+},q_{+})$, $\varphi_{+}$, $(V_{-},q_{-})$, $\varphi_{-}$ quelconques. Le m\^eme raisonnement que dans le cas impair prouve que $\gamma^G(\varphi_{+},\varphi_{-})=1$.

Revenons \`a la situation g\'en\'erale de 4.10, en supposant $G'$ d\'eploy\'e. On a calcul\'e toutes les constantes et l'\'egalit\'e de ce paragraphe se r\'eduit \`a
$$(1) \qquad m(\varphi,s;\varphi',s')=\boldsymbol{\epsilon}(s)\boldsymbol{\epsilon}'(s')(E(\varphi,\varphi')+1)/2.$$
Supposons d'abord $(V_{+},q_{+})=(V,q)$ et $(V'_{+},q'_{+})=(V',q')$. Alors $s=0$, $s'=0$ et $m(\varphi,0;\varphi',0)$ est le nombre de couples $(\sigma,\sigma')\in \Pi^G(\varphi)\times \Pi^{G'}(\varphi')$ tels que $m(\sigma,\sigma')=1$. Si $E(\varphi,\varphi')=-1$, ce nombre vaut $0$. Supposons $E(\varphi,\varphi')=1$.  Le nombre vaut $1$ c'est-\`a-dire qu'il y a un unique couple $(\sigma,\sigma')$ comme ci-dessus. Soit $(\epsilon,\epsilon')$ le couple de caract\`eres qui le param\`etre. Revenons \`a des donn\'ees endoscopiques quelconques. On a l'\'egalit\'e $m(\varphi,s;\varphi',s')=\epsilon(s)\epsilon'(s')$ et l'\'egalit\'e (1) ci-dessus est v\'erifi\'ee pour tous $s$ et $s'$. Cela entra\^{\i}ne $\epsilon=\boldsymbol{\epsilon}$ et $\epsilon'=\boldsymbol{\epsilon}'$.

\bigskip

\subsection{Le cas $G'$ non d\'eploy\'e}

Consid\'erons seulement un espace quadratique $(V',q')$ et $\varphi'\in \Phi_{temp}(G')$. Supposons $G'$ non d\'eploy\'e et $\Pi^{G'}(\varphi')$ non vide. Quitte \`a remplacer $q'$ par un multiple, on peut supposer que $V'$ est somme orthogonale 
de plans hyperboliques et d'un espace de dimension $3$ poss\'edant un \'el\'ement $v$ tel que $q(v,v)=-2\nu_{0}$. Fixons un tel \'el\'ement, notons $(V,q)$ son orthogonal. Alors $(V,q)$ et $(V',q')$ v\'erifient les hypoth\`eses de 2.1 et $\delta(q)\not=1$. Montrons que

(1) il existe $\varphi\in \Phi_{temp}(G)$   tel que l'application $(\sigma,\sigma')\mapsto m(\sigma,\sigma')$ soit non identiquement nulle sur $\Pi^G(\varphi)\times \Pi^{G'}(\varphi')$.

Fixons $\sigma'\in \Pi^{G'}(\varphi')$, munissons son espace $E_{\sigma'}$ d'un produit hermitien invariant d\'efini positif. Pour $e'_{1},e'_{2}\in E_{\sigma'}$, consid\'erons la fonction sur $G(F)$: $g\mapsto f'(g)=(e_{1}',\sigma'(g)e_{2}')$. Pour $e'_{1}$ et $e'_{2}$ convenables, elle est non nulle. D'apr\`es [W2] lemme 4.9, c'est une fonction de Schwartz-Harish-Chandra sur $G(F)$. La formule de Plancherel entra\^{\i}ne qu'il y a au moins une repr\'esentation irr\'eductible et temp\'er\'ee $\sigma$ de $G(F)$ telle que $\sigma^{\vee}(f')$ soit non nulle.  Munissons  l'espace $E_{\sigma}$  d'un produit hermitien invariant d\'efini positif. On peut trouver $e_{1},e_{2}\in E_{\sigma}$ tels que
$$\int_{G(F)}(\sigma(g)e_{1},e_{2})(e'_{1},\sigma'(g)e'_{2})dg\not=0.$$
D'apr\`es [W2] proposition 5.7, cela entra\^{\i}ne $m(\sigma,\sigma')\not=0$. Il suffit de choisir pour $\varphi$ l'\'el\'ement de $\Phi_{temp}(G)$ tel que $\sigma$ appartienne \`a $\Pi^G(\varphi)$. $\square$

Choisissons $\varphi$ v\'erifiant (1), appliquons les constructions de 4.10 \`a $(V_{+},q_{+})=(\underline{V},\underline{q})$, $\varphi_{+}=\varphi$.  Puisque $\delta(q)\not=1$, on a $\gamma^G(\varphi,0)=1$ d'apr\`es la derni\`ere remarque de 4.6. En utilisant  ce que l'on a d\'ej\`a d\'emontr\'e, l'\'egalit\'e de 4.10 se r\'eduit \`a
$$(2) \qquad \gamma^{G'}(\varphi'_{+},\varphi'_{-})m(\varphi,0;\varphi',s')=\boldsymbol{\epsilon}'(s')(E(\varphi,\varphi')-1)/2.$$
Consid\'erons d'abord le cas o\`u $(V'_{+},q'_{+})=(\underline{V}',\underline{q}')$, $\varphi'_{+}=\varphi'$.  Alors $s'=0$ et $m(\varphi,0,\varphi',0)$ est le nombre de couples $(\sigma,\sigma')\in \Pi^G(\varphi)\times \Pi^{G'}(\varphi')$ tels que $m(\sigma,\sigma')=1$. D'apr\`es le choix de $\varphi$, c'est un entier strictement positif. L'\'egalit\'e (2) entra\^{\i}ne que ce nombre est \'egal \`a $1$ et que $E(\varphi,\varphi')=-1$. Soit $(\sigma,\sigma')$ l'unique couple tel que $m(\sigma,\sigma')=1$. Quitte \`a changer le param\'etrage de $\Pi^{G'}(\varphi')$, on peut supposer que $\sigma'$ est param\'etr\'e par $\boldsymbol{\epsilon}'$. Revenons \`a un groupe endoscopique g\'en\'eral de $G'$. Alors $m(\varphi,0;\varphi',s')$ est \'egal \`a $\boldsymbol{\epsilon}'(s')$. L'\'egalit\'e (2) entra\^{\i}ne $\gamma^{G'}(\varphi'_{+},\varphi'_{-})=-1$. Cela ach\`eve de prouver le (i) du lemme 4.8. On a utilis\'e des donn\'ees auxiliaires $(V,q)$ et $\varphi$. Dans la suite, on les oublie, mais on suppose le lemme 4.8(i) v\'erifi\'e.

Consid\'erons maintenant un espace quadratique $(V,q)$ et $\varphi\in \Phi_{temp}(G)$. Supposons  que $(V,q) $ ne v\'erifie pas l'hypoth\`ese (QD) de 1.7  et que  $\Pi^{G}(\varphi)$ est non vide. Si $\delta(q)\not=1$, on d\'eduit comme en 4.6 du param\'etrage   de $\Pi^{\underline{G}}(\varphi)$ fix\'e dans le paragraphe pr\'ec\'edent un param\'etrage de $\Pi^G(\varphi)$ qui satisfait le (ii) du lemme 4.8. Supposons que $\delta(q)=1$, donc que $G$ n'est pas quasi-d\'eploy\'e. Dans ce cas, on peut d\'ecomposer $(V,q)$ en somme orthogonale d'un espace $(V',q')$ et d'une droite $(D,q_{D})$ poss\'edant un \'el\'ement $v$ tel que $q_{D}(v,v)=2\nu_{0}$. Les espaces $(V,q)$ et $(V',q')$ satisfont les hypoth\`eses de 2.1 et $G'$ n'est pas d\'eploy\'e. On d\'emontre l'analogue de (1): il existe $\varphi'\in \Phi_{temp}(G')$ tel qu'il existe $(\sigma,\sigma')\in \Pi^G(\varphi)\times \Pi^{G'}(\varphi')$ de sorte que $m(\sigma,\sigma')=1$. On fixe un tel $\varphi'$ et on applique les constructions de 4.10 \`a $(V'_{+},q'_{+})=(\underline{V}',\underline{q}')$ et $\varphi'_{+}=\varphi'$. D'apr\`es le calcul ci-dessus de $\gamma^{G'}(\varphi',0)$, l'\'egalit\'e de 4.10 devient
$$-\gamma^G(\varphi_{+},\varphi_{-})m(\varphi,s;\varphi',0)=\boldsymbol{\epsilon}(s)(E(\varphi,\varphi')-1)/2.$$
Le m\^eme raisonnement que ci-dessus montre que, quitte \`a modifier le param\'etrage de $\Pi^G(\varphi)$, le (ii) du lemme 4.8 est v\'erifi\'e.

On suppose d\'esormais le lemme 4.8 v\'erifi\'e et on revient \`a la situation g\'en\'erale de 4.10, en supposant $G'$ non d\'eploy\'e. L'\'egalit\'e de ce paragraphe se r\'eduit \`a
 $$-m(\varphi,s;\varphi',s')=\boldsymbol{\epsilon}(s)\boldsymbol{\epsilon}'(s')(E(\varphi,\varphi')-1)/2.$$
 C'est la m\^eme \'egalit\'e que 4.11(1), \`a ceci pr\`es que $E(\varphi,\varphi')$ est chang\'e en $-E(\varphi,\varphi')$. On ach\`eve la preuve du th\'eor\`eme comme dans ce paragraphe.

 \bigskip
 
 {\bf Bibliographie}

[AGRS] A. Aizenbud, D. Gourevitch, S. Rallis, G. Schiffmann: {\it Multiplicity one theorems}, pr\'epublication 2007

[A1] J. Arthur: {\it An introduction to the trace formula}, Clay Math. Proc. 4 (2005), p.1-253

[A2] -----------: {\it A note on $L$-packets}, Pure and A. Math. Quart. 2 (2006), p.199-217

[C] L. Clozel: {\it Characters of non-connected, reductive $p$-adic groups}, Can. J. Math. 39 (1987), p.149-167

[GGP] W.T. Gan, B. Gross, D. Prasad: {\it Symplectic local root numbers, central critical $L$-values, and restriction problems in the representation theory of classical groups}, pr\'epublication 2009

[GP] B. Gross, D. Prasad: {\it On the decomposition of a representation of $SO_{n}$ when restricted to $SO_{n-1}$}, Can. J. Math. 44 (1992), p.974-1002

[GP2] -------------------------:  {\it On irreducible representations of $SO_{2n+1}\times SO_{2m}$}, Can. J. Math. 46 (1994), p.930-950

[HT] M. Harris, R. Taylor: {\it On the geometry and cohomology of some simple Shimura varieties}, Ann. of Math. Studies 151, PUP 2002

[H] G. Henniart: {\it Une preuve simple des conjectures de Langlands pour $GL(n)$ sur un corps $p$-adique}, Invent. Math. 139 (2000), p. 439-456

[KS] R. Kottwitz, D. Shelstad: {\it Foundations of twisted endoscopy}, Ast\'erisque 255 (1999)

[LS] R.P. Langlands, D. Shelstad: {\it On the definition of transfer factors}, Math. Annalen 278 (1987), p.219-271

[M] C. Moeglin: {\it Repr\'esentations quadratiques unipotentes des groupes classiques $p$-adiques}, Duke Math. 84 (1996), p.267-332

[S] J.-P. Serre: {\it Corps locaux}, Hermann 1968

[V] D. Vogan: {\it The local Langlands conjecture}, Contemp. Math. 145 (1993), p.305-379

[W1] J.-L. Waldspurger: {\it Une formule int\'egrale reli\'ee \`a la conjecture locale de Gross-Prasad}, pr\'epublication 2009

[W2] ------------------------: {\it Une formule int\'egrale reli\'ee \`a la conjecture locale de Gross-Prasad, $2^{\grave{e}me}$ partie: extension aux repr\'esentations temp\'er\'ees}, pr\'epublication 2009

[W3] ------------------------: {\it Calcul d'une valeur d'un facteur $\epsilon$ par une formule int\'egrale}, pr\'epublication 2009

[W4] ------------------------: {\it Les facteurs de transfert pour les groupes classiques: un formulaire}, pr\'epublication 2009

[W5] ------------------------: {\it Une variante d'un r\'esultat de Aizenbud, Gourevitch, Rallis et Schiffmann}, pr\'epublication 2009

\end{document}